\documentclass[a4paper, 12pt, leqno]{article}
\usepackage{amsmath}
\usepackage{amscd}
\usepackage{amssymb}
\usepackage{amsthm}
\newtheorem{lemma}{{\sc Lemma}}[section]
\newtheorem{corollary}[lemma]{{\sc Corollary}}
\newtheorem{proposition}[lemma]{{\sc Proposition}}
\newtheorem{theorem}[lemma]{{\sc Theorem}}
\newtheorem{remark}[lemma]{{\sc Remark}}
\newtheorem{conjecture}[lemma]{{\sc Conjecture}}

\numberwithin{equation}{section}

\def\Gg{{\mathfrak{g}}}
\def\Gh{{\mathfrak{h}}}

\def\Gz{{\mathfrak{z}}}
\def\Gsl{{\mathfrak{sl}}}
\def\BC{{\mathbb{C}}}
\def\BF{{\mathbb{F}}}
\def\BK{{\mathbb{K}}}
\def\BP{{\mathbb{P}}}
\def\BQ{{\mathbb{Q}}}
\def\BZ{{\mathbb{Z}}}

\def\CB{{\mathcal B}}

\def\DD{{\mathcal D}}
\def\CO{{\mathcal O}}

\def\CK{{\mathcal K}}
\def\CL{{\mathcal L}}
\def\CM{{\mathcal M}}
\def\CR{{\mathcal R}}
\def\CX{{\mathcal X}}

\def\ad{\mathop{\rm ad}\nolimits}

\def\ch{\mathop{\rm ch}\nolimits}
\def\Coker{\mathop{\rm Coker}\nolimits}
\def\deru{\partial}
\def\e{\mathop{\rm{\bf e}}\nolimits}
\def\End{\mathop{\rm{End}}\nolimits}

\def\fin{{\mathop{\rm{fin}}\nolimits}}
\def\Hom{\mathop{\rm Hom}\nolimits}
\def\id{\mathop{\rm id}\nolimits}
\def\Id{\mathop{\rm Id}\nolimits}
\def\Image{\mathop{\rm Im}\nolimits}
\def\Ker{\mathop{\rm Ker\hskip.5pt}\nolimits}
\def\Lie{\mathop{\rm Lie}\nolimits}
\def\Mod{\mathop{\rm Mod}\nolimits}
\def\MMod{\mathop{{\CM}od}\nolimits}
\def\op{{\mathop{\rm op}\nolimits}}
\def\Proj{\mathop{\rm Proj}\nolimits}
\def\Spec{{\rm{Spec}}}

\def\Tor{{\rm{Tor}}}
\def\oE{{\overline{E}}}
\def\tD{{\tilde{D}}}
\def\tU{{\tilde{U}}}

\title{The Beilinson-Bernstein correspondence for quantized enveloping algebras}
\author{Toshiyuki TANISAKI
\cr\cr
{\footnotesize Department of Mathematics, 
Osaka City University,} \cr
{\footnotesize 3-3-138, Sugimoto, Sumiyoshi-ku,}\cr
{\footnotesize Osaka, 558-8585 Japan}\cr
{\footnotesize \texttt{tanisaki@sci.osaka-cu.ac.jp}}
}

\begin{document}
\date{September 22, 2003}
\maketitle
\renewcommand{\thefootnote}{}
\begin{abstract}
Theory of the quantized flag manifold as a quasi-scheme (non-commutative scheme) has been developed by Lunts-Rosenberg \cite{LR}.
They have formulated an analogue of the Beilinson-Bernstein correspondence using the $q$-differential operators introduced in their earlier paper \cite{LRD}.
In this paper we shall establish its modified version using a class of $q$-differential operators, which is (possibly) smaller than the one in \cite{LRD}.
\end{abstract}
\setcounter{section}{-1}
\section{Introduction}
\label{sec:Intoro}
Let $G$ be a connected, simply-connected semisimple algebraic group over the complex number field $\BC$, and let $B$ and $B^-$ be Borel subgroups of $G$ such that $H=B\cap B^-$ is a maximal torus of $G$.
Denote the Weyl group by $W$ and the character group of $H$ by $\Lambda$.
We choose a system of positive roots $\Delta^+\subset\Lambda$ as the weights of $\Lie([B,B])$.
Let $\Lambda^+$ be the set of dominant integral weights.
We have $\Lambda=\sum_{i=1}^\ell\BZ\varpi_i,\,\,\Lambda^+=\sum_{i=1}^\ell\BZ_{\geqq0}\varpi_i$, where $\{\varpi_1,\dots, \varpi_\ell\}$ denotes the set of fundamental weights.
For $\lambda\in\Lambda^+$ we denote the irreducible (left) $G$-module with highest weight $\lambda$ by $V^1(\lambda)$.

The algebraic variety $\CB=B^-\backslash G$ is called the flag manifold for $G$.
It has an affine open covering 
\[
\CB=\bigcup_{w\in W}U_w,\qquad
U_w=B^-\backslash B^-Bw=\Spec(R^1_w).
\]
We have a closed embedding of $\CB$ into 
$\BP(V^1(\varpi_1)^*)\times\cdots\times\BP(V^1(\varpi_\ell)^*)$ given by 
$B^-g\mapsto([v_1g],\dots,[v_\ell g])$,
where $v_i$ is a non-zero element of the (right) $G$-module $V^1(\varpi_i)^*=\Hom_\BC(V^1(\varpi_i),\BC)$ satisfying $v_ih=\varpi_i(h)v_i$ for $h\in H$.
Hence we have $\CB=\Proj_{\BZ^\ell}(A^1)$ for a $\BZ^\ell$-graded ring $A^1$.
The graded ring $A^1$ is described as follows.
Let $\BC[G]$ denote the coordinate algebra of $G$.
One has a natural $G$-bimodule structure on $\BC[G]$.
Then we have the identification
\[
A^1=\{\varphi\in\BC[G]\mid \varphi g=\varphi\,\,(g\in[B^-,B^-])\},
\]
and the grading $A^1=\bigoplus_{\lambda\in\Lambda^+}A^1(\lambda)$ 
by $\Lambda(\simeq\BZ^\ell)$ is given by
\[
A^1(\lambda)
=\{\varphi\in A^1\mid \varphi h=\lambda(h)\varphi\,\,(h\in H)\} \qquad(\lambda\in \Lambda^+).
\]

In this paper we shall be concerned with the $q$-analogue $\CB_q$ of the flag manifold $\CB$ introduced by Lunts-Rosenberg \cite{LR}.
Let $U=U_q(\Gg)$ denote the corresponding simply-connected quantized enveloping algebra.
It is a Hopf algebra over the field $\BF=\BQ(q^{1/\ell_0})$, where $q$ is transcendental and $\ell_0$ is an appropriate positive integer. 
See \S\ref{sec:QG} below for the precise definition.
We note that the Cartan part of $U$ is identified with the group algebra of the weight lattice $\Lambda$.
It is well-known that a $q$-analogue $\BC_q[G]$ of $\BC[G]$ is defined as a Hopf algebra dual to $U$.
Using $\BC_q[G]$ one can naturally define $q$-analogues $A$ and $R_w$ of $A^1$ and $R^1_w$ respectively.
$A$ is a $\Lambda$-graded $\BF$-algebra and  $R_w$ is an $\BF$-algebra; however, they are non-commutative.
Hence in order to give meanings to
\begin{equation}
 \CB_q=\Proj_\Lambda(A)=\bigcup_{w\in W}U_{w,q},\qquad
 U_{w,q}=\Spec(R_w)
 \end{equation}
we need the notion of non-commutative algebraic varieties.

A starting point of the non-commutative algebraic geometry is the general fact that 
a (commutative) scheme $X$ is uniquely determined from the category $\MMod(\CO_X)$ of quasi-coherent sheaves (Rosenberg \cite{R}).
It was Manin who proposed to consider generalization of the category $\MMod(\CO_X)$ in the non-commutative setting.
There already exist several works along the line of this Manin's idea.
Theory of non-commutative projective schemes for non-commutative graded rings are developed by Artin-Zhang \cite{AZ}, Manin \cite{Manin}, Verevkin \cite{V}, Rosenberg \cite{R:book}.
Theory of general non-commutative schemes equipped with non-commutative affine open covering is also given by Rosenberg \cite{R}.
The quantized flag manifold $\CB_q$ is a non-commutative projective scheme as well as a general non-commutative scheme (quasi-scheme) in the sense of Rosenberg.

Let us give a description of $\CB_q$ as a non-commutative projective scheme.
For a ring $R$ we denote the category of left $R$-modules by $\Mod(R)$.
For a $\Lambda$-graded ring $R=\bigoplus_{\lambda\in\Lambda}R(\lambda)$ a left $R$-module $M$ equipped with the decomposition $M=\bigoplus_{\lambda\in\Lambda}M(\lambda)$ of $R(0)$-submodules is called a $\Lambda$-graded left $R$-module if $R(\xi)M(\lambda)\subset M(\lambda+\xi)$ for any $\lambda, \xi\in\Lambda$.
We denote the category of $\Lambda$-graded left $R$-modules by $\Mod_\Lambda(R)$.
Let $\Tor_\Lambda(A)$ be the full subcategory of $\Mod_\Lambda(A)$ consisting of $M\in\Mod_\Lambda(A)$ such that for any $m\in M$ there exists some $\mu\in\Lambda^+$ satisfying $A(\xi)m=\{0\}$ for any 
$\xi\in\mu+\Lambda^+$.
Then we can define an abelian category $\MMod(\CO_{\CB_q})$ of ``quasi-coherent sheaves on $\CB_q$" as the quotient
\[
\MMod(\CO_{\CB_q})
=\Mod_\Lambda(A)/\Tor_\Lambda(A).
\]
Moreover, ``the global section functor"
\[
\Gamma:\MMod(\CO_{\CB_q})\to\Mod(\BF)
\]
is defined as follows.
The natural functor 
\[
\omega^*:\Mod_\Lambda(A)\to\MMod(\CO_{\CB_q})
=\Mod_\Lambda(A)/\Tor_\Lambda(A)
\]
admits a right adjoint functor 
\[
\omega_*:\MMod(\CO_{\CB_q})\to\Mod_\Lambda(A),
\]
and $\Gamma$ is defined by
$\Gamma(M)=(\omega_*M)(0)$.
We can also define the higher cohomology groups $H^i(M)$ for an object $M$ of $\MMod(\CO_{\CB_q})$ by $H^i(M)=(R^i\Gamma)(M)$ using the right derived functors.

Now let us consider ``$\DD$-modules on $\CB_q$".
The ring $\tD$ of $q$-differential operators acting on the graded algebra $A$ is defined by Lunts-Rosenberg \cite{LRD, LR} as follows.
For $\varphi\in A$ let $\ell_\varphi, r_\varphi\in\End_\BF(A)$ denotes the left multiplication and the right multiplication respectively; i.e. $\ell_\varphi(\psi)=\varphi\psi, r_\varphi(\psi)=\psi\varphi$ for $\psi\in A$.
Define an increasing sequence 
\[
\{0\}=F^{-1}\tD\subset F^0\tD\subset F^1\tD\subset\cdots\subset\End_\BF(A)
\]
of $\BF$-subspaces of $\End_\BF(A)$ inductively by
\[
F^p\tD=\sum_{\varphi_1,\varphi_2\in A, \lambda,\mu\in\Lambda}\ell_{\varphi_1}(F^p\tD)_{\lambda,\mu}\ell_{\varphi_2},
\]
where $(F^p\tD)_{\lambda,\mu}$ consists of $d\in\End_\BF(A)$ satisfying
\begin{itemize}
\item[(a)]
$d(A(\xi))\subset A(\xi+\lambda)$ for any $\xi\in\Lambda$,
\item[(b)]
$d\ell_{\varphi}-q^{(\mu,\xi)}\ell_\varphi d\in F^{p-1}\tD$ for any $\xi\in\Lambda$ and $\varphi\in A(\xi)$.
\end{itemize}
Here, $(\,,\,):\Lambda\times\Lambda\to\BQ$ is the restriction of a standard symmetric bilinear form on $\BQ\otimes_{\BZ}\Lambda$.
Set
\[
\tD=\bigcup_{p}F^p\tD\subset\End_\BF(A).
\]
Then $\tD$ is a subalgebra of $\End_\BF(A)$, and its $\Lambda$-grading $\tD=\bigoplus_{\lambda\in\Lambda}\tD(\lambda)$ is given by 
\[
\tD(\lambda)=\{d\in \tD\mid d(A(\xi))\subset A(\xi+\lambda)\}.
\]
It seems to be a hard task to determine how large $\tD$ is.
Anyway $\tD$ contains $\ell_\varphi, r_\varphi$ for $\varphi\in A$, operators $\partial_u\,\,(u\in U)$ given by the natural action of $U$ on $A$, and the degree operators $\sigma_\lambda\,\,(\lambda\in\Lambda)$ given by $\sigma_\lambda|A(\xi)=q^{(\lambda,\xi)}\id$.
We denote by $D$ the subalgebra of $\tD$ generated by the operators $\ell_\varphi, r_\varphi\,\,(\varphi\in A)$, $\partial_u\,\,(u\in U)$,  $\sigma_\lambda\,\,(\lambda\in\Lambda)$.
It is shown using the universal $R$-matrix that $D$ is generated by $\ell_\varphi\,\,(\varphi\in A)$, $\partial_u\,\,(u\in U)$,  $\sigma_\lambda\,\,(\lambda\in\Lambda)$.

Let $\lambda\in\Lambda$.
We define the category $\MMod(\tilde{\DD}_{\CB_q,\lambda})$ of ``modules over the sheaf of rings of twisted differential operators $\tilde{\DD}_{\CB_q,\lambda}$" by
\begin{equation}
\label{eq:MD}
\MMod(\tilde{\DD}_{\CB_q,\lambda})
=\Mod_{\Lambda,\lambda}(\tD)/\Tor_{\Lambda,\lambda}(\tD),
\end{equation}
where $\Mod_{\Lambda,\lambda}(\tD)$ denotes the category of $\Lambda$-graded left $\tD$-modules $M$ satisfying $\sigma_\mu|M(\xi)=q^{(\mu,\lambda+\xi)}\id$ for any $\mu, \xi\in\Lambda$, and $\Tor_{\Lambda,\lambda}(\tD)$ is its full subcategory consisting of objects of $\Mod_{\Lambda,\lambda}(\tD)$ which belong to $\Tor_\Lambda(A)$ as $\Lambda$-graded $A$-modules.
Here, we identify $A$ with a graded subring of $\tD$ by $A\ni\varphi\mapsto\ell_\varphi\in \tD$.

Define $\tD_\lambda\in\Mod_{\Lambda,\lambda}(\tD)$ by
\[
\tD_\lambda=\tD/\sum_{\mu\in\Lambda}\tD(\sigma_\mu-q^{(\mu,\lambda)}).
\]
Since $\sigma_\mu$ belongs to the center of $\tD(0)$ we have an $\BF$-algebra structure on $\tD_\lambda(0)$.
Then the global section functor $\Gamma:\MMod(\CO_{\CB_q})\to\Mod(\BF)$ induces a left exact functor
\[
\tilde{\Gamma}_\lambda:\MMod(\tilde{\DD}_{\CB_q,\lambda})\to\Mod(\tD_\lambda(0)).
\]
Denote the Verma module for $U$ with highest weight $\lambda$ by $T(\lambda)$ and its annihilator in $U$ by $J_\lambda$.
By Joseph \cite{J:an} the ideal $J_\lambda$ is generated by its intersection with the center of $U$.
We have canonical $\BF$-algebra homomorphisms
\[
U/J_\lambda\to \tD_\lambda(0)\to\tilde{\Gamma}_\lambda(\omega^*\tD_\lambda).
\]
Let $\rho\in\Lambda$ denote the half sum of the positive roots.
\begin{theorem}[Lunts-Rosenberg]
\label{thm:LR}
If $\lambda+\rho\in\Lambda^+$, then the functor $\tilde{\Gamma}_\lambda$ is exact.
\end{theorem}
\begin{conjecture}[Lunts-Rosenberg]
\label{conj1:LR}
If $\lambda\in\Lambda^+$, then $\tilde{\Gamma}_\lambda(M)=0$ implies $M=0$.
\end{conjecture}
\begin{conjecture}[Lunts-Rosenberg]
\label{conj2:LR}
For any $\lambda\in\Lambda$ we have $U/J_\lambda\simeq \tD_\lambda(0)\simeq\tilde{\Gamma}_\lambda(\omega^*\tD_\lambda)$.
\end{conjecture}

By a standard argument Theorem \ref{thm:LR}, Conjecture \ref{conj1:LR} and Conjecture \ref{conj2:LR} for $\lambda\in\Lambda^+$ imply the following analogue of the Beilinson-Bernstein correspondence (Beilinson-Bernstein \cite{BB}).
\begin{conjecture}[Lunts-Rosenberg]
\label{conj3:LR}
For $\lambda\in\Lambda^+$ $\tilde{\Gamma}_\lambda$ induces the equivalence of categories:
\[
\MMod(\tilde{\DD}_{\CB_q,\lambda})\simeq\Mod(U/J_\lambda).
\]
\end{conjecture}

We can define $\MMod(\DD_{\CB_q,\lambda})$, $D_\lambda$, $\Gamma_\lambda:\MMod(\DD_{\CB_q,\lambda})\to\Mod(D_\lambda(0))$ and 
$U/J_\lambda\to D_\lambda(0)\to\Gamma_\lambda(\omega^*D_\lambda)$ using $D$ instead of $\tD$.
Our main result is the following.
\begin{theorem}
\label{thm1}
Conjecture \ref{conj1:LR} is true.
\end{theorem}
\begin{theorem}
\label{thm2}
Theorem \ref{thm:LR}, Conjecture \ref{conj1:LR} and Conjecture \ref{conj2:LR} for $\lambda+\rho\in\Lambda^+$ are all true for $D$; that is,
\begin{itemize}
\item[{\rm(i)}]
If $\lambda+\rho\in\Lambda^+$, then the functor $\Gamma_\lambda$ is exact.
\item[{\rm(ii)}]
If $\lambda\in\Lambda^+$, then $\Gamma_\lambda(M)=0$ implies $M=0$.
\item[{\rm(iii)}]
For any $\lambda\in\Lambda$ we have 
$U/J_\lambda\simeq D_\lambda(0)$.
\item[{\rm(iv)}]
If $\lambda+\rho\in\Lambda^+$, then we have 
$D_\lambda(0)\simeq\Gamma_\lambda(\omega^*D_\lambda)$.
\end{itemize}
\end{theorem}

We can deduce the following from Theorem \ref{thm2}.
\begin{corollary}
\label{cor3}
For $\lambda\in\Lambda^+$ $\Gamma_\lambda$ induces the equivalence of categories:
\[
\MMod(\DD_{\CB_q,\lambda})\simeq\Mod(U/J_\lambda).
\]
\end{corollary}

In Lunts-Rosenberg \cite{LR} it is noted that a $q$-analogue of the formula
\begin{equation}
\label{eq:obvious}
\Gamma(\CB,(\CO_\CB\otimes_\BC V^1(\mu))\otimes_{\CO_\CB}M)
\simeq
V^1(\mu)\otimes_\BC \Gamma(\CB,M),
\end{equation}
implies Theorem \ref{thm1}.
We can show it using basic properties of universal $R$-matrices (Proposition \ref{prop:filt-E1} below),  from which we obtain Theorem \ref{thm1}.
The proofs of Theorem \ref{thm2} (i) and (ii) are similar to those for Theorem \ref{thm:LR} and Theorem \ref{thm1} respectively.
Our proof of Theorem \ref{thm2} (iii) and (iv) is similar to that of the corresponding fact for Lie algebras given by Borho-Brylinski \cite{BoB}.
The proof by Borho-Brylinski uses the structure of the annihilators of Verma modules and a result of N. Conze-Berline and M. Duflo.
In the quantum setting both the structure theorem of the annihilators of Verma modules and an analogue of the theorem of N. Conze-Berline and M. Duflo are already obtained by Joseph \cite{J:an}, \cite{J:book}.
Theorem \ref{thm2} (iii) follows from the result about the annihilators of Verma modules easily; 
however, unlike the Lie algebra case Joseph's theorem giving an analogue of a result by N. Conze-Berline and M. Duflo does  not immediately imply Theorem \ref{thm2} (iv) since $U_q(\Gg)$ is not locally finite with respect to the adjoint action.
We overcome the difficulty by the arguments used in the proof of Theorem \ref{thm2} (i), where the assumption $\lambda+\rho\in\Lambda^+$ is necessary.

Let us give a comment in order to justify the usage of $D$ instead of $\tD$.
Let ${\tD^1}$ and ${D^1}$ denote the subalgebras of $\End_\BC(A^1)$ corresponding to $\tD$ and $D$ in the ordinary enveloping algebra situation.
The algebra ${D^1}$ is in fact a proper subalgebra of $\tD^1$ by Bernstein-Gelfand--Gelfand \cite[Example 2]{BGG}; however, the corresponding categories  $\MMod(\tilde{\DD}_{\CB,\lambda})$ and $\MMod(\DD_{\CB,\lambda})$ defined similarly to \eqref{eq:MD} are equivalent since the corresponding rings of differential operators are isomorphic locally on $\CB$.

We note that Theorem \ref{thm2} and Corollary \ref{cor3} for $\Gg=\Gsl(2)$ is due to Hodges \cite{Hodges}.
We also note that a different approach to the Beilinson-Bernstein correspondence for the quantized enveloping algebras is given in Joseph \cite{J:emb}.

In this paper we shall use the following notation for a Hopf algebra $H$ over a field $\BK$.
The multiplication, the unit, the comultiplication, the counit, and the antipode of $H$ are denoted by
\begin{align}
&m_H:H\otimes_\BK H\to H,\\
&\eta_H:\BK\to H,\\
&\Delta_H:H\to H\otimes_\BK H,\\
&\epsilon_H:H\to\BK,\\
&S_H:H\to H
\end{align}
respectively.
The subscript $H$ will often be omitted.
For $n\in\BZ_{>0}$ we denote by
\[
\Delta_n:H\to H^{\otimes n+1}
\]
the algebra homomorphism given by 
\[
\Delta_1=\Delta,\qquad
\Delta_n=(\Delta\otimes \id_{H^{\otimes n-1}})\circ\Delta_{n-1},
\]
and write
\[
\Delta_n(h)=\sum_{(h)_n}h_{(0)}\otimes\cdots\otimes h_{(n)}.
\]

\section{Quantum groups}
\label{sec:QG}
\subsection{Quantized enveloping algebras}
Let $\Gg$ be a finite-dimensional semisimple Lie algebra over $\BC$ and let $\Gh$ be its Cartan subalgebra.
We denote by $\Delta\subset\Gh^*, \Lambda\subset\Gh^*, W\subset GL(\Gh^*)$ the set of roots, the weight lattice and the Weyl group respectively.
We fix a set of simple roots $\{\alpha_i\}_{i\in I}$.
Let $\Delta^+\subset\Gh^*, \Lambda^+\subset\Gh^*, \{\varpi_i\}_{i\in I}\subset\Gh^*, \{s_i\}_{i\in I}\subset W$ denote the corresponding set of positive roots, dominant weights, fundamental weights, and simple reflections respectively.
Set 
\[
Q^+
=\sum_{\alpha\in\Delta^+}\BZ_{\geqq0}\alpha
=\bigoplus_{i\in I}\BZ_{\geqq0}\alpha_i
\subset \Gh^*.
\]
Let $\rho\in\Gh^*$ be the half sum of positive roots.
We denote the longest element of $W$ by $w_0$.
We fix a $W$-invariant symmetric bilinear form 
\begin{equation}
(\,,\,):\Gh^*\times\Gh^*\to\BC
\end{equation}
satisfying $(\alpha_i,\alpha_i)\in\BQ_{>0}$ for any $i\in I$.
For $i\in I$ we set
\begin{equation}
\alpha_i^\vee=2\alpha_i/(\alpha_i,\alpha_i)\in\Gh^*.
\end{equation}

Take a positive integer $\ell_0$ satisfying
\begin{equation}
{\ell_0}(\alpha_i,\alpha_i)\subset2\BZ\quad (i\in I),\qquad
{\ell_0}(\Lambda,\Lambda)\subset\BZ,
\end{equation}
and let $\BF=\BQ(q^{1/\ell_0})$ be the rational function field over $\BQ$ with variable $q^{1/\ell_0}$.
In this paper $\otimes$ stands for $\otimes_{\BF}$.

For $n\in\BZ_{\geqq0}$ we set
\[
[n]_t=\frac{t^n-t^{-n}}{t-t^{-1}}\in\BZ[t,t^{-1}],
\qquad
[n]_t!=[n]_t[n-1]_t\cdots[2]_t[1]_t\in\BZ[t,t^{-1}].
\]

The simply-connected quantized enveloping algebra $U=U_q(\Gg)$ is an associative algebra over  $\BF$ with the identity element $1$ generated by the elements $k_\lambda\,(\lambda\in\Lambda),\,e_i, f_i\,(i\in I)$ satisfying the following defining relations:
\begin{align}
&k_0=1,\quad 
k_\lambda k_\mu=k_{\lambda+\mu}
\qquad(\lambda,\mu\in\Lambda),
\label{eq:def1}\\
&k_\lambda e_ik_\lambda^{-1}=q^{(\lambda,\alpha_i)}e_i,\qquad(\lambda\in\Lambda, i\in I),
\label{eq:def2a}\\
&k_\lambda f_ik_\lambda^{-1}=q^{-(\lambda,\alpha_i)}f_i
\qquad(\lambda\in\Lambda, i\in I),
\label{eq:def2b}\\
&e_if_j-f_je_i=\delta_{ij}\frac{k_i-k_i^{-1}}{q_i-q_i^{-1}}
\qquad(i, j\in I),
\label{eq:def3}\\
&\sum_{n=0}^{1-a_{ij}}(-1)^ne_i^{(1-a_{ij}-n)}e_je_i^{(n)}=0
\qquad(i,j\in I,\,i\ne j),
\label{eq:def4}\\
&\sum_{n=0}^{1-a_{ij}}(-1)^nf_i^{(1-a_{ij}-n)}f_jf_i^{(n)}=0
\qquad(i,j\in I,\,i\ne j),
\label{eq:def5}
\end{align}
where $q_i=q^{(\alpha_i,\alpha_i)/2}, k_i=k_{\alpha_i}, a_{ij}=2(\alpha_i,\alpha_j)/(\alpha_i,\alpha_i)$ for $i, j\in I$, and
\[
e_i^{(n)}=
e_i^n/[n]_{q_i}!,
\qquad
f_i^{(n)}=
f_i^n/[n]_{q_i}!
\]
for $i\in I$ and $n\in\BZ_{\geqq0}$.
Algebra homomorphisms $\Delta:U\to U\otimes U, \epsilon:U\to\BF$ and an algebra anti-automorphism $S:U\to U$ are defined by:
\begin{align}
&\Delta(k_\lambda)=k_\lambda\otimes k_\lambda,\\
&\Delta(e_i)=e_i\otimes 1+k_i\otimes e_i,\quad
\Delta(f_i)=f_i\otimes k_i^{-1}+1\otimes f_i,
\nonumber\\
&\epsilon(k_\lambda)=1,\quad
\epsilon(e_i)=\epsilon(f_i)=0,\\
&S(k_\lambda)=k_\lambda^{-1},\quad
S(e_i)=-k_i^{-1}e_i, \quad S(f_i)=-f_ik_i,
\end{align}
and $U$ is endowed with a Hopf algebra structure with the comultiplication $\Delta$, the counit $\epsilon$ and the antipode $S$.

We define subalgebras $U^0, U^{\geqq0}, U^{\leqq0}, U^{+}, U^{-}$ of $U$ by 
\begin{align}
U^0&=\langle k_\lambda\mid \lambda\in\Lambda\rangle,\\
U^{\geqq0}&=\langle k_\lambda, e_i\mid \lambda\in\Lambda, i\in I\rangle,\\
U^{\leqq0}&=\langle k_\lambda, f_i\mid \lambda\in\Lambda, i\in I\rangle,\\
U^{+}&=\langle e_i\mid i\in I\rangle,\\
U^{-}&=\langle f_i\mid i\in I\rangle.
\end{align}

Note that $U^0, U^{\geqq0}, U^{\leqq0}$ are Hopf subalgebras of $U$, while $U^+$ and $U^-$ are not Hopf subalgebras.

The following result is standard.
\begin{proposition}
\label{prop:Str-of-U}
\begin{itemize}
\item[\rm(i)]
$\{k_\lambda\mid \lambda\in\Lambda\}$ is an $\BF$-basis of $U^0$.
\item[\rm(ii)]
$U^+$ $($resp.\ $U^-$$)$ is isomorphic to the $\BF$-algebra generated by 
$\{e_i\mid i\in I\}$ $($resp.\ $\{f_i\mid i\in I\}$$)$ with defining relation 
\eqref{eq:def4} $($resp.\ \eqref{eq:def5}$)$.
\item[\rm(iii)]
$U^{\geqq0}$ $($resp.\ $U^{\leqq0}$$)$ is isomorphic to the $\BF$-algebra generated by $\{e_i, k_\lambda\mid i\in I, \lambda\in\Lambda\}$ 
$($resp.\ $\{f_i, k_\lambda\mid i\in I, \lambda\in\Lambda\}$$)$ 
with defining relations \eqref{eq:def1}, \eqref{eq:def2a}, \eqref{eq:def4} $($resp.\ \eqref{eq:def1}, \eqref{eq:def2b}, \eqref{eq:def5}$)$.
\item[\rm(iv)]
The linear maps 
\begin{gather*}
U^-\otimes U^0\otimes U^+\to U \gets U^+\otimes U^0\otimes U^-,\\
U^+\otimes U^0\to U^{\geqq0} \gets U^0\otimes U^+,\qquad
U^-\otimes U^0\to U^{\leqq0} \gets U^0\otimes U^-
\end{gather*}
induced by the multiplication are all isomorphisms.
\end{itemize}
\end{proposition}

We define Hopf algebra homomorphisms
\begin{equation}
\pi^+:U^{\geqq0}\to U^0,\qquad
\pi^-:U^{\leqq0}\to U^0
\end{equation}
by $\pi^{\pm}(k_\lambda)=k_\lambda\,\,(\lambda\in\Lambda),\,\,\pi^+(e_i)=\pi^-(f_i)=0\,\,(i\in I)$.

For $\gamma\in Q^+$ we set
\[
U^\pm_{\pm\gamma}
=\{x\in U^\pm\mid
k_\lambda xk_\lambda^{-1}=q^{\pm(\lambda,\gamma)}x\,\,(\lambda\in\Lambda)\}.
\]
We have
\[
U^\pm=\bigoplus_{\gamma\in Q^+}U^\pm_{\pm\gamma}.
\]
\subsection{Representations}
Let $H$ be a Hopf algebra over a field $\BK$.
For left $H$-modules $V_1, V_2$ we endow $V_1\otimes_\BK V_2$ with a left $H$-module structure by
\begin{equation}
h(v_1\otimes v_2)=\Delta(u)(v_1\otimes v_2)
\qquad(h\in H, v_1\in V_1, v_2\in V_2).
\end{equation}

For a left $H$-module $V$ its dual space $V^*=\Hom_\BK(V,\BK)$ is endowed with a structure of a right $H$-module (i.e., a left $H^{\op}$-module, where $H^{\op}$ denotes the algebra opposite to $H$) by
\begin{equation}
\langle v^*h,v\rangle=\langle v^*, hv\rangle\qquad
(v^*\in V^*,\, h\in H,\, v\in V),
\end{equation}
where $\langle\,,\,\rangle:H^*\times H\to\BK$ denotes the canonical pairing.

For $\lambda\in\Lambda$ we define an algebra homomorphism 
\begin{equation}
\chi_\lambda:U^0\to\BF
\end{equation}
by $\chi_\lambda(k_\mu)=q^{(\lambda,\mu)}\,(\mu\in\Lambda)$.
We can extend it to algebra homomorphisms 
\begin{equation}
\chi_\lambda^+:U^{\geqq0}\to\BF,\qquad
\chi_\lambda^-:U^{\leqq0}\to\BF
\end{equation}
by $\chi_\lambda^\pm=\chi_\lambda\circ\pi^\pm$.

For a left (resp.\ right) $U$-module $V$ and $\lambda\in\Lambda$ the subspace 
\begin{gather}
V_\lambda=\{v\in V\mid tv=\chi_\lambda(t)v\,\,\,(t\in U^0\}\\
(\mbox{resp.}\quad 
V_\lambda=\{v\in V\mid vt=\chi_\lambda(t)v\,\,\,(t\in U^0\})
\end{gather}
of $V$ is called the weight space with weight $\lambda$.
Those $\lambda\in\Lambda$ such that $V_\lambda\ne\{0\}$ are called the weights of $V$.
For a left (or right) $U$-module $V$ which is a direct sum of finite-dimensional  weight spaces we define its character by the formal sum
\[
\ch(V)=\sum_{\lambda\in \Lambda}\dim V_\lambda \e^\lambda.
\]
We denote by $\Mod^f(U)$ (resp.\ $\Mod^{f}(U^{\op})$) the abelian category whose objects are finite-dimensional left (resp.\ right) $U$-modules which are direct sums of weight spaces.

For any $\lambda\in\Lambda^+$ there exists a unique irreducible object $V(\lambda)$ of $\Mod^f(U)$ such that $\lambda$ is a weight of $V$ and any weight of $V$ belongs to $\lambda-Q^+$.
Any object of $\Mod^f(U)$ is a direct sum of this type of irreducible objects.
As in the Lie algebra case the character of $V(\lambda)$ is given by Weyl's character formula:
\begin{equation}
\label{eq:Weyl}
\ch(V(\lambda))=
\frac
{\sum_{w\in W}(-1)^{\det(w)}\e^{w(\lambda+\rho)-\rho}}
{\prod_{\alpha\in\Delta^+}(1-\e^{-\alpha})}
\qquad(\lambda\in\Lambda^+)
\end{equation}
(Lusztig \cite{L1}).
By \eqref{eq:Weyl} we obtain the following.
\begin{lemma}
\label{lem:V}
Let $\gamma\in Q^+$.
\begin{itemize}
\item[\rm(i)]
For sufficiently large $\lambda\in\Lambda^+$ the linear maps
\[
U^-_{-\gamma}\ni u\mapsto uv_\lambda\in V(\lambda)_{\lambda-\gamma},
\qquad
U^+_{\gamma}\ni u\mapsto v^*_\lambda u\in V^*(\lambda)_{\lambda-\gamma}
\]
are bijective.
Here, $v_\lambda$ and $v^*_\lambda$ are non-zero elements of $V(\lambda)_\lambda$ and $V^*(\lambda)_\lambda$ respectively.
\item[\rm(ii)]
For sufficiently large $\lambda\in\Lambda^+$ the linear maps
\[
U^+_{\gamma}
\ni u\mapsto uv_{-\lambda}\in 
V(-w_0\lambda)_{-\lambda+\gamma},
\quad
U^-_{-\gamma}
\ni u\mapsto v^*_{-\lambda}u\in 
V^*(-w_0\lambda)_{-\lambda+\gamma}
\]
are bijective.
Here, $v_{-\lambda}$ and $v^*_{-\lambda}$ are non-zero elements of 
$V(-w_0\lambda)_{-\lambda}$ and $V^*(-w_0\lambda)_{-\lambda}$ respectively.
\end{itemize}
\end{lemma}

\begin{remark}
{\rm
In this paper the expression 
``for sufficiently large $\lambda\in\Lambda^+$ ..." 
means that 
``there exists some $\xi\in\Lambda^+$ such that 
for any $\lambda\in\xi+\Lambda^+$ ...".
}
\end{remark}

For any $V\in\Mod^f(U)$ we have $(V^*)_\lambda\simeq(V_\lambda)^*$ and hence $\ch(V)=\ch(V^*)$.
For $\lambda\in\Lambda^+$ we set $V^*(\lambda)=(V(\lambda))^*$.
Any object of $\Mod^f(U^{\op})$ is a direct sum of irreducible submodules isomorphic to $V^*(\lambda)$ for $\lambda\in\Lambda^+$.
We note 
\begin{align}
V(\lambda)_\lambda&=\{v\in V(\lambda)\mid e_iv=0\,\,\,(i\in I)\},\\
V^*(\lambda)_\lambda&=\{v\in V^*(\lambda)\mid vf_i=0\,\,\,(i\in I)\}.
\label{eq:right-hw}
\end{align}

For a left (resp.\ right) $U$-module $V$ which is a direct sum of finite-dimensional weight spaces we define its restricted dual  $V^\bigstar$ by
\begin{equation}
V^\bigstar=\sum_{\lambda\in\Lambda}(V_\lambda)^*\subset V^*.
\end{equation}
It is easily seen that $V^\bigstar$ is a right (resp.\ left) $U$-submodule of $V^*$.
We have
\[
(V^\bigstar)^\bigstar\simeq V, \qquad
\ch(V^\bigstar)=\ch(V).
\]
For $\lambda\in\Lambda$ we define left $U$-modules $T(\lambda)$, $T^*(\lambda)$ and right $U$-modules 
$T_{\rm r}(\lambda)$, $T^*_{\rm r}(\lambda)$ by
\begin{align}
\label{eq:Verma}
&T(\lambda)=U/\sum_{u\in U^{\geqq0}}U(u-\chi^+_\lambda(u)),\qquad
T_{\rm r}(\lambda)
=U/\sum_{u\in U^{\leqq0}}(u-\chi^-_\lambda(u))U,\\
\label{eq:dual Verma}
&T^*(\lambda)=(T_{\rm r}(\lambda))^\bigstar,\qquad
T^*_{\rm r}(\lambda)
=(T(\lambda))^\bigstar.
\end{align}
We have
\[
\ch(T(\lambda))=\ch(T^*(\lambda))=
\ch(T_{\rm r}(\lambda))=\ch(T_{\rm r}^*(\lambda))
=\frac{\e^\lambda}
{\prod_{\alpha\in\Delta^+}(1-\e^{-\alpha})}.
\]

\subsection{Universal $R$-matrices}
\label{subsec:R}
There exists a unique bilinear form
\begin{equation}
\label{eq:Drinfeld-paring}
(\,\,,\,\,):U^{\geqq0}\times U^{\leqq0}\to\BF
\end{equation}
satisfying
\begin{align}
&(x,y_1y_2)=(\Delta(x),y_1\otimes y_2)
&(x\in U^{\geqq0},\,y_1,y_2\in U^{\leqq0}),\\
&(x_1x_2,y)=(x_2\otimes x_1,\Delta(y))
&(x_1, x_2\in U^{\geqq0},\,y\in U^{\leqq0}),\\
&(k_\lambda,k_\mu)=q^{-(\lambda,\mu)}
&(\lambda,\mu\in\Lambda),\\
&(k_\lambda, f_i)=(e_i,k_\lambda)=0
&(\lambda\in\Lambda,\,i\in I),\\
&(e_i,f_j)=\delta_{ij}/(q_i^{-1}-q_i)
&(i,j\in I)
\end{align}
(see Tanisaki \cite{T}).
For any $\beta\in Q^+$ the restriction of \eqref{eq:Drinfeld-paring} to $U^+_\beta\times U^-_{-\beta}$ is non-degenerate, and we denote the corresponding canonical element of $U^+_\beta\otimes U^-_{-\beta}$ by $\Xi_\beta$.
Set
\begin{equation}
\Xi=\sum_{\beta\in Q^+}q^{(\beta,\beta)}
(k_\beta^{-1}\otimes k_\beta)\Xi_\beta.
\end{equation}
It belongs to a completion of $U\otimes U$.

Let $V, V'\in\Mod^f(U)$.
We define $\tau_{V,V'}\in\Hom_\BF(V\otimes V',V'\otimes V)$ and $\kappa_{V,V'}\in GL(V\otimes V')$ by $\tau_{V,V'}(v\otimes v')=v'\otimes v$ and $\kappa_{V,V'}(v\otimes v')=q^{(\lambda,\mu)}v\otimes v'$ for $v\in V_\lambda, v'\in V'_\mu$.
Set
\begin{align}
\CR_{V,V'}&=\kappa_{V,V'}^{-1}\circ\Xi
\in \End_\BF(V\otimes V'),\\
\CR_{V,V'}^\vee&=\tau_{V,V'}\circ\CR_{V,V'}
\in\Hom_\BF(V\otimes V',V'\otimes V).
\end{align}
For morphisms $f:V_1\to V_2$, $f':V_1'\to V_2'$ in $\Mod^f(U)$ we have
\begin{equation}
\label{eq:R}
(f\otimes f')\circ\CR_{V_1,V_1'}=\CR_{V_2,V_2'}\circ(f\otimes f'),
\quad
(f'\otimes f)\circ\CR_{V_1,V_1'}^\vee=\CR_{V_2,V_2'}^\vee\circ(f\otimes f')
\end{equation}
by definition.

We shall also use the following properties of $\CR_{V,V'}$ (see Drinfeld \cite{D}, Lusztig \cite{Lbook}, Tanisaki \cite{T}).
\begin{proposition}
\label{prop:R}
Let $V, V', V''\in\Mod^f(U)$.
\begin{itemize}
\item[{\rm(i)}]
$\CR_{V,V'}$ is invertible, ant its inverse is given by
\[
\CR_{V,V'}^{-1}=
\left(\sum_{\beta\in Q^+}q^{(\beta,\beta)}(1\otimes k_\beta)(S\otimes\id)(\Xi_\beta)\right)\circ\kappa_{V,V'}.
\]
\item[{\rm(ii)}]
$\CR_{V,V'}^\vee$ is an isomorphism of $U$-modules.
\item[{\rm(iii)}]
The composition of 
\[
V\otimes V'\otimes V''
\xrightarrow{\id_V\otimes\CR_{V',V''}^\vee}
V\otimes V''\otimes V'
\xrightarrow{\CR_{V,V''}^\vee\otimes\id_{V'}}
V''\otimes V\otimes V'
\]
coincides with $\CR_{V\otimes V',V''}^\vee$.
\end{itemize}
\end{proposition}
\subsection{Center}
\label{subsection:center}
We denote by $\Gz$ the center of $U$.
Let $\BF[\Lambda]=\bigoplus_{\lambda\in\Lambda}\BF\e(\lambda)$
be the group algebra of $\Lambda$, and define a linear map
$\zeta:\Gz\to\BF[\Lambda]$ as the composition of 
\[
\Gz\hookrightarrow U\simeq U^-\otimes U^0\otimes U^+
\xrightarrow{\epsilon^-\otimes\kappa\otimes\epsilon^+}\BF\otimes\BF[\Lambda]\otimes\BF\simeq\BF[\Lambda],
\]
where $U\simeq U^-\otimes U^0\otimes U^+$ is given by the multiplication of $U$, 
$\epsilon^\pm:U^\pm\to\BF$ are the restrictions of the comultiplication $\epsilon:U\to\BF$, and $\kappa:U^0\to\BF[\Lambda]$ is the isomorphism given by $\kappa(k_\lambda)=\e(\lambda)$ for any $\lambda\in\Lambda$.
Define shifted actions of the Weyl group $W$ on $\Lambda$ and $\BF[\Lambda]$ by
\begin{align}
&w\circ\lambda=w(\lambda+\rho)-\rho
\qquad(w\in W,\,\,\lambda\in\Lambda),\\
&w\circ\e(\lambda)=q^{(w\lambda-\lambda,\rho)}\e(w\lambda)
\qquad(w\in W,\,\,\lambda\in\Lambda)
\end{align}
respectively.
Note that the action of $W$ on $\Lambda$ is an affine action and the one on $\BF[\Lambda]$ is a linear action.
\begin{proposition}
[see Tanisaki \cite{T}]
\label{prop:center1}
$\zeta:\Gz\to\BF[\Lambda]$ is an injective algebra homomorphism, and its image coincides with 
\[
\BF[\Lambda]^{W\circ}=
\{x\in\BF[\Lambda]\mid w\circ x=x\,\,(w\in W)\}.
\]
\end{proposition}
The isomorphism $\Gz\simeq\BF[\Lambda]^{W\circ}$ induced by $\zeta$ is called the Harish-Chandra isomorphism.

For $\lambda\in\Lambda$ we define an algebra homomorphism
\begin{equation}
\zeta_\lambda:\Gz\to\BF
\end{equation}
as the composition of $\Gz\hookrightarrow\BF[\Lambda]\to\BF$ where $\BF[\Lambda]\to\BF$ is given by $\e(\mu)\mapsto q^{(\lambda,\mu)}$.
The following is result is standard.
\begin{proposition}
\label{prop:center2}
\begin{itemize}
\item[\rm(i)]
For $\lambda_1, \lambda_2\in \Lambda$ we have $\zeta_{\lambda_1}=\zeta_{\lambda_2}$ if and only if $\lambda_2\in W\circ\lambda_1$.
\item[\rm(ii)]
For $\lambda\in\Lambda^+$ and $z\in\Gz$ we have $z|V(\lambda)=\zeta_\lambda(z)\id$.
\item[\rm(iii)]
For $\lambda\in\Lambda$ and $z\in\Gz$ the action of $z$ on $T(\lambda)$ and $T^*(\lambda)$ are given by the multiplication by $\zeta_\lambda(z)$.
\end{itemize}
\end{proposition}

\subsection{Braid group actions}
We set
\[
\exp_t(x)=\sum_{n=0}^\infty\frac{t^{n(n-1)/2}}{[n]_t!}x^n
\in\BQ(t)[[x]].
\]
Note that 
\[
\exp_t(x)^{-1}=\exp_{t^{-1}}(-x).
\]
For $i\in I$ define an operator $T_i$ on $V\in\Mod^f(U)$ by
\begin{align*}
T_i
&=
\exp_{q_i^{-1}}(q_ik_if_i)
\exp_{q_i^{-1}}(-e_i)
\exp_{q_i^{-1}}(q_i^{-1}k_i^{-1}f_i)H_i\\
&=
\exp_{q_i^{-1}}(-q_ik_i^{-1}e_i)
\exp_{q_i^{-1}}(f_i)
\exp_{q_i^{-1}}(-q_i^{-1}k_ie_i)H_i,
\end{align*}
where $H_i$ is the operator on $V$ which acts by $q_i^{(\lambda,\alpha_i^\vee)((\lambda,\alpha_i^\vee)+1)/2}\id$ on $V_\lambda$ for each $\lambda\in\Lambda$.
This operator coincides with Lusztig's operator $T''_{i,1}$ in \cite[5.2]{Lbook}.

We shall use the following result later (see Lusztig \cite[5.3]{Lbook}).
\begin{lemma}
\label{lem:Ti}
Let $V_1, V_2\in\Mod^f(U)$.
As an operator on $V_1\otimes V_2\in\Mod^f(U)$ we have
\begin{align*}
T_i&=
\exp_{q_i}(q_i^{-2}(q_i-q_i^{-1})e_ik_i^{-1}\otimes f_ik_i)
(T_i\otimes T_i)\\
&=(T_i\otimes T_i)
\exp_{q_i}((q_i-q_i^{-1})f_i\otimes e_i)
\end{align*}
\end{lemma}

For $w\in W$ we choose a minimal expression $w=s_{i_1}\cdots s_{i_n}$ and set 
\[
T_w=T_{i_1}\cdots T_{i_n}.
\]
It is known that $T_w$ does not depend on the choice of a minimal expression and that 
\[
T_w(V_\lambda)=V_{w\lambda}\qquad
(V\in\Mod^f(U),\,\,\lambda\in\Lambda).
\]

For $i\in I$ we can define an algebra automorphism $T_i$ of $U$ by
\begin{align*}
&T_i(k_\mu)=k_{s_i\mu}\qquad(\mu\in\Lambda),\\
&T_i(e_j)=
\begin{cases}
\sum_{k=0}^{-a_{ij}}(-1)^kq_i^{-k}e_i^{(-a_{ij}-k)}e_je_i^{(k)}\qquad
&(j\in I,\,\,j\ne i),\\
-f_ik_i\qquad
&(j=i),
\end{cases}\\
&T_i(f_j)=
\begin{cases}
\sum_{k=0}^{-a_{ij}}(-1)^kq_i^{k}f_i^{(k)}f_jf_i^{(-a_{ij}-k)}\qquad
&(j\in I,\,\,j\ne i),\\
-k_i^{-1}e_i\qquad
&(j=i).
\end{cases}
\end{align*}
For $w\in W$ we define an algebra automorphism $T_w$ of $U$ by $T_w=T_{i_1}\cdots T_{i_n}$ where $w=s_{i_1}\cdots s_{i_n}$ is a minimal expression.
The automorphism $T_w$ does not depend on the choice of a minimal expression.
It is known that
\begin{equation}
T_w(uv)=T_w(u)T_w(v)
\qquad(w\in W, u\in U, v\in V\in\Mod^f(U)).
\end{equation}

Let $V\in\Mod^f(U^\op)$.
By $V^*\in\Mod^f(U)$ we can define an operator ${}^tT_w$ on $V$ by
\[
\langle {}^tT_w(v),v^*\rangle=\langle v,T_w(v^*)\rangle\qquad
(v\in V,\,\, v^*\in V^*).
\]
\begin{lemma}
\label{lem:W-action-on-tensor}
Let $w\in W$ and let $\lambda\in\Lambda^+$.
\begin{itemize}
\item[\rm(i)]
Let $V\in\Mod^f(U)$.
For any $v\in V$ and any $\ell\in V(\lambda)_\lambda$ we have 
\[
T_w^{-1}(\ell\otimes v)=T_w^{-1}(\ell)\otimes T_w^{-1}(v),\qquad
T_w(v\otimes\ell)=T_w(v)\otimes T_w(\ell).
\]
\item[\rm(ii)]
Let $V\in\Mod^f(U^\op)$.
For any $v\in V$ and any $\ell\in V^*(\lambda)_\lambda$ we have 
\[
{}^tT_w^{-1}(\ell\otimes v)
={}^tT_w^{-1}(\ell)\otimes {}^tT_w^{-1}(v),\qquad
{}^tT_w(v\otimes \ell)={}^tT_w(v)\otimes {}^tT_w(\ell).
\]
\end{itemize}
\end{lemma}
\begin{proof}
We can easily reduce the proof to the rank one case.
In the rank one case they follow from Lemma \ref{lem:Ti}.
\end{proof}

\subsection{Dual Hopf algebras}
Let $H$ be a Hopf algebra over a field $\BK$.
The dual space $H^*=\Hom_\BK(H,\BK)$ is endowed with a structure of an $H$-bimodule by
\[
\langle h_1f h_2, h\rangle=\langle f, h_2hh_1\rangle\qquad
(h, h_1, h_2\in H, f\in H^*).
\]
The linear maps
$m_H, \eta_H, \Delta_H, \epsilon_H, S_H$ 
induce linear maps
\begin{align}
&m_{H^*}={}^t\Delta_H:(H\otimes_\BK H)^*\to {H^*},
\label{eq:H*1}\\
&\eta_{H^*}={}^t\epsilon_H:\BK\to {H^*},\\
&\Delta_{H^*}={}^tm_H:{H^*}\to (H\otimes_\BK H)^*,\\
&\epsilon_{H^*}={}^t\eta_H:{H^*}\to\BK,\\
&S_{H^*}={}^tS_H:{H^*}\to {H^*}.\label{eq:H*5}
\end{align}
Note that we have $H^*\otimes_\BK H^*\subset(H\otimes_\BK H)^*$.

Let $T$ be a Hopf subalgebra of $H$.
We assume that $T$ is commutative and cocommutative.
The set $\Hom_{\rm alg}(T,\BK)$ of algebra homomorphisms from $T$ to $\BK$ is endowed with a structure of an abelian group by 
\[
(\varphi\psi)(t)=\sum_{(t)_1}\varphi(t_{(1)})\psi(t_{(2)})\qquad
(\varphi, \psi\in\Hom_{\rm alg}(T,\BK),  t\in T).
\]
Assume that we are given a subgroup $\Upsilon$ of $\Hom_{\rm alg}(T,\BK)$.

Denote by $\Mod _\Upsilon(T)$ (resp.\ 
$\Mod _{\Upsilon\times\Upsilon}(T\otimes_\BK T)$) the subcategory of $\Mod(T)$ (resp.\ $\Mod(T\otimes_\BK T$) consisting of finite-dimensional semisimple $T$-modules (resp.\ $T\otimes_\BK T$-modules) whose irreducible factors are contained in $\Upsilon$ (resp.\ $\Upsilon\times\Upsilon$).
Here elements of $\Upsilon$ (resp.\ $\Upsilon\times\Upsilon$) are identified with isomorphism classes of objects of the category $\Mod(T)$ (resp.\ $\Mod(T\otimes_\BK T)$).
\begin{proposition}
\label{prop:dHopf1}
The following conditions on $f\in H^*$ are equivalent.
\begin{itemize}
\item[{\rm(a)}]
$Hf\in\Mod _\Upsilon(T)$.
\item[{\rm(b)}]
$fH\in\Mod _\Upsilon(T)$.
\item[{\rm(c)}]
$HfH\in\Mod _{\Upsilon\times\Upsilon}(T\otimes_\BK T)$.
\item[{\rm(d)}]
There exists a two-sided ideal $I$ of $H$ such that $\langle f, I\rangle=\{0\}$ and 
$H/I\in\Mod _{\Upsilon\times\Upsilon}(T\otimes_\BK T)$.
\end{itemize}
\end{proposition}
\begin{proof}
We have obviously (c)$\Rightarrow$(a).
We obtain (a)$\Rightarrow$(d) by setting $I=\Ker(H\to\End_\BK(Hf))$.
From (d) we obtain (c) as a consequence of $HfH\subset(H/I)^*$.
The implications (c)$\Rightarrow$(b)$\Rightarrow$(d) is proved similarly.
\end{proof}
We denote by $H^*_{T,\Upsilon}$ the set of $f\in H^*$ satisfying the equivalent conditions in Proposition \ref{prop:dHopf1}.
\begin{proposition}
\label{prop:dHopf2(ii)}
A Hopf algebra structure on $H^*_{T,\Upsilon}$ is induced by the linear maps \eqref{eq:H*1},\dots,\eqref{eq:H*5}.
\end{proposition}
\begin{proof}
We need to show 
$m_{H^*}(H^*_{T,\Upsilon}\otimes_\BK H^*_{T,\Upsilon})\subset H^*_{T,\Upsilon}$, 
$\eta_{H^*}(1)\subset H^*_{T,\Upsilon}$, 
$\Delta_{H^*}(H^*_{T,\Upsilon})
\subset H^*_{T,\Upsilon}\otimes_\BK H^*_{T,\Upsilon}$, 
$S_{H^*}(H^*_{T,\Upsilon})\subset H^*_{T,\Upsilon}$.
They are consequences of our assumptions on $T$ and $\Upsilon$.
Details are omitted.
\end{proof}

We  denote by $\Mod_{T,\Upsilon}(H)$ the full subcategory of $\Mod(H)$ consisting of $H$-modules which belong to $\Mod_\Upsilon(T)$ as a $T$-module.
For $M\in\Mod_{T,\Upsilon}(H)$ we define a homomorphism of $H$-bimodules
\begin{equation}
\label{eq:dmatrix-coef}
\Phi_M:M\otimes_\BK M^*\to H^*_{T,\Upsilon}
\end{equation}
by $\langle\Phi_M(v\otimes v^*),h\rangle=\langle v^*, hv\rangle$ for $h\in H, v\in M, v^*\in M^*$.
Elements of $\Image(\Phi_M)$ are called matrix coefficients of the $H$-module $M$.

We denote by $\Mod^{\rm irr}_{T,\Upsilon}(H)$ the set of isomorphism classes of irreducible $H$-modules contained in $\Mod_{T,\Upsilon}(H)$.
\begin{proposition}
\label{prop:dHopf2(i)}
\begin{itemize}
\item[\rm(i)]
We have $H^*_{T,\Upsilon}=\sum_{M\in\Mod_{T,\Upsilon}(H)}\Image(\Phi_M)$.
\item[\rm(ii)]
Assume that $\Mod_{T,\Upsilon}(H)$ is a semisimple category and that $\End_\BK(M)=\BK\id$ for any $M\in\Mod^{\rm irr}_{T,\Upsilon}(H)$.
Then the homomorphism
\[
\bigoplus_{M\in\Mod^{\rm irr}_{T,\Upsilon}(H)}\Phi_M:
\bigoplus_{M\in\Mod^{\rm irr}_{T,\Upsilon}(H)}M\otimes_\BK M^*
\to H^*_{T,\Upsilon}
\]
of $H$-bimodules is an isomorphism.
\end{itemize}
\end{proposition}
\begin{proof}
(i) We have obviously $\Image(\Phi_M)\subset H^*_{T,\Upsilon}$ for $M\in\Mod_{T,\Upsilon}(H)$.
Let $f\in H^*_{T,\Upsilon}$.
Set $M=Hf$. 
Let $\{v_j\}_j$ be a basis of $M$ and let $\{v_j^*\}$ be the dual basis of $M^*$.
For $h\in H$ we have
\[
\langle f,h\rangle=\langle hf,1\rangle=
\sum_j\langle v_j^*,hf\rangle\langle v_j,1\rangle=
\sum_j\langle\Phi_M(f\otimes v_j^*),h\rangle \langle v_j,1\rangle
\]
and hence 
$f=\sum_j\langle v_j,1\rangle\Phi_M(f\otimes v_j^*)
\in\Image(\Phi_M)$.

(ii) By (i) we have 
$H^*_{T,\Upsilon}=\sum_{M\in\Mod ^{\rm irr} _{T,\Upsilon}(H)}\Image(\Phi_M)$.
For $M\in\Mod ^{\rm irr} _{T,\Upsilon}(H)$ $\Image(\Phi_M)$ is a sum of left $H$-modules isomorphic to $M$, and hence we have 
$H^*_{T,\Upsilon}=\bigoplus_{M\in\Mod ^{\rm irr} _{T,\Upsilon}(H)}\Image(\Phi_M)$.
Let $M\in\Mod ^{\rm irr}_{T,\Upsilon}(H)$.
We see that $M\otimes_\BK M^*$ is an irreducible $H$-bimodule from $\End_H(M)=\BK\id$, and hence $\Phi_M$ is injective.
\end{proof}
\subsection{Quantized coordinate algebras}
Set
\begin{equation}
F=U^*_{U^0,\Lambda},\qquad
F^{\geqq0}=\left(U^{\geqq0}\right)^*_{U^0,\Lambda},\qquad
F^0=\left(U^0\right)^*_{U^0,\Lambda},
\end{equation}
where $\Lambda$ is regarded as a subgroup of  $\Hom_{\rm alg}(U^0,\BF)$ by $\lambda\mapsto\chi_\lambda$.
The Hopf algebras $F$, $F^{\geqq0}$, $F^0$ are $q$-analogues of the coordinate algebras of $G$, $B$, $H$ (in the notation of  Section \ref{sec:Intoro}) respectively.

By Proposition \ref{prop:dHopf2(i)} (ii) we have
\begin{align}
&F^0=\bigoplus_{\lambda\in\Lambda}\BF\chi_\lambda,\\
&F\cong\bigoplus_{\lambda\in\Lambda^+}V(\lambda)\otimes V^*(\lambda).
\label{eq:PW}
\end{align}
We see easily that in $F^0$ we have
\[
\chi_\lambda\chi_\mu=\chi_{\lambda+\mu}
\quad(\lambda, \mu\in\Lambda), \qquad
\chi_0=1.
\]
In particular, $F^0$ is isomorphic to the group algebra of $\Lambda$.

The Hopf algebra homomorphisms
\[
U^0\hookrightarrow U^{\geqq0}\hookrightarrow U,\quad
U^{\geqq0}\xrightarrow{\pi^+}U^0
\]
induce Hopf algebra homomorphisms
\begin{equation}
F\xrightarrow{r_+}F^{\geqq0}\xrightarrow{r^{+}_0}F^0,
\qquad
F^0\xrightarrow{i^0_{+}}F^{\geqq0}.
\end{equation}
Since the composition of $U^0\hookrightarrow U^{\geqq0}$ and ${\pi^+}$ is the identity of $U^0$, we have 
\begin{equation}
r^{+}_0\circ i^0_{+}=\id.
\end{equation}
In particular, $r^{+}_0$ is surjective and $i^0_{+}$ is injective.
By $\chi_\lambda^+=i^0_+(\chi_\lambda)$ we have
\begin{equation}
\label{eq:chi}
\chi^+_{\lambda}\chi^+_{\mu}=\chi^+_{\lambda+\mu}\quad(\lambda,\mu\in\Lambda),\qquad
\chi^+_{0}=1
\end{equation}
in $F^{\geqq0}$.

For $\lambda\in\Lambda$ we set
\begin{equation}
F^{\geqq0}(\lambda)=
\{f\in F^{\geqq0}\mid ft=\chi_\lambda(t)f\,\,\,(t\in U^0)\}.
\end{equation}
We see easily that
\[
F^{\geqq0}(\lambda)F^{\geqq0}(\mu)\subset F^{\geqq0}(\lambda+\mu)
\quad(\lambda,\mu\in\Lambda),\qquad
\chi_\lambda\in F^{\geqq0}(\lambda).
\]
In particular, $F^{\geqq0}(0)$ is a subalgebra of $F^{\geqq0}$.

Set
\[
(U^+)^\bigstar=\bigoplus_{\gamma\in Q^+}(U^+_\gamma)^*\subset(U^+)^*.
\]
\begin{proposition}
\label{prop:CqB}
\begin{itemize}
\item[\rm(i)]
The linear map $F^{\geqq0}(0)\otimes F^0\to F^{\geqq0}$ given by $f\otimes f'\mapsto f \,i^0_+(f')$ is an isomorphism.
\item[\rm(ii)]
For $f\in F^{\geqq0}(0)$ we have $f|U^+\in(U^+)^\bigstar$.
Moreover, the linear map $F^{\geqq0}(0)\to(U^+)^\bigstar\,\,(f\mapsto f|U^+)$ is an isomorphism.
\end{itemize}
\end{proposition}
\begin{proof}
For $f\in F^{\geqq0}(0),\,\,f'\in F^0$ we have
\begin{equation}
\label{eq:CqB1}
\langle f\,i^0_+(f'),tu\rangle=\langle f',t\rangle\langle f,u\rangle
\qquad(u\in U^+, t\in U^0).
\end{equation}
Indeed, for $u\in U^+, \lambda\in\Lambda$ we have 
\begin{align*}
&\langle f\,i^0_+(f'),k_\lambda u\rangle=
\langle f\otimes i^0_+(f'),\Delta(k_\lambda u)\rangle=
\langle f\otimes f',(\id\otimes\pi^+)\Delta(k_\lambda u)\rangle\\
=&\langle f\otimes f',k_\lambda u\otimes k_\lambda\rangle
=\langle f',k_\lambda\rangle\langle f,u\rangle.
\end{align*}
by 
\[
(\id\otimes\pi^+)\Delta(u)=u\otimes 1\,\,(u\in U^+),\qquad
(\id\otimes\pi^+)\Delta(k_\lambda)=k_\lambda\otimes k_\lambda\,\,(\lambda\in\Lambda).
\]

Recall that $F^{\geqq0}$ was defined as a subspace of $\left(U^{\geqq0}\right)^*$ .
By \eqref{eq:CqB1} it is sufficient to show that $F^{\geqq0}$ coincides with 
$F^0\otimes(U^+)^\bigstar$ under the identification $\left(U^{\geqq0}\right)^*=(U^0\otimes U^+)^*$.
It is easy to show $F^{\geqq0}\subset F^0\otimes(U^+)^\bigstar$ and 
$F^{\geqq0}\supset F^0\otimes 1$.
Hence it is sufficient to show $F^{\geqq0}\supset \chi_\gamma^+\otimes(U^+_\gamma)^*$ for any $\gamma\in Q^+$.
Define $M\in\Mod_{U^0,\Lambda}(U^{\geqq0})$ by
\[
N=U^{\geqq0}/\sum_{\lambda\in\Lambda}U^{\geqq0}(k_\lambda-1),\qquad
M=N/\sum_{\gamma-\gamma'\notin Q^+}N_{\gamma'}.
\]
Then the elements of $\chi_\gamma^+\otimes(U^+_\gamma)^*$ are obtained as matrix coefficients of the $U^{\geqq0}$-module $M$.
\end{proof}
By Proposition \ref{prop:CqB} (i) we obtain
\begin{equation}
\label{eq:CqB}
F^{\geqq0}=\bigoplus_{\lambda\in\Lambda}F^{\geqq0}(\lambda),\qquad
F^{\geqq0}(\lambda)=F^{\geqq0}(0)\chi^+_\lambda\quad(\lambda\in\Lambda).
\end{equation}
\begin{proposition}
\label{prop:CqGCqB}
The Hopf algebra homomorphism $r_+:F\to F^{\geqq0}$ is surjective.
\end{proposition}
\begin{proof}
Identify $F^{\geqq0}$ with $F^0\otimes(U^+)^\bigstar$ by Proposition \ref{prop:CqB}.
We see easily that $\chi_\lambda\otimes1\in\Image(r_+)$ for any $\lambda\in\Lambda$.
Hence it is sufficient to show that for any $\gamma\in Q^+$  we have
$\Image(r_+)\supset \chi_{-\lambda+\gamma}\otimes(U^+_\gamma)^*$ for sufficiently large $\lambda\in \Lambda^+$.
Let $v_{-\lambda}$ be a non-zero element of $V(w_0\lambda)_{-\lambda}$.
By Lemma \ref{lem:V} (ii) the linear map 
$U^+_{\gamma}\ni u\mapsto uv_{-\lambda}\in V(-w_0\lambda)_{-\lambda+\gamma}$ is bijective when $\lambda$ is sufficiently large.
Hence elements of $\chi_{-\lambda+\gamma}\otimes(U^+_\gamma)^*$ are obtained as the matrix coefficients of the $U^{\geqq0}$-module $V(-w_0\lambda)$.
\end{proof}

\section{Quantized flag manifolds}
\subsection{Homogeneous coordinate algebras}
We set
\begin{equation}
A=\{\varphi\in F\mid\varphi u=\epsilon(u)\varphi\,\,(u\in U^-)\}.
\end{equation}
\begin{lemma}
\label{lem:prpperties of A}
\begin{itemize}
\item[{\rm(i)}]
$A$ is a subalgebra of $F$.
\item[{\rm(ii)}]
$A$ is a left $U$-submodule of $F$.
\item[{\rm(iii)}]
We have $\Delta_{F}(A)\subset A\otimes F$.
\item[{\rm(iv)}]
The multiplication $A\otimes A\to A$ is a homomorphism of $U$-modules; i.e.
\begin{equation}
u(\varphi_0\varphi_1)
=\sum_{(u)_1}(u_{(0)}\varphi_0)(u_{(1)}\varphi_1)\qquad
(\varphi_0, \varphi_1\in A,\,u\in U).
\end{equation}
\end{itemize}
\end{lemma}
\begin{proof}
(i) For $\varphi\in F$ we have $\varphi\in A$ if and only if $\varphi f_i=0$ for any $i\in I$.
We have
\[
\langle 1f_i, u\rangle
=\langle 1, f_iu\rangle
=\epsilon(f_iu)=0
\]
for any $u\in U$, and hence $1\in A$.
For $\varphi_1, \varphi_2\in A$ we have
\begin{align*}
\langle(\varphi_1\varphi_2)f_i,u\rangle
&=\langle\varphi_1\varphi_2,f_iu\rangle
=\langle\varphi_1\otimes\varphi_2,\Delta(f_iu)\rangle\\
&=\langle\varphi_1\otimes\varphi_2,
(f_i\otimes k_i^{-1}+1\otimes f_i)\Delta(u)\rangle\\
&=\langle\varphi_1f_i\otimes\varphi_2k_i^{-1},\Delta(u)\rangle
+\langle\varphi_1\otimes\varphi_2f_i,\Delta(u)\rangle=0
\end{align*}
for any $i\in I$ and $u\in U$, and hence $\varphi_1\varphi_2\in A$.

The statement (ii) is obvious from the definition.

For $\varphi\in A, u\in U^+, u_1, u_2\in U$ we have
\begin{align*}
&\langle(\Delta(\varphi))(u\otimes 1),u_1\otimes u_2\rangle
=\langle\Delta(\varphi),uu_1\otimes u_2\rangle
=\langle\varphi,uu_1u_2\rangle
=\langle\varphi u,u_1u_2\rangle\\
=&\epsilon(u)\langle\varphi,u_1u_2\rangle
=\epsilon(u)\langle\Delta(\varphi), u_1\otimes u_2\rangle
\end{align*}
and hence $(\Delta(\varphi))(u\otimes 1)=\epsilon(u)\Delta(\varphi)$.
It follows that $\Delta(\varphi)\in A\otimes F$.
The statement (iii) is proved.

(iv) For $u'\in U$ we have
\[
\begin{split}
&\langle u(\varphi_0\varphi_1),u'\rangle
=\langle \varphi_0\varphi_1,u'u\rangle
=\langle \varphi_0\otimes\varphi_1,\Delta(u')\Delta(u)\rangle\\
=&\sum_{(u)_1}
\langle u_{(0)}\varphi_0\otimes u_{(1)}\varphi_1,\Delta(u')\rangle
=\sum_{(u)_1}
\langle (u_{(0)}\varphi_0)(u_{(1)}\varphi_1),u'\rangle.
\end{split}
\]
\end{proof}

By Lemma \ref{lem:prpperties of A} (iii) we obtain an algebra homomorphism
\begin{equation}
\overline{\Delta}:A\to A\otimes F
\end{equation}
by restricting $\Delta_{F}$ to $A$.

For $\lambda\in\Lambda$ we set
\begin{equation}
A(\lambda)=\{\varphi\in F\mid\varphi u=\chi_\lambda^-(u)\varphi\,\,(u\in U^{\leqq0})\}.
\end{equation}
By \eqref{eq:right-hw} we have 
\[
A\simeq 
\bigoplus_{\lambda\in\Lambda^+}V(\lambda)\otimes V^*(\lambda)_\lambda
\]
under the isomorphism \eqref{eq:PW}.
Hence we have
\begin{equation}
A(\lambda)\simeq
\begin{cases}
V(\lambda)\quad&(\lambda\in\Lambda^+)\\
0&(\mbox{otherwise})
\end{cases}
\end{equation}
as a $U$-module, and
\begin{equation}
A=\bigoplus_{\lambda\in\Lambda^+}A(\lambda).
\end{equation}

For each $\lambda\in\Lambda^+$ we fix 
\begin{equation}
v^*_\lambda\in V^*(\lambda)_\lambda\setminus\{0\}.
\end{equation}
Then an isomorphism 
\begin{equation}
f_\lambda:V(\lambda)\to A(\lambda)
\end{equation}
of $U$-modules is defined by 
\[
\langle f_\lambda(v),u\rangle=\langle v^*_\lambda, uv\rangle\quad(v\in V(\lambda),\, u\in U).
\]

Let $\lambda, \mu\in\Lambda^+$.
Then the $U$-module $V(\lambda)\otimes V(\mu)$ 
(resp.\ the right $U$-module $V^*(\lambda)\otimes V^*(\mu)$) contains $V(\lambda+\mu)$ (resp.\ $V^*(\lambda+\mu)$) with multiplicity one. 
Let
\[
i_{\lambda,\mu}:V^*(\lambda+\mu)\to V^*(\lambda)\otimes V^*(\mu)
\]
be the embedding of $U$-modules such that 
$i_{\lambda,\mu}(v^*_{\lambda+\mu})=v^*_\lambda\otimes v^*_\mu$, and denote the corresponding projection by 
\begin{equation}
p_{\lambda,\mu}:V(\lambda)\otimes V(\mu)\to V(\lambda+\mu).
\end{equation}

\begin{lemma}
\label{lem:grading}
For $\lambda,\mu\in\Lambda^+$ we have
\[
f_\lambda(v_0)f_\mu(v_1)=
f_{\lambda+\mu}(p_{\lambda,\mu}(v_0\otimes v_1))
\qquad(v_0\in V(\lambda),\,\,v_1\in V(\mu)).
\]
In particular, the multiplication of the algebra $A$ induces a surjective homomorphism
$A(\lambda)\otimes A(\mu)\to A(\lambda+\mu)$ of $U$-modules.
\end{lemma}
\begin{proof}
For $u\in U$ we have
\begin{align*}
&\langle f_\lambda(v_0)f_\mu(v_1), u\rangle\\
=&\langle f_\lambda(v_0) \otimes f_\mu(v_1), \Delta(u)\rangle
=\sum_{(u)_1}\langle f_\lambda(v_0), u_{(0)}\rangle\langle f_\mu(v_1),  u_{(1)}\rangle\\
=&\sum_{(u)_1}\langle v^*_\lambda, u_{(0)}v_0\rangle\langle v^*_\mu, u_{(1)}v_1\rangle
=\langle v^*_\lambda\otimes v^*_\mu,u(v_0\otimes v_1)\rangle\\
=&\langle i_{\lambda,\mu}(v^*_{\lambda+\mu}),u(v_0\otimes v_1)\rangle
=\langle v^*_{\lambda+\mu},p_{\lambda,\mu}(u(v_0\otimes v_1))\rangle
=\langle v^*_{\lambda+\mu},up_{\lambda,\mu}(v_0\otimes v_1)\rangle\\
=&\langle f_{\lambda+\mu}(p_{\lambda,\mu}(v_0\otimes v_1)),u\rangle.
\end{align*}
Hence we have 
$f_\lambda(v_0)f_\mu(v_1)
=f_{\lambda+\mu}(p_{\lambda,\mu}(v_0\otimes v_1))$.
\end{proof}

Hence $A$ is a $\Lambda$-grade $\BF$-algebra with $A(0)=\BF1$.

By Joseph \cite{J:emb} we have the following.
\begin{proposition}
[Joseph]
\label{prop:property-A}
\begin{itemize}
\item[{\rm(i)}]
$A$ is a domain, i.e., if $\varphi\psi=0$ for $\varphi,\psi\in A$, then we have $\varphi=0$ or $\psi=0$.
\item[{\rm(ii)}]
$A$ is left and right noetherian.
\end{itemize}
\end{proposition}
For a ring (resp.\ $\Lambda$-graded ring) $R$ we denote by $\Mod(R)$ (resp.\ $\Mod_\Lambda(R)$) the category of left $R$-modules (resp.\ $\Lambda$-graded left $R$-modules).
For $M\in\Mod_\Lambda(R)$ and $\nu\in\Lambda$ we define 
$M[\nu]\in\Mod_\Lambda(R)$ by
\[
(M[\nu])(\lambda)=M(\lambda+\nu).
\]

\subsection{Category of quasi-coherent sheaves}
For $M\in\Mod_\Lambda(A)$ we denote by $\Tor(M)$ the graded $A$-submodule of $M$ consisting of $m\in M$ such that $A(\lambda)m=\{0\}$ for sufficiently large $\lambda\in\Lambda^+$.
Let $\Tor_\Lambda(A)$ be the full subcategory of $\Mod_\Lambda(A)$ consisting of $M\in\Mod_\Lambda(A)$ satisfying $M=\Tor(M)$.
Note that $\Tor_\Lambda(A)$ is closed under taking subquotients and extensions in $\Mod_\Lambda(A)$.
Let $\Sigma$ denote the collection of morphisms $f$ in $\Mod_\Lambda(A)$ satisfying $\Ker(f), \Coker(f)\in\Tor_\Lambda(A)$.
Then we define the abelian category $\CM(A)$ of ``quasi-coherent sheaves" on the ``quantized flag manifold $\CB_q$" by
\begin{equation}
\CM(A)
=\frac{\Mod_\Lambda(A)}{\Tor_\Lambda(A)}
=\Sigma^{-1}\Mod_\Lambda(A)
\end{equation}
(see Gabriel-Zisman \cite{GZ} and Popescu \cite{P} for the notion of localization of categories).
For $\nu\in\Lambda$ we denote by
\[
\CM(A)\ni M\mapsto M[\nu]\in\CM(A)
\]
the exact functor induced by $\Mod_\Lambda(A)\ni M\mapsto M[\nu]\in\Mod_\Lambda(A)$.

Let
\begin{equation}
\omega^*:\Mod_\Lambda(A)\to\CM(A)
\end{equation}
be the canonical localization functor.
We have the following by the definition.
\begin{lemma}
\begin{itemize}
\item[\rm(i)]
$\omega^*$ is an exact functor.
\item[\rm(ii)]
Let $f$ be a morphism in $Mod_\Lambda(A)$.
Then $\omega^*f$ is an isomorphism if and only if the kernel and the cokernel of $f$ belong to $\Tor_\Lambda(A)$.
\end{itemize}
\end{lemma}

By Popescu \cite[Ch4, Corollary 6.2]{P} we have the following.

\begin{proposition}
\label{prop:injectives}
The abelian category $\CM(A)$ has enough injectives.
\end{proposition}

It is shown using Proposition \ref{prop:injectives} that there exists an additive functor 
\begin{equation}
\omega_*:\CM(A)\to\Mod_\Lambda(A)
\end{equation}
which is right adjoint to $\omega^*$ (see Popescu \cite[Ch4, Proposition 5.2]{P}).
Note that $\omega_*$ is left exact.
By Popescu \cite[Ch4, Proposition 4.3]{P} we have the following.
\begin{proposition}
\label{prop:omega-omega}
The canonical morphism
$\omega^*\circ\omega_*\to\Id$ is an isomorphism.
\end{proposition}
\begin{corollary}
\label{cor:ch-of-omega-omega}
Let $M\in\Mod_\Lambda(A)$.
Set $N=\omega_*\omega^*M$ and let $f:M\to N$ be the canonical morphism.
Then $N$ and $f$ are uniquely characterized by the following properties.
\begin{itemize}
\item[{\rm(a)}]
$\Ker(f)$ and $\Coker(f)$ belong to $\Tor_\Lambda(A)$.
\item[{\rm(b)}]
$\Tor(N)=\{0\}$.
\item[{\rm(c)}]
Any monomorphism $N\to L$ with $L/N\in\Tor_\Lambda(A)$ is a split morphism.
\end{itemize}
\end{corollary}
\begin{proof}
Let $f:M\to N=\omega_*\omega^*M$ be the canonical morphism.
We have also a canonical morphism $g:\omega^*\omega_*(\omega^*M)\to\omega^*M$, and the composition $g\circ \omega^*f:\omega^*M\to\omega^*M$ is equal to $\id_{\omega^*M}$ by the definition of the adjoint functors.
Since $g$ is an isomorphism by Proposition \ref{prop:omega-omega}, $\omega^*f$ is also an isomorphism. This implies (a).
If $T$ is a subobject of $N=\omega_*\omega^*M$ belonging to $\Tor_\Lambda(A)$, then we have
\[
\Hom(T,N)=\Hom(T,\omega_*\omega^*M)\simeq
\Hom(\omega^*T,\omega^*M)=\Hom(0,\omega^*M)=0,
\]
and hence $T=0$.
The statement (b) is proved.
To show (c) it is sufficient to show that the homomorphism $\Hom(L,N)\to \Hom(N,N)$ induced by $N\to L$ is surjective.
By $\omega^*M\simeq\omega^*N\simeq\omega^*L$ we have 
\[
\Hom(L,N)\simeq\Hom(\omega^*L,\omega^*M)
\simeq\Hom(\omega^*M,\omega^*M)
\]
and similarly $\Hom(N,N)\simeq\Hom(\omega^*M,\omega^*M)$.
Hence we have $\Hom(L,N)\simeq\Hom(N,N)$.

Assume that the conditions (a), (b), (c) are satisfied for some $f:M\to N$.
By (a) $\omega^*f$ is an isomorphism.
Let $h:N\to\omega_*\omega^*N$ be the canonical morphism, and define
$g:N\to\omega_*\omega^*M$ by $g=(\omega_*\omega^*f)^{-1}\circ h$.
Then the kernel and the cokernel of $g$ belong to $\Tor_\Lambda(A)$.
Hence by (b) we have $\Ker(g)=0$.
By applying (c) to $g$ we see that $\Coker(g)$ is isomorphic to a subobject of $\omega_*\omega^*M$, and hence $\Coker(g)=0$.
\end{proof}

We define the ``global section functor "
\begin{equation}
\Gamma:\CM(A)\to\Mod(\BF).
\end{equation}
by $\Gamma(M)=(\omega_*M)(0)$.
Note that $\Gamma$ is left exact.
Its right derived functors are denoted by
\begin{equation}
H^k=R^k\Gamma:\CM(A)\to\Mod(\BF).
\end{equation}

\subsection{Affine open covering}
For each $w\in W$ and $\lambda\in\Lambda^+$ we fix a non-zero element $c^w_\lambda$ of $A(\lambda)_{w^{-1}\lambda}$.
Note that we have $\dim A(\lambda)_{w^{-1}\lambda}=1$.
By Lemma \ref{lem:grading} we have $c^w_\lambda c^w_\mu\in\BF^\times c^w_{\lambda+\mu}$ for any $w\in W$ and $\lambda, \mu\in\Lambda^+$, and hence 
\begin{equation}
S_w=\bigcup_{\lambda\in\Lambda^+}\BF^\times c^w_\lambda
\end{equation}
is a multiplicative subset of $A$ for any $w\in W$.
Moreover, we have the following.
\begin{proposition}
[Joseph \cite{J:emb}]
\label{prop:OreA}
Let $w\in W$.
\begin{itemize}
\item[{\rm(i)}]
$S_w$ satisfies the left and right Ore conditions in $A$.
\item[{\rm(ii)}]
The canonical homomorphism $A\to S_w^{-1}A$ is injective.
\end{itemize}
\end{proposition}
We shall give a proof of Proposition \ref{prop:OreA} different from the one in \cite{J:emb}.
We need the following.
\begin{lemma}
\label{lem:lem}
Let $w\in W, \,\,\mu\in\Lambda^+$, and fix 
$v_{w^{-1}\mu}\in V(\mu)_{w^{-1}\mu}\setminus\{0\}.$
\begin{itemize}
\item[\rm(i)]
For any $\lambda\in\Lambda^+$ 
\[
p_{\lambda,\mu}|V(\lambda)\otimes v_{w^{-1}\mu}:
V(\lambda)\otimes v_{w^{-1}\mu}\to V(\lambda+\mu)
\]
is injective.
\item[\rm(ii)]
Let $\gamma\in Q^+$.
For sufficiently large $\lambda\in\Lambda^+$ we have 
\[
p_{\lambda,\mu}
(V(\lambda)_{w^{-1}(\lambda-\gamma)}\otimes v_{w^{-1}\mu})
=V(\lambda+\mu)_{w^{-1}(\lambda+\mu-\gamma)}.
\]
\end{itemize}
\end{lemma}
\begin{proof}
By Lemma \ref{lem:W-action-on-tensor} we may assume that $w=1$.

(i) Let $M$ and $N$ be $U$-submodules of $V(\lambda)\otimes V(\mu)$ such that $M\simeq V(\lambda+\mu)$ and $V(\lambda)\otimes V(\mu)=M\oplus N$.
It is sufficient to show 
$(V(\lambda)_{\lambda-\gamma}\otimes v_\mu)
\cap N=\{0\}$.
for any $\gamma\in Q^+$.
We shall show this by induction on $\gamma$.
It is obvious for $\gamma=0$.
Assume $\gamma\in Q^+\setminus\{0\}, v\in V(\lambda)_{\lambda-\gamma}, v\otimes v_\mu\in N$.
Then we have $e_iv\otimes v_\mu=e_i(v\otimes v_\mu)\in N$ for any $i\in I$.
Thus we have $e_iv=0$ for any $i\in I$ by the hypothesis of induction.
Hence we obtain $v=0$.

(ii) By (i) the linear map
\[
p_{\lambda,\mu}|V(\lambda)_{\lambda-\gamma}\otimes v_{\mu}:
V(\lambda)_{\lambda-\gamma}\otimes v_{\mu}
\to V(\lambda+\mu)_{\lambda+\mu-\gamma}.
\]
is injective.
Hence it is sufficient to show $\dim V(\lambda)_{\lambda-\gamma}=\dim V(\lambda+\mu)_{\lambda+\mu-\gamma}$ when $\lambda$ is sufficiently large.
This follows from Lemma \ref{lem:V}.
\end{proof}

\noindent
{\sc Proof of Proposition \ref{prop:OreA}.}
It is sufficient to show the following.
\begin{itemize}
\item[(a)]
For any $\varphi\in A$ and $s\in S_w$ there exists some $t\in S_w$ and $\psi\in A$ satisfying $t\varphi=\psi s$.
\item[(b)]
For any $\varphi\in A$ and $s\in S_w$ there exists some $t\in S_w$ and $\psi\in A$ satisfying $\varphi t=s\psi$.
\item[(c)]
If $\varphi s=0$ for $\varphi\in A$ and $s\in S_w$, then we have $\varphi=0$.
\item[(d)]
If $s\varphi=0$ for $\varphi\in A$ and $s\in S_w$, then we have $\varphi=0$.
\end{itemize}
Let us first show (a).
We may assume $s=c^w_\lambda$ and $\varphi=f_\eta(m)$ for some $\lambda, \eta\in\Lambda^+$ and $m\in V(\eta)$.
We may further assume that $m\in V(\eta)_{w^{-1}(\eta-\gamma)}$ for some $\eta\in Q^+$.
Then we need to find $\mu, \xi\in\Lambda^+$ and $n\in V(\xi)$ such that 
$c^w_\mu f_\eta(m)=f_\xi(n)c^w_\lambda$.
By Lemma \ref{lem:grading} it is sufficient to show the existence of 
$\mu, \xi\in\Lambda^+$ and $n\in V(\xi)$ such that 
$\mu+\eta=\xi+\lambda$ and 
$p_{\mu,\eta}(v_{w^{-1}\mu}\otimes m)=
p_{\xi,\lambda}(n\otimes v_{w^{-1}\lambda})$,
where $v_{w^{-1}\lambda}$ and $v_{w^{-1}\mu}$ are non-zero elements of 
$V(\lambda)_{w^{-1}\lambda}$ and $V(\mu)_{w^{-1}\mu}$ respectively.
Take sufficiently large $\xi\in\Lambda^+$ and set 
$\mu=\xi+\lambda-\eta\in\Lambda^+$.
By $p_{\mu,\eta}(v_{w^{-1}\mu}\otimes m)\in V(\xi+\lambda)_{w^{-1}(\xi+\lambda-\gamma)}$ the assertion follows from Lemma \ref{lem:lem} (ii).

Let us show (c).
We may assume that there exists 
$\lambda, \mu\in\Lambda^+, \gamma\in Q^+$ and 
$m\in V(\mu)_{w^{-1}(\mu-\gamma)}$ such that
$s=c^w_\lambda$ and $\varphi=f_\mu(m)$.
By Lemma \ref{lem:grading} we have
$p_{\mu,\lambda}(m\otimes v_{w^{-1}\lambda})=0$, 
where $v_{w^{-1}\lambda}$ is a non-zero element of 
$V(\lambda)_{w^{-1}\lambda}$.
Then we obtain $m=0$ by Lemma \ref{lem:lem} (i).
Hence $\varphi=0$

The statements (b) and (d) are proved similarly.
\hfill$\Box$
\medskip

Since $S_w$ consists of homogeneous elements, $S_w^{-1}A$ is a $\Lambda$-graded $\BF$-algebra.
We define an $\BF$-algebra $R_w$ by 
\begin{equation}
R_w=(S_w^{-1}A)(0).
\end{equation}
Note that for any $\lambda\in\Lambda$ $S_w^{-1}A(\lambda)$ is a free left (right) $R_w$-module generated by $c^w_{\lambda_1}(c^w_{\lambda_2})^{-1}$ where $\lambda=\lambda_1-\lambda_2$ with $\lambda_1,\lambda_2\in\Lambda^+$.
In particular, we have an equivalence of categories 
\[
\Mod_\Lambda(S_w^{-1}A)\simeq \Mod(R_w)
\]
given by 
\begin{align*}
&\Mod_\Lambda(S_w^{-1}A)
\ni M\mapsto M(0)\in
\Mod(R_w),\\
&\Mod(R_w)
\ni N\mapsto S_w^{-1}A\otimes_{R_w}N\in
\Mod_\Lambda(S_w^{-1}A).
\end{align*}
For any $M\in\Tor_\Lambda(A)$ we have $S_w^{-1}M=0$ by the definition of $\Tor_\Lambda(A)$, and hence the localization functor $\Mod_\Lambda(A)\to\Mod_\Lambda(S_w^{-1}A)$ induces an exact functor 
\begin{equation}
j_w^*:\CM(A)\to\Mod(R_w)
\end{equation}
for any $w\in W$.
$\Mod(R_w)$ is regarded as the category of ``quasi-coherent sheaves on the affine open subset $U_{w,q}=\Spec(R_w)$ of $\CB_q$".

\begin{proposition}[Lunts-Rosenberg \cite{LR}]
\label{prop:quasi-scheme}
$\CM(A)$ is a quasi-scheme with Zariski cover $\{j_w^*\mid w\in W\}$ in the sense of Rosenberg \cite{R}.
In particular, a morphism $f$ in $\CM(A)$ is an isomorphism if and only if $j_w^*f$ is an isomorphism for any $w\in W$.
\end{proposition}
\begin{remark}
{\rm
An essential part in the proof of Proposition \ref{prop:quasi-scheme} is to show the following fact:
For any $\mu\in\Lambda^+$ one has 
\[
\sum_{w\in W}A(\lambda)c^w_\mu=A(\lambda+\mu)
\]
for any sufficiently large $\lambda\in\Lambda^+$.
The proofs of this fact given in Joseph \cite{J:emb} and Lunts-Rosenberg \cite{LR} both use the reduction to the case $q=1$.
One needs another proof in order to define the quantized flag manifold at roots of unity as a quasi-scheme.
}
\end{remark}
\begin{corollary}
\label{cor:cor-quasi-scheme}
For $M\in\Mod_\Lambda(A)$ the canonical homomorphism
\[
\omega_*\omega^*M\to\bigoplus_{w\in W}S_w^{-1}M
\]
is injective.
\end{corollary}

By Proposition \ref{prop:quasi-scheme} one can use the general results in Rosenberg \cite{R} for quasi-schemes.
Especially we have a description of the cohomology groups in terms of certain \v{C}ech cohomology groups.

Using the \v{C}ech cohomology groups and the arguments involving the reduction to the case $q=1$ Lunts-Rosenberg \cite[III, \S4]{LR} obtained the following analogues of Serre's theorem and the Borel-Weil theorem.
\begin{proposition}
[Lunts-Rosenberg \cite{LR}]
\label{prop:Serre}
Let $f:M\to N$ be an epimorphism in $\Mod_\Lambda(A)$.
Assume that $M, N$ are finitely-generated as $A$-modules.
Then the homomorphism
\[
\Gamma(\omega^*M[\lambda])\to\Gamma(\omega^*N[\lambda])
\]
is surjective for sufficiently large $\lambda\in\Lambda^+$.
\end{proposition}
\begin{proposition}
[Lunts-Rosenberg \cite{LR}]
\label{prop:Borel-Weil}
For $\lambda\in\Lambda$ we have
\[
\Gamma(A[\lambda])=
\begin{cases}
A(\lambda)\quad&(\lambda\in\Lambda^+)\\
0&(\lambda\notin\Lambda^+).
\end{cases}
\]
In other words we have $A\simeq\omega_*\omega^*A$.
\end{proposition}

\subsection{Schubert varieties}
The contents of this subsection will not be used in the sequel.

Let $w\in W$.
Since $F$ is a sum of finite-dimensional right $U$-submodules contained in $\Mod^f(U^\op)$, we can define the operator
${}^tT_w:F\to F$.
Let $\tilde{\epsilon}_w:F\to\BF$ be the linear map defined by 
$\tilde{\epsilon}_w(\varphi)=\langle{}^tT_w(\varphi),1\rangle$ for $\varphi\in F$, and define 
\begin{equation}
{\epsilon}_w:A\to\BF
\end{equation}
to be the composition of the inclusion $A\hookrightarrow F$ with ${\tilde{\epsilon}_w}$.
\begin{lemma}
\label{lem:ew}
$\epsilon_w$ is an algebra homomorphism, and we have $\epsilon_w(S_w)\subset\BF^\times$.
\end{lemma}
\begin{proof}
By ${}^tT_w(1)=1$ we have $\epsilon_w(1)=1$.
Let us show $\epsilon_w(\varphi\psi)=\epsilon_w(\varphi)\epsilon_w(\psi)$ for $\varphi, \psi\in A$.
We may assume that $\varphi\in A(\lambda), \psi\in A(\mu)$ for $\lambda, \mu\in\Lambda^+$.
Then we have $\varphi=f_\lambda(v_1),\,\,\psi=f_\mu(v_2)$ for some $v_1\in V(\lambda), v_2\in V(\mu)$.
For $u\in U$ we have 
\[
\langle\varphi,u\rangle
=\langle v^*_\lambda, uv_1\rangle,\quad
\langle\psi,u\rangle
=\langle v^*_\mu, uv_2\rangle,\quad
\langle\varphi\psi,u\rangle
=\langle v^*_\lambda\otimes v^*_\mu, u(v_1\otimes v_2)\rangle,
\]
and hence
\begin{align*}
&\epsilon_w(\varphi)=\langle {}^tT_w(v^*_\lambda), v_1\rangle, \qquad
\epsilon_w(\psi)=\langle {}^tT_w(v^*_\mu), v_2\rangle,\\
&\epsilon_w(\varphi\psi)=\langle {}^tT_w(v^*_\lambda\otimes v^*_\mu), v_1\otimes v_2\rangle.
\end{align*}
Hence we obtain $\epsilon_w(\varphi\psi)=\epsilon_w(\varphi)\epsilon_w(\psi)$ by Lemma \ref{lem:W-action-on-tensor}.
For $\lambda\in\Lambda^+$ take 
$v_{w^{-1}\lambda}\in V(\lambda)_{w^{-1}\lambda}\setminus\{0\}$ 
such that $c_\lambda^w=f_\lambda(v_{w^{-1}\lambda})$.
Then we have 
$\langle c^w_\lambda, u\rangle=
\langle v^*_\lambda, uv_{w^{-1}\lambda}\rangle$
for $u\in U$, 
and hence 
\[
\epsilon_w(c^w_\lambda)=
\langle{}^tT_w(v^*_\lambda), v_{w^{-1}\lambda}\rangle
=\langle v^*_\lambda, T_w (v_{w^{-1}\lambda})\rangle
\ne0.
\]
\end{proof}
By $\epsilon_w(S_w)\subset\BF^\times$ the algebra homomorphism $\epsilon_w$ is uniquely extended to 
\begin{equation}
{\epsilon}_w:S_w^{-1}A\to\BF.
\end{equation}

We define an algebra homomorphism 
\begin{equation}
\Phi_w:A\to F^{\geqq0}
\end{equation}
as the composition of 
\[
A\xrightarrow{\overline{\Delta}}A\otimes F
\xrightarrow{{\epsilon}_w\otimes r_+}\BF\otimes F^{\geqq0}=F^{\geqq0}.
\]
\begin{lemma}
\label{lem:SV}
\begin{itemize}
\item[\rm(i)]
Elements of $\Phi_w(S_w)$ are invertible in $F^{\geqq0}$.
\item[\rm(ii)]
We have $\Phi_w(A(\lambda))\subset F^{\geqq0}(w^{-1}\lambda)$ for any $\lambda\in\Lambda$.
\end{itemize}
\end{lemma}
\begin{proof}
(i) Let $\lambda\in\Lambda^+$.
Take 
$v_{w^{-1}\lambda}\in V(\lambda)_{w^{-1}\lambda}\setminus\{0\}$ 
such that $c_\lambda^w=f_\lambda(v_{w^{-1}\lambda})$.
Then for $u\in U$ and $x\in U^{\geqq0}$ we have
\[
\langle
(\id\otimes r_+)\circ\overline{\Delta}(c^w_\lambda),
u\otimes x\rangle
=\langle c^w_\lambda, ux\rangle
=\langle v^*_\lambda, uxv_{w^{-1}\lambda}\rangle.
\]
Therefore, we have
\[
\langle
\Phi_w(c^w_\lambda),
x\rangle
=\langle {}^tT_w(v^*_\lambda), xv_{w^{-1}\lambda}\rangle
=\chi^+_{w^{-1}\lambda}(x)
\langle{}^tT_w(v_\lambda^*),v_{w^{-1}\lambda}\rangle.
\]
By ${}^tT_w(v^*_\lambda)\in V^*(\lambda)_{w^{-1}\lambda}\setminus\{0\}$ we obtain 
$\Phi_w(c^w_\lambda)
=\chi^+_{w^{-1}\lambda}$ up to a non-zero constant multiple.
The statement (i) is proved.

The proof of (ii) is similar and omitted.
\end{proof}

By Lemma \ref{lem:SV} (ii) $\Ker(\Phi_w)$ is a graded ideal of $A$.
The $\Lambda$-graded $\BF$-algebra
\begin{equation}
A_w=A/\Ker(\Phi_w)
\end{equation}
is a $q$-analogue of the homogeneous coordinate algebra of the Schubert variety corresponding to $w\in W$.

\section{Quasi-coherent sheaves with $U$-actions}
\subsection{The algebra $\tU$}
Recall that $A$ is a left $U$-module by Lemma \ref{lem:prpperties of A} (ii).
We sometimes write this action of $U$ on $A$ by
\[
U\otimes A\to A\qquad(u\otimes \varphi\mapsto\deru_u(\varphi)).
\]
We have also a left $A$-module structure on $A$ given by the left multiplication. 
By Lemma \ref{lem:prpperties of A} (iv) we have
\[
\deru_u(\varphi\psi)=\sum_{(u)_1}\deru_{u_{(0)}}(\varphi)\deru_{u_{(1)}}(\psi)\qquad
(u\in U,\,\,\varphi, \psi\in A),
\]
and hence $A$ is a module over the $\BF$-algebra $\tU$  generated by the elements 
$\{\overline{a}\mid a\in A\}\cup\{\overline{u}\mid u\in U\}$ 
satisfying the  fundamental relations:
\begin{align}
&
\overline{\varphi_1}\,\overline{\varphi_2}=\overline{\varphi_1\varphi_2}
\qquad
&(\varphi_1, \varphi_2\in A),
\label{eq:tU1}\\
&
\overline{u_1}\,\overline{u_2}=\overline{u_1u_2}
\qquad
&(u_1, u_2\in U),
\label{eq:tU2}\\
&
\overline{u}\,\overline{\varphi}=
\sum_{(u)_1}
\overline{\deru_{u_{(0)}}(\varphi)}\,\overline{u_{(1)}}
\qquad
&(u\in U,\varphi\in A).
\label{eq:tU3}
\end{align}

For $u\in U$ and $\varphi\in A$ we have
\[
\sum_{(u)_1}\overline{u_{(1)}}\overline{\deru_{S^{-1}u_{(0)}}(\varphi)}
=\sum_{(u)_2}\overline{\deru_{u_{(1)}S^{-1}u_{(0)}}(\varphi)}
\overline{u_{(2)}}
=\sum_{(u)_1}\epsilon(u_{(0)})\,\overline{\varphi}
\overline{u_{(1)}}
=\overline{\varphi}\,
\overline{u}
\]
by \eqref{eq:tU3}, and hence
\begin{equation}
\label{eq:tU10}
\overline{\varphi}\,\overline{u}=
\sum_{(u)_1}\overline{u_{(1)}}\overline{\deru_{S^{-1}u_{(0)}}(\varphi)}
\qquad(u\in U,\,\,\varphi\in A).
\end{equation}
By a similar calculation we see that \eqref{eq:tU10} implies \eqref{eq:tU3}.
Hence we can replace \eqref{eq:tU10} with \eqref{eq:tU3} in defining $\tU$.
Rewriting \eqref{eq:tU3} and \eqref{eq:tU10} in terms of generators of $U$ we obtain
\begin{align}
&\overline{k_\lambda}\,\overline{\varphi}
=\overline{\deru_{k_\lambda}(\varphi)}\,\overline{k_\lambda}
\qquad
&(\lambda\in\Lambda,\,\,\varphi\in A),
\label{eq:tU4}\\
&\overline{e_i}\,\overline{\varphi}
=\overline{\deru_{k_i}(\varphi)}\,\overline{e_i}
+\overline{\deru_{e_i}(\varphi)}
\qquad
&(i\in I,\,\,\varphi\in A),
\label{eq:tU5}\\
&\overline{f_i}\,\overline{\varphi}
=\overline{\varphi}\,\overline{f_i}
+\overline{\deru_{f_i}(\varphi)}\,\overline{k_i^{-1}}
\qquad
&(i\in I,\,\,\varphi\in A),
\label{eq:tU6}\\
&\overline{\varphi}\,\overline{k_\lambda}
=\overline{k_\lambda}\,\overline{\deru_{k_\lambda^{-1}}(\varphi)}
\qquad
&(\lambda\in\Lambda,\,\,\varphi\in A),
\label{eq:tU7}\\
&\overline{\varphi}\,\overline{e_i}
=\overline{e_i}\,\overline{\deru_{k_i^{-1}}(\varphi)}
-\overline{\deru_{e_ik_i^{-1}}(\varphi)}
\qquad
&(i\in I,\,\,\varphi\in A),
\label{eq:tU8}\\
&\overline{\varphi}\,\overline{f_i}
=\overline{f_i}\,\overline{\varphi}
-\overline{k_i^{-1}}\,\overline{\deru_{k_if_i}(\varphi)}
\qquad
&(i\in I,\,\,\varphi\in A).
\label{eq:tU9}
\end{align}
\begin{proposition}
\label{prop:tU1}
The linear maps
\[
i_1:A\otimes U\to\tU\quad
(\varphi\otimes u\mapsto\overline{\varphi}\,\overline{u}),
\qquad
i_2:U\otimes A\to\tU\quad
(u\otimes\varphi\mapsto\overline{u}\,\overline{\varphi})
\]
are bijective.
\end{proposition}
\begin{proof}
We can define an $\BF$-algebra structure on $A\otimes U$ by
\[
(\varphi\otimes u)(\varphi'\otimes u')
=
\sum_{(u)_1}
\varphi\deru_{u_{(0)}}(\varphi')\otimes u_{(1)}u'.
\]
Then an algebra homomorphism $j:\tU\to A\otimes U$ is defined by
\[
j(\overline{\varphi})=\varphi\otimes1\quad(\varphi\in A),\qquad
j(\overline{u})=1\otimes u\quad(u\in U).
\]
Moreover,  $i_1$ is an algebra homomorphism, and we have $j\circ i_1=\id, i_1\circ j=\id$.
Thus $i_1$ is an isomorphism.
Similarly, $i_2$ is an isomorphism.
\end{proof}

We shall regard $A$ and $U$ as subalgebras of $\tU$ by the embeddings
\[
A\hookrightarrow\tU\,\,(\varphi\mapsto\overline{\varphi}),\qquad
U\hookrightarrow\tU\,\,(u\mapsto\overline{u}),
\]
and we sometimes write $\varphi$ and $u$ for $\overline{\varphi}$ and 
$\overline{u}$.

Note that $\tU$ is naturally a $\Lambda$-graded $\BF$-algebra by
\[
\tU(\lambda)=A(\lambda)U=UA(\lambda)\qquad
(\lambda\in\Lambda),
\]
and $A$ is an object of $\Mod_\Lambda(\tU)$.

\begin{proposition}
\label{prop:OretU}
Let $w\in W$.
\begin{itemize}
\item[{\rm(i)}]
$S_w$ satisfies the left and right Ore conditions in $\tU$.
\item[{\rm(ii)}]
The canonical homomorphism $\tU\to S_w^{-1}\tU$ is injective.
\end{itemize}
\end{proposition}
\begin{proof}
It is sufficient to show the following.
\begin{itemize}
\item[(a)]
For any $d\in \tU$ and $s\in S_w$ there exists some $t\in S_w$ and $d'\in \tU$ satisfying $td=d's$.
\item[(b)]
For any $d\in \tU$ and $s\in S_w$ there exists some $t\in S_w$ and $d'\in \tU$ satisfying $dt=sd'$.
\item[(c)]
If $sd=0$ for $d\in \tU$ and $s\in S_w$, then we have $d=0$.
\item[(d)]
If $ds=0$ for $d\in \tU$ and $s\in S_w$, then we have $d=0$.
\end{itemize}
In proving (a) and (b) we only need to deal with the cases $d=\varphi\in A$ and $d=u\in U$.
The case $d=\varphi\in A$ is already known since $S_w$ satisfies the left and right Ore conditions in $A$.
The case $d=u\in U$ is a consequence of \eqref{eq:tU4},\dots, \eqref{eq:tU9} and the case $d=\varphi\in A$.
The statements (c) and (d) follow from Proposition \ref{prop:tU1} since $S_w$ satisfies the left and right Ore conditions in $A$ and since $A\to S_w^{-1}A$ is injective.
\end{proof}
By Proposition \ref{prop:OretU} (i) we have
\begin{equation}
S_w^{-1}A\otimes_A\tU\simeq 
S_w^{-1}\tU\simeq 
\tU \otimes_A S_w^{-1}A.
\end{equation}
Moreover, for any $M\in\Mod_\Lambda(\tU)$ we have
\begin{equation}
S_w^{-1}M=
S_w^{-1}A\otimes_A M=
S_w^{-1}\tU\otimes_\tU M
\in\Mod_\Lambda(S_w^{-1}\tU).
\end{equation}
In particular, $S_w^{-1}A$ is a $U$-module.
We write the action of $U$ on $S_w^{-1}A$ by
\begin{equation}
\label{eq:s0}
U\otimes S_w^{-1}A\to S_w^{-1}A\qquad
(u\otimes \varphi\mapsto\deru_u(\varphi)).
\end{equation}
\begin{lemma}
\label{lem:UonSA}
For any $w\in W$ we have
\begin{equation}
\label{eq:UonSA}
u\varphi=\sum_{(u)_1}\deru_{u_{(0)}}(\varphi){u_{(1)}}
\qquad(u\in U,\,\,\varphi\in S_w^{-1}A)
\end{equation}
in $S_w^{-1}\tU$.
\end{lemma}
\begin{proof}
Set $\CX=\{k_\lambda,\,\,e_i\,\,,f_i\mid\lambda\in\Lambda,\,\,i\in I\}$.
For $u\in U$ and $s\in S_w$ we have $us=\sum_{(u)_1}\deru_{u_{(0)}}(s)u_{(1)}$ in $S_w^{-1}\tU$, and hence 
\[
\epsilon(u)1=\deru_u(1)=\deru_u(ss^{-1})=(us)(s^{-1})=
\sum_{(u)_1}\deru_{u_{(0)}}(s)\deru_{u_{(1)}}(s^{-1}).
\]
Considering the case $u\in\CX$ we obtain
\begin{align*}
&\deru_{k_\lambda}(s)\deru_{k_\lambda}(s^{-1})=1
&\qquad(\lambda\in\Lambda),\\
&\deru_{e_i}(s)s^{-1}+\deru_{k_i}(s)\deru_{e_i}(s^{-1})=0
&\qquad(i\in I),\\
&\deru_{f_i}(s)\deru_{k^{-1}_i}(s^{-1})+s\deru_{f_i}(s^{-1})=0
&\qquad(i\in I).
\end{align*}
We can rewrite them as
\begin{align*}
&\deru_{k_\lambda}(s^{-1})\deru_{k_\lambda}(s)=1
&\qquad(\lambda\in\Lambda),\\
&\deru_{e_i}(s^{-1})s+\deru_{k_i}(s^{-1})\deru_{e_i}(s)=0
&\qquad(i\in I),\\
&\deru_{f_i}(s^{-1})\deru_{k^{-1}_i}(s)+s^{-1}\deru_{f_i}(s)=0
&\qquad(i\in I),
\end{align*}
which gives
\begin{equation}
\label{eq:s1}
\sum_{(u)_1}\deru_{u_{(0)}}(s^{-1})\deru_{u_{(1)}}(s)
=\epsilon(u)1\qquad(u\in\CX).
\end{equation}
Hence for $u\in\CX$ and $s\in S_w$ we have
\[
\left(\sum_{(u)_1}\deru_{u_{(0)}}(s^{-1})u_{(1)}\right)s
=\sum_{(u)_2}\deru_{u_{(0)}}(s^{-1})\deru_{u_{(1)}}(s)u_{(2)}
=\sum_{(u)_1}\epsilon(u_{(0)})u_{(1)}=u
\]
in $S_w^{-1}\tU$.
Here, we  have used \eqref{eq:s1} and the fact that for any $u\in\CX$ $\Delta(u)$ is a sum of the elements of $\CX\otimes\CX$.
Hence we have 
\[
us^{-1}=\sum_{(u)_1}\deru_{u_{(0)}}(s^{-1})u_{(1)}
\qquad(u\in\CX).
\]
Let $\varphi\in S_w^{-1}A$.
We can write it as $\varphi=s^{-1}\psi$ for $s\in S_w,\,\,\psi\in A$.
Then we obtain
\begin{equation}
\label{eq:s2}
u\varphi=us^{-1}\psi=
\sum_{(u)_1}\deru_{u_{(0)}}(s^{-1})u_{(1)}\psi
=\sum_{(u)_2}\deru_{u_{(0)}}(s^{-1})\deru_{u_{(1)}}(\psi)u_{(2)}
\end{equation}
for any $u\in\CX$.
In particular, we have
\begin{equation}
\label{eq:s3}
\deru_u(\varphi)=(u\varphi)(1)
=\sum_{(u)_2}\deru_{u_{(0)}}(s^{-1})\deru_{u_{(1)}}(\psi)\deru_{u_{(2)}}(1)
=\sum_{(u)_1}\deru_{u_{(0)}}(s^{-1})\deru_{u_{(1)}}(\psi)
\end{equation}
for any $u\in\CX$.
By \eqref{eq:s2}, \eqref{eq:s3} we obtain \eqref{eq:UonSA} for $u\in\CX$. 
Since $\CX$ generates $U$, \eqref{eq:UonSA} holds for any $u\in U$.
\end{proof}
\begin{remark}
{\rm
We can give a more conceptual proof of Lemma \ref{lem:UonSA} using the right $A$-module structures on $A$ and $S_w^{-1}A$ given by the right multiplications.
}
\end{remark}

\begin{proposition}
\label{prop:U-action}
Let $w\in W$.
The $U$-module structure on $S_w^{-1}A$ is characterized by the following conditions:
\begin{itemize}
\item[\rm(a)]
The inclusion $A\hookrightarrow S_w^{-1}A$ is a homomorphism of $U$-modules.
\item[\rm(b)]
For $u\in U, \varphi_0, \varphi_1\in S_w^{-1}A$ we have 
$\deru_u(\varphi_0\varphi_1)=\sum_{(u)_1}\deru_{u_{(0)}}(\varphi_0)\deru_{u_{(1)}}(\varphi_1)$.
\end{itemize}
\end{proposition}
\begin{proof}
The $U$-module structure \eqref{eq:s0} obviously satisfies the condition (a), and it satisfies (b) by Lemma \ref{lem:UonSA}.

Assume that we are given a $U$-module structure 
\[
U\otimes S_w^{-1}A\to S_w^{-1}A\qquad
(u\otimes \varphi\mapsto\deru_u(\varphi))
\]
satisfying (a) and (b).
For $s\in S_w$ and $u\in U$ we have
\[
\epsilon(u)1=\deru_u(1)=\deru_u(s^{-1}s)
=\sum_{(u)_1}\deru_{u_{(0)}}(s^{-1})\deru_{u_{(1)}}(s)
\]
by (a) and (b). 
In particular, we obtain
\begin{align*}
&\deru_{k_\lambda}(s^{-1})\deru_{k_\lambda}(s)=1
&\qquad(\lambda\in\Lambda),\\
&\deru_{e_i}(s^{-1})s+\deru_{k_i}(s^{-1})\deru_{e_i}(s)=0
&\qquad(i\in I),\\
&\deru_{f_i}(s^{-1})\deru_{k^{-1}_i}(s)+s^{-1}\deru_{f_i}(s)=0
&\qquad(i\in I).
\end{align*}
Hence we have
\begin{align*}
&\deru_{k_\lambda}(s^{-1})=(\deru_{k_\lambda}(s))^{-1}
&\qquad(\lambda\in\Lambda),\\
&\deru_{e_i}(s^{-1})=-(\deru_{k_i}(s))^{-1}\deru_{e_i}(s)s^{-1}
&\qquad(i\in I),\\
&\deru_{f_i}(s^{-1})
=-s^{-1}\deru_{f_i}(s)(\deru_{k^{-1}_i}(s))^{-1}
&\qquad(i\in I).
\end{align*}
By (b) we have
\begin{align*}
&\deru_{k_\lambda}(s^{-1}\psi)=
(\deru_{k_\lambda}(s))^{-1}\deru_{k_\lambda}(\psi),\\
&\deru_{e_i}(s^{-1}\psi)=
-(\deru_{k_i}(s))^{-1}\deru_{e_i}(s)s^{-1}\psi
+(\deru_{k_i}(s))^{-1}\deru_{e_i}(\psi),\\
&\deru_{f_i}(s^{-1}\psi)
=-s^{-1}\deru_{f_i}(s)(\deru_{k^{-1}_i}(s))^{-1}\deru_{k_i^{-1}}(\psi)
+s^{-1}\deru_{f_i}(\psi).
\end{align*}
for any $\psi\in A$.
It implies that the action of the generators of $U$ on $S_w^{-1}A$ is uniquely determined from the $U$-module structure of $A$.
The uniqueness of the $U$-module structure on $S_w^{-1}A$ satisfying (a) and (b) is verified.
\end{proof}

\subsection{Global section functors}
\label{subsec:GS}
In this subsection $C$ denotes a $\Lambda$-graded $\BF$-algebra satisfying 
\begin{align}
&\text{$C$ contains $A$ as a $\Lambda$-graded subalgebra,}
\label{eq:condC1}\\
&\text{$S_w$ satisfies the left and right Ore conditions in $C$ for any }
\label{eq:condC2}\\
&\text{$w\in W$,}
\nonumber\\
&\text{$A(\lambda)C(\mu)\subset C(\mu)A(\lambda)$ for any $\lambda, \mu\in\Lambda$.}
\label{eq:condC3}
\end{align}
Note that $A$ and $\tU$ satisfy the conditions \eqref{eq:condC1}, \eqref{eq:condC2}, \eqref{eq:condC3}.
Another example will be the $\Lambda$-graded $\BF$-algebra $D$, which will be introduced in Section \ref{subsec:QD}.

We have an obvious exact functor
\begin{equation}
F_*:\Mod_\Lambda(C)\to \Mod_\Lambda(A)
\end{equation}
given by restricting the action of $C$ to $A$.
Its left adjoint functor is given by
\begin{equation}
F^*:\Mod_\Lambda(A)\to\Mod_\Lambda(C)\qquad(M\mapsto C\otimes_AM).
\end{equation}

The condition \eqref{eq:condC3} implies the following.
\begin{lemma}
Let $M\in\Mod_\Lambda(C)$. 
Then $\Tor(M)$, which is a priori an object of $\Mod_\Lambda(A)$, is in fact an object of $\Mod_\Lambda(C)$.
\end{lemma}

Set 
\begin{equation}
\CM(C):=\frac{\Mod_\Lambda(C)}{\Tor_\Lambda(C)}
=\Sigma_C^{-1}\Mod_\Lambda(C),
\end{equation}
where $\Tor_\Lambda(C)$ denotes the full subcategory of $\Mod_\Lambda(C)$ consisting of $M\in\Mod_\Lambda(C)$ with $F_*M\in\Tor_\Lambda(A)$, and $\Sigma_C$ is the collection of morphisms $f$ in $\Mod_\Lambda(C)$ satisfying $F_*f\in\Sigma$.
Let 
\begin{equation}
\omega_C^*:\Mod_\Lambda(C)\to\CM(C)
\end{equation}
be the localization functor.
Similarly to the case $C=A$, the abelian category $\CM(C)$ has enough injectives and we have a left exact functor 
\begin{equation}
\omega_{C*}:\CM(C)\to\Mod_\Lambda(C)
\end{equation}
which is right adjoint to $\omega_C^*$.
The canonical morphism $\omega_C^*\circ\omega_{C*}\to\Id$ is an isomorphism and we have the following.
\begin{lemma}
\label{lem:ch-of-omega-omega2}
Let $M\in\Mod_\Lambda(C)$.
Set $N=\omega_{C*}\omega_C^*M$ and let $f:M\to N$ be the canonical morphism.
Then $N$ and $f$ are uniquely characterized by the following properties.
\begin{itemize}
\item[{\rm(a)}]
$\Ker(f)$ and $\Coker(f)$ belong to $\Tor_\Lambda(C)$.
\item[{\rm(b)}]
$\Tor(N)=\{0\}$.
\item[{\rm(c)}]
Any monomorphism $N\to L$ with $L/N\in\Tor_\Lambda(C)$ is a split morphism.
\end{itemize}
\end{lemma}
We define a left exact functor
\begin{equation}
\Gamma_C:\CM(C)\to\Mod(C(0))
\end{equation}
by $\Gamma_C(M)=(\omega_{C*}M)(0)$.

By the universal property of the localization functor there exists uniquely an exact functor 
\begin{equation}
\tilde{F}_*:\CM(C)\to\CM(A)
\end{equation}
such that the following diagram is commutative:
\[
\begin{CD}
\Mod_\Lambda(C)@>{F_*}>>\Mod_\Lambda(A)
\\
@V{\omega_C^*}VV@VV{\omega^*}V
\\
\CM(C)@>>{\tilde{F}_*}>\CM(A).
\end{CD}
\]
\begin{lemma}
There exists uniquely an additive functor 
\begin{equation}
\tilde{F}^*:\CM(A)\to\CM(C)
\end{equation}
such that the following diagram is commutative:
\[
\begin{CD}
\Mod_\Lambda(A)@>{F^*}>>\Mod_\Lambda(C)
\\
@V{\omega^*}VV@VV{\omega_C^*}V
\\
\CM(A)@>>{\tilde{F}^*}>\CM(C).
\end{CD}
\]
\end{lemma}
\begin{proof}
By the universal property of the localization functor it is sufficient to show 
$F^*(\Sigma)\subset\Sigma_C$.
Let $f:M\to N$ be a morphism in $\Sigma$.
The corresponding homomorphism $S_w^{-1}A\otimes_AM\to S_w^{-1}A\otimes_AN$ is an isomorphism for any $w\in W$.
By Proposition \ref{prop:quasi-scheme} we have only to show that the corresponding morphism
$S_w^{-1}A\otimes_AC\otimes_AM\to S_w^{-1}A\otimes_AC\otimes_AN$ is an isomorphism for any $w\in W$.
This follows from $S_w^{-1}C\simeq S_w^{-1}A\otimes_AC\simeq C\otimes_AS_w^{-1}A$.
\end{proof}
\begin{lemma}
\label{lem:compatiGamma}
The following diagram is commutative:
\[
\begin{CD}
\CM(C)@>{\tilde{F}_*}>>\CM(A)
\\
@V\omega_{C*}VV@VV\omega_{*}V
\\
\Mod_\Lambda(C)@>>F_*>\Mod_\Lambda(A).
\end{CD}
\]
\end{lemma}
\begin{proof}
We have a sequence of morphisms
\[
F_*\circ\omega_{C*}
\to
\omega_*\circ\omega^*\circ F_*\circ\omega_{C*}
=
\omega_*\circ\tilde{F}_*\circ\omega_C^*\circ\omega_{C*}
=
\omega_*\circ\tilde{F}_*.
\]
Hence it is sufficient to show that $g:F_*\omega_{C*}M\to\omega_*\omega^*F_*\omega_{C*}M$ is an isomorphism for any $M$.
By 
\[
\omega^*\circ F_*\circ\omega_{C*}
=\tilde{F}_*\circ\omega_C^*\circ\omega_{C*}
=\tilde{F}_*
=\omega^*\circ\omega_*\circ\tilde{F}_*
\]
$\omega^*g$ is an isomorphism, and hence 
$\Ker(g), \Coker(g)\in\Tor_\Lambda(A)$.

By $\Tor(F_*\omega_{C*}M)=0$ and $\Ker(g)\in\Tor_\Lambda(A)$ $g$ is a monomorphism.
Hence it is sufficient to show that the injective homomorphism
\[
g_L:\Hom(L,F_*\omega_{C*}M)\to\Hom(L,\omega_*\omega^*F_*\omega_{C*}M)
\]
is surjective for any $L$.
By 
\begin{gather*}
\Hom(L,\omega_*\omega^*F_*\omega_{C*}M)
\simeq
\Hom(\omega^*L,\omega^*F_*\omega_{C*}M),
\\
\Hom(L,F_*\omega_{C*}M)
\simeq
\Hom(\omega_C^*F^*L,M)
\simeq
\Hom(\tilde{F}^*\omega^*L,M)
\end{gather*} 
$\Hom(L,F_*\omega_{C*}M)$ and 
$\Hom(L,\omega_*\omega^*F_*\omega_{C*}M)$ 
depend only on $\omega^*L$
Let $a\in\Hom(L,\omega_*\omega^*F_*\omega_{C*}M)$.
Set $L_1=\Ker(L\to\Coker(g))$.
Since $\Coker(g)$ belongs to $\Tor_\Lambda(A)$ we have 
$L/L_1\in\Tor_\Lambda(A)$, and hence we have 
$\omega^*L_1\simeq\omega^*L$.
Therefore, we can replace $L$ with $L_1$.
Then $a|L_1$ is clearly contained in the image of $g_{L_1}$.
\end{proof}
In view of Lemma \ref{lem:compatiGamma} we shall often drop the subscript $C$ in $\omega_{C*}$ and $\Gamma_C$ and write them simply as $\omega_*$ and $\Gamma$.

Let $K$ be a full subcategory of $\Mod_\Lambda(C)$ closed under taking subquotients in $\Mod_\Lambda(C)$, and set
\[
\CK=
{K}/
{\Tor_\Lambda(C)\cap K}
=
\overline{\Sigma}^{-1}K,
\]
where $\overline{\Sigma}$ is the collection of morphisms in 
$K$ whose kernel and cokernel belong to 
$\Tor_\Lambda(C)\cap K$,
and denote by
$\overline{\omega}^*:K\to\CK$
the localization functor.
Let $j:K\to\Mod_{\Lambda}(C)$ be the embedding.
By the universal property of the localization functor we have a functor  
$\overline{j}:\CK\to\CM(C)$
such that the following diagram commutes:
\[
\begin{CD}
K@>{j}>>\Mod_\Lambda(C)
\\
@V{\overline{\omega}^*}VV@VV{\omega_C^*}V
\\
\CK@>>{\overline{j}}>\CM(C).
\end{CD}
\]
\begin{lemma}
\label{lem:FF}
$\overline{j}$ is fully faithful.
\end{lemma}
\begin{proof}
We need to  show that the canonical homomorphism
\begin{equation}
\label{eq:FF}
\Hom_{\CK}(\overline{\omega}^*M,\overline{\omega}^*N)
\to
\Hom_{\CM(C)}(\omega_C^*M,\omega_C^*N)
\end{equation}
is an isomorphism for $M, N\in K$.
We may assume $\Tor(M)=\Tor(N)=\{0\}$.

Assume that $\overline{\omega}^*f\circ(\overline{\omega}^*s)^{-1}$ belongs to the kernel of \eqref{eq:FF}, where $f:L\to N$ is a morphism in 
$K$ and $s:L\to M$ belongs to $\overline{\Sigma}$.
By $\Tor(M)=\Tor(N)=\{0\}$ we have $\Tor(L)=\Ker(s)\subset\Ker(f)$, and hence we may assume that $\Tor(L)=\{0\}$ by replacing $L$ with $L/\Tor(L)$.
Since $\overline{\omega}^*f\circ(\overline{\omega}^*s)^{-1}$ belongs to the kernel of \eqref{eq:FF}, we have $\omega_C^*f\circ(\omega_C^*s)^{-1}=0$, and hence $\omega_C^*f=0$.
It means that there exists $t:R\to L$, which belongs to $\Sigma_C$, such that $f\circ t=0$.
By $\Tor(L)=\{0\}$ we have $\Tor(R)=\Ker(t)$.
Hence we can assume $\Ker(t)=\{0\}$ by replacing $R$ with $R/\Tor(R)$.
Then $R$ is a subobject of $L$, and hence $t$ belongs to $\overline{\Sigma}$.
It follows that $\overline{\omega}^*f=0$, and hence 
$\overline{\omega}^*f\circ(\overline{\omega}^*s)^{-1}=0$.

Take a morphism $\omega_C^*f\circ(\omega_C^*s)^{-1}$ in $\CM(C)$, where $f:L\to N$ is a morphism in 
$\Mod_{\Lambda}(C)$ and $s:L\to M$ belongs to $\Sigma_{C}$.
By $\Tor(M)=\Tor(N)=\{0\}$ we have $\Tor(L)=\Ker(s)\subset\Ker(f)$, and hence we may assume that $\Ker(s)=\{0\}$ by replacing $L$ with $L/\Tor(L)$.
Then $L$ is a subobject of $M$ and hence belongs to $K$.
Hence $\omega_C^*f\circ(\omega_C^*s)^{-1}$ is in the image of $\eqref{eq:FF}$.
\end{proof}

\subsection{The vector bundle $E^\mu$ associated to $V(\mu)$}
\label{subsec:Emu}
The contents of Sections \ref{subsec:Emu}, \ref{subsec:filt} except for Proposition \ref{prop:filt-E1} below are due to Lunts-Rosenberg \cite{LR}.
Some of the proofs are also included for the convenience of the readers.

Let $\mu\in\Lambda^+$.
Following Lunts-Rosenberg \cite{LR} we shall define an $A$-bimodule $E^\mu$, which is a $q$-analogue of the vector bundle $\CO_\CB\otimes_\BC V^1(\mu)$.

Set
\begin{equation}
E^\mu=V(\mu)\otimes A.
\end{equation}
We endow $E^\mu$ with a right $A$-module structure by $(v\otimes\varphi)\psi=v\otimes\varphi \psi$ for $v\in V(\mu)$ and $\varphi, \psi\in A$.
We identify $E^\mu$ with $A\otimes V(\mu)$ via the isomorphism
\[
\eta=\CR^\vee_{A,V(\mu)}:A\otimes V(\mu)\to V(\mu)\otimes A=E^\mu,
\]
and define a left $A$-module structure on $E^\mu=A\otimes V(\mu)$ by 
$\psi(\varphi\otimes v)=\psi\varphi\otimes v$ for $v\in V(\mu)$ and $\varphi, \psi\in A$.
Note that $\eta$ is well-defined since $A$ is a sum of $U$-submodules belonging to $\Mod^f(U)$.
\begin{lemma}
\label{lem:Edef}
\begin{itemize}
\item[\rm(i)]
For $\varphi, \psi\in A$ and $e\in E^\mu$ we have $(\varphi e)\psi=\varphi(e\psi)$.
\item[\rm(ii)]
We have a commutative diagram
\[
\begin{CD}
V(\mu) @= V(\mu)\\
@VVV @VVV\\
A\otimes V(\mu) @>{\eta}>> V(\mu)\otimes A,
\end{CD}
\]
where the left $($right$)$ vertical arrow is given by $v\mapsto 1\otimes v$ $($resp.\ $v\mapsto v\otimes 1$$)$.
\end{itemize}
\end{lemma}
\begin{proof}
(i) Let $m:A\otimes A\to A$ be the multiplication of the algebra $A$.
It is sufficient to show that the composition of
\[
A\otimes V(\mu)\otimes A
\xrightarrow{\eta\otimes\id_A}
V(\mu)\otimes A\otimes A
\xrightarrow{\id_{V(\mu)}\otimes m}
V(\mu)\otimes A
\]
gives the left action of $A$ on $E^\mu=V(\mu)\otimes A$.
This is equivalent to showing that 
\[
(\id_{V(\mu)}\otimes m)\circ
(\eta\otimes\id_A)\circ
(\id_A\otimes\eta):
A\otimes A\otimes V(\mu)\to V(\mu)\otimes A
\]
coincides with 
\[
\eta\circ
(m\otimes\id_{V(\mu)}):
A\otimes A\otimes V(\mu)\to V(\mu)\otimes A.
\]
Indeed we have
\begin{align*}
&(\id_{V(\mu)}\otimes m)\circ(\eta\otimes\id_A)\circ(\id_A\otimes\eta)\\
=&(\id_{V(\mu)}\otimes m)\circ(\CR^\vee_{A,V(\mu)}\otimes\id_A)\circ(\id_A\otimes\CR^\vee_{A,V(\mu)})\\
=&(\id_{V(\mu)}\otimes m)\circ\CR^\vee_{A\otimes A,V(\mu)}\\
=&\CR^\vee_{A,V(\mu)}\circ(m\otimes\id_{V(\mu)})\\
=&\eta\circ(m\otimes\id_{V(\mu)}).
\end{align*}
Here the second equality is a consequence of Proposition \ref{prop:R} (iii), and the third equality follows from \eqref{eq:R}.

(ii) Note that $A(0)=\BF 1$ is isomorphic to the trivial $U$-modules $V(0)$.
We need to show that $\CR^\vee_{V(0),V(\mu)}:V(0)\otimes V(\mu)\to V(\mu)\otimes V(0)$ is equal to $\id_{V(\mu)}$ under the identification $V(0)\otimes V(\mu)\simeq V(\mu)\otimes V(0)\simeq V(\mu)$ of $U$-modules.
This follows from the definition of $\CR^\vee_{V(0),V(\mu)}$.
\end{proof}
By Lemma \ref{lem:Edef} (i) $E^\mu$ is an $A$-bimodule.
By Lemma \ref{lem:Edef} (ii) we have a canonical embedding 
$V(\mu)\hookrightarrow E^\mu$ 
such that the actions of $A$ on $E^\mu$ from the left and the right induce 
$A\otimes V(\mu)\simeq E^\mu$ and $V(\mu)\otimes A\simeq E^\mu$ respectively.

For $\lambda\in\Lambda$ we set
\begin{equation}
E^\mu(\lambda)=V(\mu)\otimes A(\lambda)\subset E^\mu.
\end{equation}
By $\eta(A(\lambda)\otimes V(\mu))=V(\mu)\otimes A(\lambda)$ 
we have $E^\mu(\lambda)=A(\lambda)\otimes V(\mu)$ under the identification $E^\mu=A\otimes V(\mu)$.
Moreover, we have
\[
E^\mu=\bigoplus_{\lambda\in\Lambda^+}E^\mu(\lambda),\quad
A(\xi)E^\mu(\lambda)\subset E^\mu(\xi+\lambda),\quad
E^\mu(\lambda)A(\xi)\subset E^\mu(\xi+\lambda).
\]

Let $M\in\Mod_\Lambda(A)$.
We have a natural left $A$-module structure on $E^\mu\otimes_AM\simeq V(\mu)\otimes M$ induced from the left action of $A$ on $E^\mu$.
Moreover, we have $E^\mu\otimes_AM\in\Mod_\Lambda(A)$ by 
$(E^\mu\otimes_AM)(\lambda)=V(\mu)\otimes M(\lambda)$, and we obtain an exact functor
\begin{equation}
\label{eq:Emu1}
E^\mu\otimes_A(\bullet):\Mod_\Lambda(A)\to\Mod_\Lambda(A)
\end{equation}
sending $M$ to $E^\mu\otimes_A M$.
\begin{lemma}
The functor \eqref{eq:Emu1} induces
\[
E^\mu\otimes_A(\bullet):\CM(A)\to\CM(A).
\]
\end{lemma}
\begin{proof}
It is sufficient to show that for any $M\in\Tor_\Lambda(A)$ we have $E^\mu\otimes M\in\Tor_\Lambda(A)$.
This follows from $A(\xi)(V(\mu)\otimes m)\subset V(\mu)\otimes A(\xi)m$ for any $m\in M$.
\end{proof}
\begin{proposition}
\label{prop:filt-E1} 
For any $\overline{M}\in\CM(A)$ we have
\[
E^\mu\otimes_A\omega_*\overline{M}
\simeq\omega_*(E^\mu\otimes_A\overline{M}).
\]
\end{proposition}
\begin{proof}
Choose a filtration
\[
V(\mu)=V^n\supset V^{n-1}\supset\cdots\supset V^1\supset V^0=\{0\}
\]
of $V(\mu)$ consisting of $U^{\leqq0}$-submodules $V^k$ satisfying 
$\dim V^k/V^{k-1}=1$, and consider the corresponding filtration
\[
E^\mu=E^n\supset E^{n-1}\supset\cdots\supset E^1\supset E^0=\{0\}
\]
of the right $A$-module $E^\mu=V(\mu)\otimes A$ given by $E^k=V^k\otimes A$.
By the definition of the left $A$-module structure on $E^\mu$,  especially by the fact that $\Xi$ belongs to a completion of $U\otimes U^{\leqq0}$, we see that $E^k$ is a left $A$-submodule for any $k$.
Let $\nu_k\in\Lambda$ be the weight of $V^k/V^{k-1}$ and take $v_k\in V(\mu)_{\nu_k}$ such that $V^k=\BF v_k\oplus V^{k-1}$.
Let $\overline{v}_k$ be the corresponding element of the $A$-bimodule $E^k/E^{k-1}$.
For any $\varphi\in A(\lambda)_\xi$ we have $\varphi\overline{v}_k=q^{-(\nu_k,\xi)}\overline{v}_k\varphi$ by the explicit form of $\Xi$.

For $\nu\in\Lambda$ we define an automorphism $h_\nu$ of the graded $\BF$-algebra $A$ by $h_\nu(\varphi)=q^{-(\nu,\xi)}\varphi$ for $\varphi\in A(\lambda)_\xi$.
For $N\in\Mod_\Lambda(A)$ we define $h_\nu^\bullet N\in\Mod_\Lambda(A)$ by
\begin{align*}
&h_\nu^\bullet N\simeq N\quad(h_\nu^\bullet(n)\leftrightarrow n)
\qquad\text{as grade $\BF$-modules}\\
&\varphi(h_\nu^\bullet(n))=h_\nu^\bullet((h_\nu(\varphi))n)
\qquad(\varphi\in A,\,\, n\in N).
\end{align*}
Then we have $E^k/E^{k-1}\otimes_AN\simeq h_{\nu_k}^\bullet N$ as a graded $A$-module.
Note that $h_\nu^\bullet$ induces exact functors
\[
h_\nu^\bullet:\Mod_\Lambda(A)\to\Mod_\Lambda(A),\qquad
\tilde{h}_\nu^\bullet:\CM(A)\to\CM(A)
\]
satisfying $\omega_*\tilde{h}_\nu^\bullet=h_\nu^\bullet\omega_*$.

Let $\overline{M}\in\CM(A)$.
We have morphisms
\begin{align*}
&\Phi^k:E^k\otimes_A\omega_*\overline{M}
\to\omega_*(E^k\otimes_A\overline{M}),\\
&\Psi^k:E^k/E^{k-1}\otimes_A\omega_*\overline{M}
\to\omega_*(E^k/E^{k-1}\otimes_A\overline{M})
\end{align*}
in $\Mod_\Lambda(A)$ functorial with respect to $\overline{M}$ such that
\begin{align*}
&\omega^*\Phi^k:
\omega^*(E^k\otimes_A\omega_*\overline{M})
\to\omega^*\omega_*(E^k\otimes_A\overline{M}),\\
&\omega^*\Psi^k:
\omega^*(E^k/E^{k-1}\otimes_A\omega_*\overline{M})
\to\omega^*\omega_*(E^k/E^{k-1}\otimes_A\overline{M})
\end{align*}
are isomorphisms.
Note that $E^k\otimes_A(\bullet)$ and $E^k/E^{k-1}\otimes_A(\bullet)$ on $\CM(A)$ are defined similarly to $E^\mu\otimes_A(\bullet)$ on $\CM(A)$.

Note that $\Psi^k$ is an isomorphism by 
\[
E^k/E^{k-1}\otimes_A\omega_*\overline{M}
\simeq
h_{\nu_k}^\bullet\omega_*\overline{M}
\simeq
\omega_*\tilde{h}_\nu^\bullet \overline{M}
\simeq
\omega_*(E^k/E^{k-1}\otimes_A\overline{M}).
\]

Let us show that $\Phi^k$ is an isomorphism.
The surjectivity is proved by induction on $k$ using the following commutative diagram whose rows are exact.
{\scriptsize
\[
\begin{CD}
0 @>>> 
E^{k-1}\otimes_A\omega_*\overline{M} @>>> 
E^{k}\otimes_A\omega_*\overline{M}
@>>>E^{k}/E^{k-1}\otimes_A\omega_*\overline{M} 
@>>> 0
\\
@. @V{\Phi^{k-1}}VV @V{\Phi^k}VV @VV{\Psi^k}V@.
\\
0 @>>> \omega_*(E^{k-1}\otimes_A\overline{M}) 
@>>> \omega_*(E^{k}\otimes_A\overline{M})
@>>>\omega_*(E^k/E^{k-1}\otimes_A\overline{M})@.
\end{CD}
\]
}
Since $\omega^*\Phi^k$ is an isomorphism, $\Ker(\Phi^k)$ belongs to $\Tor_\Lambda(A)$.
Hence, in order to prove that $\Phi^k$ is injective, we have only to show $\Tor(E^k\otimes_A\omega_*\overline{M})=\{0\}$.
By $\Tor(\omega_*\overline{M})=\{0\}$ (see Corollary \ref{cor:ch-of-omega-omega}) it is sufficient to show $\Tor(E^k\otimes_AN)=\{0\}$ for any $N$ with $\Tor(N)=\{0\}$.
The image of $\Tor(E^k\otimes_AN)$ under $E^k\otimes_AN\to E^k/E^{k-1}\otimes_AN$ is contained in $\Tor(E^k/E^{k-1}\otimes_AN)\simeq\Tor(h_{\nu_k}^\bullet N)=\{0\}$, and hence $\Tor(E^k\otimes_AN)\subset\Tor(E^{k-1}\otimes_AN)=\{0\}$ by induction on $k$.

We have obtained the desired result by considering the case $k=n$.
\end{proof}

\begin{remark}{\rm
By Proposition \ref{prop:filt-E1} we obtain
\[
\Gamma(V(\mu)\otimes \overline{M})
=\Gamma(E^\mu\otimes_A\overline{M})=
E^\mu\otimes_A\Gamma(\overline{M})=
V(\mu)\otimes\Gamma(\overline{M})
\]
for any object $M$ of $\CM(A)$.
As noted in Lunts-Rosenberg \cite[IV, 6.6 Remark]{LR} this implies Theorem \ref{thm1}.
}
\end{remark}
\subsection{Filtration of $E^\mu$}
\label{subsec:filt}
Define a left $U$-module structure on $E^\mu=V(\mu)\otimes A$ by 
\[
u(v\otimes\varphi)=\sum_{(u)_1}u_{(0)}v\otimes u_{(1)}\varphi
\qquad(u\in U, v\in V(\mu), \varphi\in A).
\]
Since $\eta$ is an isomorphism of $U$-modules, we have 
\begin{equation}
u(\varphi e\psi)=
\sum_{(u)_2}(u_{(0)}\varphi)(u_{(1)}e)(u_{(2)}\psi)\qquad
(u\in U,\,\,\varphi, \psi\in A,\,\, e\in E^\mu).
\end{equation}
For $M\in\Mod_\Lambda(\tilde{U})$ a left $U$-module structure on $E^\mu\otimes_AM$ is defined by
\begin{equation}
u(e\otimes m)=\sum_{(u)_1}u_{(0)}e\otimes u_{(1)}m\qquad
(u\in U, \,\,e\in E^\mu,\,\, m\in M),
\end{equation}
and it gives a $\Lambda$-graded left $\tilde{U}$-module structure on $E^\mu\otimes_AM$.
Moreover $E^\mu\otimes_A(\bullet)$ induces an exact functor
\begin{equation}
E^\mu\otimes_A(\bullet):\CM(\tilde{U})\to\CM(\tilde{U}).
\end{equation}

We fix $\lambda_0\in\Lambda^+$ such that 
\begin{equation}
\lambda_0+\nu\in\Lambda^+ \qquad\text{ for any weight $\nu$ of $V(\mu)$}, 
\end{equation}
and set
\begin{equation}
\overline{E}^\mu=V(\mu)\otimes
(
\bigoplus_{\lambda\in\lambda_0+\Lambda^+}A(\lambda)
)
\subset V(\mu)\otimes A=E^\mu.
\end{equation}
The following is obvious from the definition.
\begin{lemma}
\label{lem:filt0}
$\overline{E}^\mu$ is a graded $(A,A)$-submodule and a $U$-submodule of $E^\mu$.
Moreover, $E^\mu/\overline{E}^\mu$ belongs to $\Tor_\Lambda(A)$ as a graded left $A$-module.
\end{lemma}

We fix a labeling $\{\nu_1,\dots,\nu_r\}$ of the set of distinct weights of $V(\mu)$ such that $\nu_i-\nu_j\in \Lambda^+$ implies $i\geqq j$.
In particular, $\nu_1$ is the lowest weight $w_0\mu$ and $\nu_r$ is the highest weight $\mu$.
Set $m_j=\dim V(\mu)_{\nu_j}$.
For $\lambda\in\lambda_0+\Lambda^+$ we have 
\[
\overline{E}^\mu(\lambda)\simeq\bigoplus_{j=1}^rV(\lambda+\nu_j)^{\oplus m_j}
\]
as $U$-modules since
\[
\ch(\overline{E}^\mu(\lambda))
=\ch(V(\mu)\otimes A(\lambda))
=\sum_{j=1}^rm_j\ch(V(\lambda+\nu_j))
\]
by Weyl's character formula.
Define a filtration 
\begin{equation}
\{0\}=\overline{E}^\mu_0\subset
\overline{E}^\mu_1\subset
\cdots\subset
\overline{E}^\mu_r=\overline{E}^\mu
\end{equation}
of $\overline{E}^\mu$ consisting of graded $U$-submodules by
\[
\overline{E}^\mu_k(\lambda)\simeq
\bigoplus_{j=1}^kV(\lambda+\nu_j)^{\oplus m_j}
\qquad(\lambda\in\lambda_0+\Lambda^+).
\]
\begin{lemma}
$\overline{E}^\mu_k$ is an $(A,A)$-submodule of $\overline{E}^\mu$.
\end{lemma} 
\begin{proof}
Let $\lambda\in\lambda_0+\Lambda^+$ and let $T$ be a $U$-submodule of $\overline{E}^\mu(\lambda)=V(\mu)\otimes A(\lambda)$ isomorphic to $V(\lambda+\nu_t)$.
Let $\xi\in\Lambda^+$.
Then $TA(\xi)$ is the image of $T\otimes A(\xi)$ under the homomorphism $\overline{E}^\mu(\lambda)\otimes A(\xi)
\ni e\otimes\varphi\mapsto e\varphi\in
\overline{E}^\mu(\lambda+\xi)$ of $U$-modules.
Hence $TA(\xi)$ is a $U$-submodule of $\overline{E}^\mu(\lambda+\xi)$ whose weights are contained in $\lambda+\xi+\nu_t-\Lambda^+$.
It follows that $TA(\xi)\subset\bigoplus_{j=1}^tV(\lambda+\xi+\nu_j)^{\oplus m_j}$, and hence $\overline{E}^\mu_k$ is a right $A$-submodule.
The assertion about the left module structure is proved similarly.
\end{proof}

For $k=1,\dots,r$ we set 
\[
A^{(k)}=\bigoplus_{\lambda\in\lambda_0+\nu_k+\Lambda^+}A(\lambda).
\]
It is a graded $(A,A)$-submodule and a $U$-submodule of $A$.
\begin{lemma}
\label{lem:filt2}
\begin{itemize}
\item[\rm(i)]
There exists an isomorphism 
\[
\Phi:\oE^\mu_k/\oE^\mu_{k-1}\to A^{(k)}[\nu_k]^{\oplus m_k}
\]
of graded right $A$-modules and $U$-modules.
\item[\rm(ii)]
Identify $A^{(k)}[\nu_k]^{\oplus m_k}$ with $A^{(k)}[\nu_k]\otimes\BF^{m_k}$.
Then there exists a group homomorphism $\tau:\Lambda\to GL_{m_k}(\BF)$ such that $\Phi(\varphi v)=(\id\otimes \tau(\xi))\varphi\Phi(v)$ for any $\varphi\in A(\xi), v\in\oE_k^\mu/\oE_{k-1}^\mu$.
\end{itemize}
\end{lemma}
\begin{proof}
For simplicity we set $M=\oE^\mu_k/\oE^\mu_{k-1}$ and $N=A^{(k)}[\nu_k]^{\oplus m_k}$.

(i) As $U$-modules we have
\[
M(\lambda)\simeq
N(\lambda)
\simeq
\begin{cases}
V(\lambda+\nu_k)^{\oplus m_k}\qquad&(\lambda\in\lambda_0+\Lambda^+),\\
0&(\text{otherwise}).
\end{cases}
\]

For $\lambda\in\lambda_0+\Lambda^+, \xi\in\Lambda^+, 
c\in A(\xi)_\xi\setminus\{0\}$ we have linear maps
\begin{align}
&M(\lambda)_{\lambda+\nu_k}\ni v\mapsto vc\in M(\lambda+\xi)_{\lambda+\xi+\nu_k},
\label{eq:Efilt1}\\
&N(\lambda)_{\lambda+\nu_k}\ni v\mapsto vc\in N(\lambda+\xi)_{\lambda+\xi+\nu_k}.
\label{eq:Efilt2}
\end{align}
Let us show that they are isomorphisms.
Considering the dimensions it is sufficient to show that they are injective.
Since $A$ is a domain, \eqref{eq:Efilt2} is injective.
Set $\tilde{M}=\oE^\mu_k$.
Then the projection $\oE^\mu_k\to\oE^\mu_k/\oE^\mu_{k-1}$ induces 
$\tilde{M}(\lambda)_{\lambda+\nu_k}\simeq M(\lambda)_{\lambda+\nu_k}$.
Hence the injectivity of \eqref{eq:Efilt1} follows from the injectivity of $\tilde{M}(\lambda)_{\lambda+\nu_k}\to \tilde{M}(\lambda+\xi)_{\lambda+\xi+\nu_k}$, which is a consequence of $\oE^\mu_k\subset V(\mu)\otimes A$.

Hence there exists a family $\beta_{\lambda}:M(\lambda)_{\lambda+\nu_k}\to N(\lambda)_{\lambda+\nu_k}$ \,\,($\lambda\in\lambda_0+\Lambda^+$) of linear isomorphisms satisfying 
$\beta_\lambda(v)c=\beta_{\lambda+\xi}(vc)$ for any $\lambda\in\lambda_0+\Lambda^+,\,\, \xi\in\Lambda^+,\,\,v\in M(\lambda)_{\lambda+\nu_k},\,\, c\in A(\xi)_\xi$.
From this we obtain an isomorphism $\Phi:M\to N$ of graded $U$-modules given by
\[
\Phi(uv)=u\beta_\lambda(v)\qquad
(u\in U,\,\,\lambda\in\lambda_0+\Lambda^+,\,\,v\in M(\lambda)_{\lambda+\nu_k}).
\]

It remains to show the commutativity of the diagram:
\[
\begin{CD}
M(\lambda)\otimes A(\xi)
@>>>M(\lambda+\xi)\\
@V{\Phi\otimes\id}VV @VV{\Phi}V\\
N(\lambda)\otimes A(\xi)
@>>>
N(\lambda+\xi)
\end{CD}
\]
for $\lambda\in\lambda_0+\Lambda^+, \,\,\xi\in\Lambda^+$.
Here the horizontal arrows are given by the right $A$-module structures.
In particular, they are homomorphisms of $U$-modules.
Therefore the assertion follows from the commutativity of 
\[
\begin{CD}
M(\lambda)_{\lambda+\nu_k}\otimes A(\xi)_\xi
@>>>M(\lambda+\xi)_{\lambda+\xi+\nu_k}\\
@V{\Phi\otimes\id}VV @VV{\Phi}V\\
N(\lambda)_{\lambda+\nu_k}\otimes A(\xi)_\xi
@>>>N(\lambda+\xi)_{\lambda+\xi+\nu_k}.
\end{CD}
\]

(ii) Similarly to the proof of (i) it is sufficient to show that there exists a group homomorphism $\tau:\Lambda\to GL_{m_k}(\BF)$ such that 
$\Phi(cv)=(\id\otimes \tau(\xi))c\Phi(v)$ for any 
$\xi\in\Lambda^+,\,\,
\lambda\in\lambda_0+\Lambda^+,\,\,
c\in A(\xi)_\xi, v\in M(\lambda)_{\lambda+\nu_k}$.

Let $\xi\in\Lambda^+$.
Considering the linear isomorphisms \eqref{eq:Efilt1}, \eqref{eq:Efilt2} for $\lambda=\lambda_0$ we obtain
$\tau(\xi)\in GL_{m_k}(\BF)$ such that $\Phi(cm)=(\id\otimes\tau(\xi))c\Phi(m)$ for any 
$c\in A(\xi)_\xi,\,\,m\in M(\lambda_0)_{\lambda_0+\nu_k}$.
Let $c\in A(\xi)_\xi,\,\,c'\in A(\eta)_\eta,\,\,m\in M(\lambda_0)_{\lambda_0+\nu_k}$.
Then we have
\begin{align*}
\Phi(c(mc'))&=\Phi((cm)c')=\Phi(cm)c'=((\id\otimes\tau(\xi))c\Phi(m))c'\\
&=(\id\otimes\tau(\xi))c(\Phi(m)c')
=(\id\otimes\tau(\xi))c(\Phi(mc')),
\end{align*}
and hence we obtain $\Phi(cm)=(\id\otimes\tau(\xi))c\Phi(m)$ for any $\xi\in\Lambda^+,\,\,\lambda\in\lambda_0+\Lambda^+,\,\,c\in A(\xi)_\xi,\,\,m\in M(\lambda)_{\lambda+\nu_k}$.
We have $\tau(0)=\id$ and $\tau(\xi)\tau(\xi')=\tau(\xi+\xi')$ for any $\xi, \xi'\in\Lambda^+$, and hence $\tau$ is extended to a group homomorphism $\tau:\Lambda\to GL_{m_k}(\BF)$.
\end{proof}

\begin{lemma}
Let $f:E_1\to E_2$ be a morphism of $\Lambda$-graded $A$-bimodules such that $\Ker(f), \Coker(f)\in\Tor_\Lambda(A)$ as graded left $A$-modules.
Let $M\in\Mod_\Lambda(A)$, and let $\tilde{f}:E_1\otimes_AM\to E_2\otimes_AM$ be the corresponding morphism in $\Mod_\Lambda(A)$.
Then $\Ker(\tilde{f}), \Coker(\tilde{f})\in\Tor_\Lambda(A)$.
\end{lemma}
\begin{proof}
By the assumption we have $S_w^{-1}E_1\simeq S_w^{-1}E_2$ for any $w\in W$.
Thus
\[
S_w^{-1}(E_1\otimes_AM)\simeq S_w^{-1}E_1\otimes_AM\simeq
S_w^{-1}E_2\otimes_AM\simeq S_w^{-1}(E_2\otimes_AM).
\]
Hence we have 
$\omega^*(E_1\otimes_AM)\simeq\omega^*(E_2\otimes_AM)$ by Proposition \ref{prop:quasi-scheme}.
This is equivalent to $\Ker(\tilde{f}), \Coker(\tilde{f})\in\Tor_\Lambda(A)$.
\end{proof}

Hence we have the following.
\begin{lemma}
\label{lem:EoE}
$E^\mu\otimes_A(\bullet)$ and $\overline{E}^\mu\otimes_A(\bullet)$ are isomorphic as functors from $\CM(A)$ to $\CM(A)$ and from $\CM(\tilde{U})$ to $\CM(\tilde{U})$.
\end{lemma}

Let $M\in\Mod_\Lambda(A)$ (resp.\ $\Mod_\Lambda(\tU)$).
By Lemma \ref{lem:EoE} we have 
$\omega^*(\overline{E}^\mu\otimes_AM)=
\omega^*(E^\mu\otimes_AM)$.
Hence the filtration 
\[
\{0\}=\overline{E}^\mu_0\subset
\overline{E}^\mu_1\subset
\cdots\subset
\overline{E}^\mu_r=\overline{E}^\mu
\subset E^\mu
\]
of $E^\mu$ induces the sequence 
\[
\{0\}=\omega^*(\overline{E}^\mu_0\otimes_AM)\to
\omega^*(\overline{E}^\mu_1\otimes_AM)\to
\cdots\to
\omega^*(\overline{E}^\mu_r\otimes_AM)=
\omega^*(E^\mu\otimes_AM)
\]
of morphisms in $\CM(A)$ (resp.\ $\CM(\tU)$).
\begin{lemma}
\label{lem:Filtration}
Let $M\in\Mod_\Lambda(A)$ $($resp.\ $\Mod_\Lambda(\tU)$$)$.
\begin{itemize}
\item[\rm(i)]
We have the exact sequence:
\[
0\to\omega^*(\overline{E}^\mu_{k-1}\otimes_AM)
\to\omega^*(\overline{E}^\mu_k\otimes_AM)
\to\omega^*((\overline{E}^\mu_k/\overline{E}^\mu_{k-1})\otimes_AM)
\to 0.
\]
in  $\CM(A)$ $($resp.\ $\CM(\tU)$$)$.
\item[\rm(ii)]
We have an isomorphism
\[
\omega^*((\overline{E}^\mu_k/\overline{E}^\mu_{k-1})\otimes_AM)
\simeq\omega^*M[\nu_k]^{\oplus m_k}.
\]
in  $\CM(A)$ $($resp.\ $\CM(\tU)$$)$.
\end{itemize}
\end{lemma}
\begin{proof}
The statements for $\tU$ easily follow from those for $A$, and hence we shall only consider the case $M\in\Mod_\Lambda(A)$.

(i) Since $(\bullet)\otimes_AM$ is right exact, we have an exact sequence
\[
0\to K
\to\overline{E}^\mu_{k-1}\otimes_AM
\to\overline{E}^\mu_k\otimes_AM
\to(\overline{E}^\mu_k/\overline{E}^\mu_{k-1})\otimes_AM
\to 0
\]
for some $K\in\Mod_\Lambda(A)$.
Since $\omega^*$ is exact, it is sufficient to show $\omega^*K=0$.
This is equivalent to $S_w^{-1}K=0$ for any $w\in W$, which is also equivalent to the injectivity of $S_w^{-1}(\overline{E}^\mu_{k-1}\otimes_AM)
\to S_w^{-1}(\overline{E}^\mu_k\otimes_AM)$ for any $w\in W$.
Thus we have only to show $S_w^{-1}\Tor_A^1(\overline{E}^\mu_{k}/\overline{E}^\mu_{k-1},M)=0$ for any $w\in W$.
By Lemma \ref{lem:filt2} there exists an exact sequence
\[
0\to\overline{E}^\mu_{k}/\overline{E}^\mu_{k-1}
\to F\to C\to 0
\]
of graded $A$-bimodules, where $F$ is isomorphic to $A[\nu_k]^{\oplus m_k}$ as a graded right $A$-module and $C$ belongs to $\Tor_\Lambda(A)$ as a graded left $A$-module.
By $C\in\Tor_\Lambda(A)$ we have $S_w^{-1}C=0$ and hence
\begin{align*}
&S_w^{-1}\Tor_A^n(\overline{E}^\mu_{k}/\overline{E}^\mu_{k-1},M)
=\Tor_A^n(S_w^{-1}(\overline{E}^\mu_{k}/\overline{E}^\mu_{k-1}),M)
=\Tor_A^n(S_w^{-1}F,M)\\
=&S_w^{-1}\Tor_A^n(F,M)=0
\end{align*}
for any $n\ne0$.

(ii) By the proof of (i) we have 
$S_w^{-1}((\overline{E}^\mu_{k}/\overline{E}^\mu_{k-1})\otimes_AM)\simeq S_w^{-1}(F\otimes_AM)$ for any $w\in W$, and hence 
$\omega^*((\overline{E}^\mu_{k}/\overline{E}^\mu_{k-1})\otimes_AM)\simeq 
\omega^*(F\otimes_AM)$.
Hence we have only to show that $F\otimes_AM$ is isomorphic to $A[\nu_k]^{\oplus m_k}\otimes_AM$ as a grade left $A$-module.
Note that under the identification $F=A[\nu_k]\otimes\BF^{m_k}$ of right $A$-modules the left $A$-modules structure on $A[\nu_k]\otimes\BF^{m_k}$ induced by the one on $F$ is given by $\varphi(\psi\otimes v)=\varphi \psi\otimes\tau(\lambda)v$ for $\varphi\in A(\lambda), \psi\in A[\nu_k], v\in\BF^{m_k}$, where $\tau:\Lambda\to GL_{m_k}(\BF)$ is a group homomorphism.
Thus we have an isomorphism 
\[
\delta:F\otimes_AM(=(A[\nu_k]\otimes\BF^{m_k})\otimes_AM)
\to A[\nu_k]^{\oplus m_k}\otimes_AM(=\BF^{m_k}\otimes M[\nu_k])
\]
of grade left $A$-modules given by 
\[
\delta(1\otimes v\otimes m)=\tau(\lambda)^{-1}v\otimes m
\qquad(v\in\BF^{m_k},\,\,m\in M[\nu_k](\lambda)).
\]
\end{proof}

Considering the cases $k=1$ and $k=n$ we have obtained, for $M\in\Mod_\Lambda(A)$ (resp.\ $\Mod_\Lambda(\tU)$),  the canonical monomorphism
\begin{equation}
\label{eq:i-mu}
i_\mu:\omega^*M\to\omega^*(E^\mu\otimes_AM[-w_0\mu])
\end{equation}
and the canonical epimorphism
\begin{equation}
\label{eq:p-mu}
p_\mu:\omega^*(E^\mu\otimes_AM)\to\omega^*M[\mu]
\end{equation}
in $\CM(A)$ (resp.\ $\CM(\tU)$), which are functorial with respect to $M$.

\subsection{The affine open subset $U_{1,q}$}
In this subsection we shall investigate $S_1^{-1}A$.
We denote by $\iota:A\hookrightarrow S_1^{-1}A$ the canonical algebra homomorphism.

For each $\lambda\in\Lambda^+$ take 
$v_\lambda\in V(\lambda)_\lambda$ such that 
$\langle v^*_\lambda,v_\lambda\rangle=1$.
Then we have $\langle c^1_\lambda,u\rangle=\langle v^*_\lambda,uv_\lambda\rangle$ for any $u\in U$.
For simplicity we write $c_\lambda$ instead of $c^1_\lambda$.

Define $r:A\to F^{\geqq0}$ as the composition of 
\[
A\hookrightarrow F\xrightarrow{r_+}F^{\geqq0}.
\]
By $r(c_\lambda)=\chi_\lambda^+$ there exists a unique algebra homomorphism
\begin{equation}
\theta:S_1^{-1}A\to F^{\geqq0}
\end{equation}
such that $\theta\circ\iota=r$ (see \eqref{eq:chi}).
\begin{proposition}
\label{prop:theta}
\begin{itemize}
\item[\rm(i)]
$\theta$ is an isomorphism of $\BF$-algebras.
\item[\rm(ii)]
$\theta$ is an isomorphism of $U^{\geqq0}$-modules, where the $U^{\geqq0}$-module structure on $S_1^{-1}A$ is the restriction of the $U$-module structure given in Proposition \ref{prop:U-action}.
\item[\rm(iii)]
$\theta((S_1^{-1}A)(\lambda))=F^{\geqq0}(\lambda)$ for any $\lambda\in\Lambda$.
\end{itemize}
\end{proposition}
\begin{proof}
By definition we have
\begin{align*}
&e_i(c_\lambda^{-1}\psi)=q^{-(\lambda,\alpha_i)}c_\lambda^{-1}(e_i\psi),
\qquad
k_\mu(c_\lambda^{-1}\psi)=q^{-(\lambda,\mu)}c_\lambda^{-1}(k_\mu \psi),\\
&e_i(\chi_{-\lambda}^+\varphi)=
q^{-(\lambda,\alpha_i)}\chi_{-\lambda}^+(e_i\varphi),
\qquad
k_\mu(\chi_{-\lambda}^+\varphi)=q^{-(\lambda,\mu)}\chi_{-\lambda}^+(k_\mu \varphi),
\end{align*}
for $i\in I, \mu\in\Lambda, \lambda\in\Lambda^+, \psi\in A, \varphi\in F^{\geqq0}$.
Moreover, $r$ is a homomorphism of $U^{\geqq0}$-modules.
Hence $\theta$ is a homomorphism of $U^{\geqq0}$-modules.

By definition we have $r(A(\lambda))\subset F^{\geqq0}(\lambda)$ for any $\lambda\in\Lambda^+$.
In particular, we have $r(\chi_\lambda^+)\in F^{\geqq0}(\lambda)$.
Since $\theta$ is a ring homomorphism, we obtain 
\begin{equation}
\label{eq:theta}
\theta((S_1^{-1}A)(\lambda))\subset F^{\geqq0}(\lambda) 
\qquad(\lambda\in\Lambda).
\end{equation}

By \eqref{eq:CqB} and \eqref{eq:theta} it is sufficient to show that 
\[
\theta_0=\theta|(S_1^{-1}A)(0):(S_1^{-1}A)(0)\to F^{\geqq0}(0) 
\]
is an isomorphism.
Assume that $x\in\Ker(\theta_0)$.
There exists some $\lambda\in\Lambda^+$ and $\varphi\in A(\lambda)$ such that $x=\varphi c_\lambda^{-1}$. 
Then we have $\varphi\in\Ker(\theta)\cap A(\lambda)$.
Take $v\in V(\lambda)$ such that 
$\varphi=f_\lambda(v)$.
By $f\in\Ker(\theta)$ we have
\[
\langle V^*(\lambda),v\rangle=
\langle v^*_\lambda U^{\geqq0},v\rangle=
\langle v^*_\lambda,U^{\geqq0}v\rangle=\{0\}.
\]
This implies $v=0$. 
Hence $x=0$.

It remains to show the surjectivity of $\theta_0$.
By Proposition \ref{prop:CqB} (ii) it is sufficient to show that for any $\gamma\in Q^+$ there exists some $\lambda\in\Lambda^+$ such that the linear map 
\[
V(\lambda)_{\lambda-\gamma}\to(U^+_\gamma)^*\quad
(v\mapsto \theta(f_\lambda(v)c_\lambda^{-1})|U^+_\gamma)
\]
is surjective.
By definition we have
\begin{align*}
&\langle\theta(f_\lambda(v)c_\lambda^{-1}),u\rangle
=\langle(r_+(f_\lambda(v)))\chi^+_{-\lambda},u\rangle
=\langle(r_+(f_\lambda(v)))\otimes\chi^+_{-\lambda},\Delta(u)\rangle\\
=&\langle(r_+(f_\lambda(v))),u\rangle
=\langle v^*_\lambda u, v\rangle.
\end{align*}
Hence we obtain the desired result by Lemma \ref{lem:V}.
\end{proof}
\begin{proposition}
\label{prop:S-1A}
For any $\lambda\in\Lambda$ $(S_1^{-1}A)(\lambda)$ is isomorphic to $T^*(\lambda)$ as a $U$-module.
\end{proposition}
\begin{proof}
By Proposition \ref{prop:CqB} and Proposition \ref{prop:theta} we have
$\ch((S_1^{-1}A)(\lambda))=\ch(T(\lambda))=\ch(T^*(\lambda))$.
Moreover, the restricted dual $((S_1^{-1}A)(\lambda))^\bigstar$ is a rank one free right $U^+$-module.
Hence $((S_1^{-1}A)(\lambda))^\bigstar$ is isomorphic to $T_{\rm r}(\lambda)$ as a right $U$-module.
It follows that we have $(S_1^{-1}A)(\lambda)\simeq(T_{\rm r}(\lambda))^\bigstar=T^*(\lambda)$.
\end{proof}

\subsection{Coherent sheaves with $U$-actions}
Let $\Mod^{\diamond}_{\Lambda}(\tU)$ be the full subcategory of 
$\Mod_\Lambda(\tU)$ consisting of objects $M$ of $\Mod_\Lambda(\tU)$ such that
$M=\bigoplus_{\lambda\in\Lambda}M_\lambda$, and let 
$\Mod^{f}_{\Lambda}(\tU)$ be the full subcategory of
$\Mod^{\diamond}_{\Lambda}(\tU)$ consisting of 
objects $M$ of $\Mod^{\diamond}_{\Lambda}(\tU)$ such that
$M$ is finitely generated as an $A$-module.

Let $M\in\Mod^{f}_{\Lambda}(\tU)$.
Since $M$ is a finitely generated $A$-module, we have $\dim M(\xi)<\infty$ for any $\xi\in\Lambda$.
Hence we have $M(\xi)\in\Mod^f(U)$ for any $\xi\in\Lambda$.

Let
\begin{align}
&\CM^{\diamond}(\tU)=
\frac{\Mod^{\diamond}_{\Lambda}(\tU)}
{\Tor_\Lambda(\tU)\cap\Mod^{\diamond}_{\Lambda}(\tU)}
=
(\Sigma^{\diamond}_{\tU})^{-1}\Mod^{\diamond}_{\Lambda}(\tU),\\
&\CM^{f}(\tU)=
\frac{\Mod^{f}_{\Lambda}(\tU)}
{\Tor_\Lambda(\tU)\cap\Mod^{f}_{\Lambda}(\tU)}
=
(\Sigma^{f}_{\tU})^{-1}\Mod^{f}_{\Lambda}(\tU),
\end{align}
where $\Sigma^{\diamond}_{\tU}$ and $\Sigma^{f}_{\tU}$ are the collection of morphisms in 
$\Mod^{\diamond}_{\Lambda}(\tU)$ and $\Mod^{f}_{\Lambda}(\tU)$ whose kernel and cokernel belong to 
$\Tor_\Lambda(\tU)\cap\Mod^{\diamond}_{\Lambda}(\tU)$
and $\Tor_\Lambda(\tU)\cap\Mod^{f}_{\Lambda}(\tU)$ respectively.
Denote by
\[
\omega_\diamond^{*}:\Mod^{f}_{\Lambda}(\tU)\to\CM^{\diamond}(\tU), 
\qquad
\omega_f^{*}:\Mod^{f}_{\Lambda}(\tU)\to\CM^{f}(\tU)
\]
the localization functors.

We shall regard $\CM^{\diamond}(\tU)$ and $\CM^{f}(\tU)$ as full subcategories of $\CM(\tU)$ by Lemma \ref{lem:FF}.
We sometimes write $\omega^*$ instead of $\omega_\diamond^*$ and $\omega_f^*$.

\begin{proposition}
\label{prop:irred-obj}
$\omega^*A[\lambda]$ is an irreducible object of $\CM^{f}(\tU)$ for any $\lambda\in\Lambda$, and any irreducible object of $\CM^{f}(\tU)$ is isomorphic to $\omega^*A[\lambda]$ for some $\lambda\in\Lambda$.
\end{proposition}
\begin{proof}
Let us show that $\omega^*A[\lambda]$ is irreducible.
We may assume that $\lambda=0$.
Assume that $\overline{M}$ is a non-zero subobject of $\omega^*A$.
Set $M=\omega_*\overline{M}$.
Then $M$ is a subobject of $A=\omega_*\omega^*A$ 
(see  Proposition \ref{prop:Borel-Weil}), and we have $\overline{M}=\omega^*M$ by $\omega^*\circ\omega_*=\Id$.
Since $M$ is non-zero, there exists some $\mu\in\Lambda^+$ such that $M(\mu)\ne\{0\}$.
Since $A(\mu)$ is an irreducible $U$-module, we have $M(\mu)=A(\mu)$.
It follows that $M\supset\bigoplus_{\xi\in\mu+\Lambda^+}A(\xi)$ by Lemma \ref{lem:grading}, and hence $A/M\in\Tor_\Lambda(A)$.
Thus we have $\overline{M}=\omega^*M=\omega^*A$.

Assume that $\overline{M}$ is an irreducible object of $\CM^{f}(\tU)$.
Take $M\in\Mod^{f}_{\Lambda}(\tU)$ such that $\omega^*M=\overline{M}$.
We may assume $\Tor(M)=\{0\}$.
By $M\ne\{0\}$ there exists some $\xi$ such that $M(\xi)\ne\{0\}$.
Note that $M(\xi)\in\Mod^f(U)$ as a $U$-module.
Take an irreducible $U$-submodule $V$ of $M(\xi)$ and set $N=\tU V=AUV=AV\subset M$.
By $N/\Tor(N)\simeq N\ne\{0\}$ we have $\omega^*N\ne\{0\}$.
Hence the irreducibility of $\overline{M}=\omega^*M$ implies $\omega^*M=\omega^*N$.
Thus we may assume that $M=\tU M(\xi)$ and that $M(\xi)$ is isomorphic to $V(\mu)$ as a $U$-module.
Hence $M$ is isomorphic to a quotient of  $E^\mu[-\xi]=A[-\xi]\otimes V(\mu)$.
It follows that $\overline{M}$ is isomorphic to a quotient of  $\omega^*E^\mu[-\xi]$.
By Lemma \ref{lem:filt0} and Lemma \ref{lem:Filtration} there exists an increasing sequence
\[
0=N_0\subset N_1\subset\cdots\subset N_n=\omega^*E^\mu
\]
of subobjects of $\omega^*E^\mu\in\CM(\tU)$ such that 
for each $k$ we have 
$N_k/N_{k-1}\simeq\omega^*A[\lambda_k]$ 
for some $\lambda_k\in\Lambda$.
This implies that $\overline{M}$ is isomorphic to 
$\omega^*A[\lambda_k-\xi]$ 
for some $k$.
\end{proof}

\begin{lemma}
\label{lem:PP}
For $M\in\Mod^{\diamond}_{\Lambda}(\tU)$ $($resp.\ $M\in\Mod^{f}_{\Lambda}(\tU)$$)$ we have $\omega_*\omega^*M\in\Mod^{\diamond}_{\Lambda}(\tU)$ $($resp.\ $M\in\Mod^{f}_{\Lambda}(\tU)$$)$.
\end{lemma}
\begin{proof}
Assume that $M\in\Mod^{\diamond}_{\Lambda}(\tU)$.
By Corollary \ref{cor:cor-quasi-scheme} $\omega_*\omega^*M$ is a subobject of $\bigoplus_{w\in W}S_w^{-1}M$. 
Since $\bigoplus_{w\in W}S_w^{-1}M\in\Mod^{\diamond}_{\Lambda}(\tU)$, we have $\omega_*\omega^*M\in\Mod^{\diamond}_{\Lambda}(\tU)$.

It remains to show that $\omega_*\overline{M}$ is a finitely generated $A$-module for any $\overline{M}\in\CM^f(\tU)$.
By Proposition \ref{prop:property-A} and Proposition \ref{prop:irred-obj} we may assume that $\overline{M}=\omega^*A[\lambda]$ for some $\lambda\in\Lambda$.
In this case the assertion follows from Proposition \ref{prop:Borel-Weil}.
\end{proof}

We need the following result later.
\begin{lemma}
\label{lem:fin-inj}
Let $w\in W$.
Let $M\in\Mod^{\diamond}_{\Lambda}(\tU)$ and set
\[
\Gamma(\omega^*M)^{\fin}=
\{m\in \Gamma(\omega^*M)\mid \dim_\BF Um<\infty
\}.
\]
Then the canonical homomorphism 
$\Gamma(\omega^*M)^{\fin}\to (S_w^{-1}M)(0)$ is injective.
\end{lemma}
\begin{proof}
By Lemma \ref{lem:PP} we may assume that $M\simeq\omega_*\omega^*M$.
Take $m\in\Gamma(\omega^*M)^{\fin}=M(0)^{\fin}$.
Let $N$ be the $\tU$-submodule of $M$ generated by $m$.
Then $N$ belongs to $\Mod_\Lambda^f(\tU)$, and the canonical morphisms $\Gamma(\omega^*N)\to\Gamma(\omega^*M)$ and $(S_w^{-1}N)(0)\to  (S_w^{-1}M)(0)$ are injective.
Hence we may assume that $M\in\Mod^{f}_{\Lambda}(\tU)$.
In this case we have $\Gamma(\omega^*M)^{\fin}=\Gamma(\omega^*M)$ by Lemma \ref{lem:PP}.
Assume that there exists a short exact sequence
\[
0\to \omega^*{M}_1\to \omega^*{M}_2\to\omega^*{M}_3\to 0
\]
in $\CM^f(\tU)$.
Then we have a commutative diagram:
\[
\begin{CD}
0@>>>\Gamma(\omega^*{M}_1)@>>>\Gamma(\omega^*{M}_2)@>>>\Gamma(\omega^*{M}_3)
\\
@.@VVV@VVV@VVV
\\
0@>>>S_w^{-1}M_1@>>>S_w^{-1}M_2@>>>S_w^{-1}M_3@>>>0
\end{CD}
\]
whose rows are exact.
Hence we may assume that $M=A[\lambda]$ for some $\lambda\in\Lambda$ by Proposition \ref{prop:irred-obj}.
In this case the assertion is a consequence of Proposition \ref{prop:Borel-Weil} and the injectivity of $A\to S_w^{-1}A$.
\end{proof}

\section{$D$-modules}
\subsection{$q$-differential operators}
\label{subsec:QD}
For $\varphi\in A, u\in U, \lambda\in\Lambda$ we define $\ell_\varphi, r_\varphi, \deru_u, \sigma_\lambda\in\End_\BF(A)$ by
\[
\ell_\varphi(\psi)=\varphi\psi,\quad
r_\varphi(\psi)=\psi\varphi,\quad
\deru_u(\psi)=u\psi,\quad
\sigma_\lambda(\psi)=q^{(\lambda,\mu)}\psi
\]
for $\psi\in A(\mu)$.
We define a subalgebra $D$ of $\End_\BF(A)$ by
\begin{equation}
D=\langle \ell_\varphi, r_\varphi, \deru_u, \sigma_\lambda\mid
\varphi\in A, u\in U, \lambda\in\Lambda\rangle.
\end{equation}
For $\lambda\in\Lambda$ we set 
\[
D(\lambda)=\{d\in D\mid d(A(\xi))\subset A(\lambda+\xi)\,\,(\xi\in\Lambda)\}.
\]
Since $\ell_\varphi, r_\varphi\in D(\lambda)$ for $\varphi\in A(\lambda)$ and $\deru_u, \sigma_\lambda\in D(0)$, we have
\begin{equation}
D=\bigoplus_{\lambda\in\Lambda^+}D(\lambda),\qquad
D(0)=\langle \deru_u, \sigma_\lambda\mid u\in U, \lambda\in\Lambda\rangle.
\end{equation}
In particular, $D$ is a $\Lambda$-graded $\BF$-algebra.
\begin{lemma}
\label{lem:rl}
Write
\[
\sum_{\beta\in Q^+}q^{(\beta,\beta)}(1\otimes k_\beta)(S\otimes\id)(\Xi_\beta)
=\sum_px_p\otimes y_p
\]
$($see Section \ref{subsec:R} for the notation$)$.
Then we have
\begin{align}
&r_\psi=\sum_p
\ell_{x_p\psi}\deru_{y_pk_\eta}\sigma_{-\mu}\qquad
&(\psi\in A(\mu)_\eta),
\label{eq:rl1}\\
&\ell_\varphi=\sum_p
r_{y_p\varphi}\deru_{x_pk_\xi}\sigma_{-\lambda}
&(\varphi\in A(\lambda)_\xi).
\label{eq:rl2}
\end{align}
\end{lemma}
\begin{proof}
Let $\varphi\in A(\lambda)_\xi, \psi\in A(\mu)_\eta$.
Take $v_0\in V(\lambda)_\xi, v_1\in V(\mu)_\eta$ such that $\varphi=f_\lambda(v_0),\,\, \psi=f_\mu(v_1)$.
Set $\CR^\vee=\CR_{V(\mu),V(\lambda)}^\vee$.
By Proposition \ref{prop:R} (ii) we have
\begin{align*}
&\langle\varphi\psi,u\rangle
=\langle\varphi\otimes\psi,\Delta(u)\rangle
=\langle v_\lambda^*\otimes v_\mu^*,u(v_0\otimes v_1)\rangle\\
=&\langle v_\lambda^*\otimes v_\mu^*,\CR^\vee u(\CR^\vee)^{-1}(v_0\otimes v_1)\rangle
=\langle{}^t(\CR^\vee)(v_\lambda^*\otimes v_\mu^*), u(\CR^\vee)^{-1}(v_0\otimes v_1)\rangle
\end{align*}
By $v^*_\lambda U^-_{-\beta}=\{0\}$ for $\beta\in Q^+\setminus\{0\}$ 
we have
${}^t(\CR^\vee)(v_\lambda^*\otimes v_\mu^*)
=q^{-(\lambda,\mu)}v_\mu^*\otimes v_\lambda^*.$
By
\[
\CR_{V(\mu),V(\lambda)}^{-1}=
(\sum_px_p\otimes y_p)\kappa_{V(\mu),V(\lambda)}
\]
we have
$(\CR^\vee)^{-1}(v_0\otimes v_1)
=q^{(\xi,\eta)}\sum_px_pv_1\otimes y_pv_0.$
Hence
\begin{align*}
\langle\varphi\psi,u\rangle
=&q^{(\xi,\eta)-(\lambda,\mu)}\sum_p
\langle v_\mu^*\otimes v_\lambda^*,u(x_pv_1\otimes y_pv_0)\rangle\\
=&q^{(\xi,\eta)-(\lambda,\mu)}
\sum_p
\langle(x_p\psi)(y_p\varphi),u\rangle.
\end{align*}
It follows that we have
\[
\varphi\psi=
q^{(\xi,\eta)-(\lambda,\mu)}\sum_p
(\deru_{x_p}\psi)(\deru_{y_p}\varphi).
\]
\end{proof}
\begin{corollary}
We have
\[
D=\langle \ell_\varphi, \deru_u, \sigma_\lambda\mid
\varphi\in A, u\in U, \lambda\in\Lambda\rangle=
\langle r_\varphi, \deru_u, \sigma_\lambda\mid
\varphi\in A, u\in U, \lambda\in\Lambda\rangle.
\]
\end{corollary}
We can easily check the following.
\begin{align}
&\ell_\varphi\ell_\psi=\ell_{\varphi\psi}
&(\varphi, \psi\in A),
\label{eq:comm1}
\\
&\sigma_\lambda\sigma_\mu=\sigma_{\lambda+\mu}
&(\lambda, \mu\in\Lambda),
\label{eq:comm2}
\\
&\deru_{u}\deru_{u'}=\deru_{uu'}
&(u,u'\in U),
\label{eq:comm3}
\\
&\sigma_\lambda\ell_\varphi=q^{(\lambda,\mu)}\ell_\varphi\sigma_\lambda
&(\lambda, \mu\in\Lambda, \varphi\in A(\mu)),
\label{eq:comm4}
\\
&\sigma_\lambda\deru_u=\deru_u\sigma_\lambda
&(\lambda\in\Lambda, u\in U),
\label{eq:comm5}
\\
&\deru_u\ell_\varphi=\sum_{(u)_1}\ell_{u_{(0)}\varphi}\deru_{u_{(1)}}
&(u\in U,\varphi\in A).
\label{eq:comm6}
\end{align}
In particular we have homomorphisms 
\begin{align*}
&\deru:U\to D\quad(u\mapsto\deru_u),\qquad
\ell:A\to D\quad(\varphi\mapsto\ell_\varphi),\\
&\sigma:\BF[\Lambda]\to D\quad(\Lambda\ni\lambda\mapsto\sigma_\lambda)
\end{align*}
of $\BF$-algebras.
Moreover, $\deru$ and $\ell$ induces a homomorphism 
\begin{equation}
\label{eq:tUD}
\tU\to D\quad(\varphi u\mapsto \ell_{\varphi}\deru_u)
\end{equation}
of $\Lambda$-graded $\BF$-algebras.
In particular, we have
\begin{equation}
\label{eq:tUD2}
\ell_{\varphi}\,\deru_{u}=
\sum_{(u)_1}\deru_{u_{(1)}}\ell_{{S^{-1}u_{(0)}}\varphi}
\qquad(u\in U,\,\,\varphi\in A).
\end{equation}
by \eqref{eq:tU10}.

By Proposition \ref{prop:center2} (ii) we have
\begin{equation}
\label{equ:deruz}
\deru_z=\sigma\circ\zeta(z)\qquad(z\in\Gz)
\end{equation}
where $\zeta:\Gz\to\BF[\Lambda]$ is as in Section \ref{subsection:center}.

We shall identify $A$ with a subalgebra of $D$ via the injective $\BF$-algebra homomorphism $\ell:A\to D$  (the injectivity follows from Proposition \ref{prop:property-A} (i)).

\begin{proposition}
\label{prop:OreD}
Let $w\in W$.
\begin{itemize}
\item[{\rm(i)}]
$S_w$ satisfies the left and right Ore conditions in $D$.
\item[{\rm(ii)}]
The canonical homomorphism $D\to S_w^{-1}D$ is injective.
\end{itemize}
\end{proposition}
\begin{proof}
It is sufficient to show the following.
\begin{itemize}
\item[(a)]
For any $d\in D$ and $s\in S_w$ there exists some $t\in S_w$ and $d'\in D$ satisfying $td=d's$.
\item[(b)]
For any $d\in D$ and $s\in S_w$ there exists some $t\in S_w$ and $d'\in D$ satisfying $dt=sd'$.
\item[(c)]
If $sd=0$ for $d\in D$ and $s\in S_w$, then we have $d=0$.
\item[(d)]
If $ds=0$ for $d\in D$ and $s\in S_w$, then we have $d=0$.
\end{itemize}
The statements (a), (b) is proved similarly to Proposition \ref{prop:OretU}.
The statement (c) follows from Proposition \ref{prop:property-A} (i).
The statement (d) is equivalent to the injectivity of $D\to D\otimes_AS_w^{-1}A$.
By Lunts-Rosenberg \cite[Section 1.2]{LRD} $\tilde{D}\to\tilde{D}\otimes_AS_w^{-1}A$ is injective for a ring $\tilde{D}$ containing $D$, and the assertion for $D$ follows from this.
\end{proof}
Hence $S_w^{-1}D$ is a $\Lambda$-graded $\BF$-algebra, and we have
\[
S_w^{-1}D\simeq S_w^{-1}A\otimes_AD\simeq D\otimes_AS_w^{-1}A
\simeq S_w^{-1}A\otimes_AD\otimes_AS_w^{-1}A.
\]
Moreover, for $M\in\Mod_\Lambda(D)$ we have
\[
S_w^{-1}M=S_w^{-1}D\otimes_DM\in\Mod_\Lambda(S_w^{-1}D).
\]
In particular $S_w^{-1}A$ is a graded $S_w^{-1}D$-module.

Since $A$ is a sum of $U$-submodules contained in $\Mod^f(U)$, we have the operator $T_w:A\to A$.
For $d\in\End_\BF(A)$ we set 
\[
Z_w(d)=T_w^{-1}\circ d\circ T_w\in\End_\BF(A).
\]
\begin{proposition}
\label{prop:D-local}
$Z_w$ induces an algebra automorphisms of $D$ and an algebra isomorphism $S_1^{-1}D\to S_w^{-1}D$.
\end{proposition}
\begin{proof}
We see easily that $Z_{s_i}(\deru_u)=\deru_{T_i^{-1}(u)}, Z_{s_i}(\sigma_\mu)=\sigma_\mu$ for any $u\in U$ and $\mu\in\Lambda$.
For $\varphi, \psi\in A$ we have
\[
Z_{s_i}(\ell_\varphi)(\psi)
=T_i^{-1}(\varphi T_i(\psi))
=m(T_i^{-1}(\varphi\otimes T_i(\psi)))
\]
since the multiplication $m:A\otimes A\to A$ is a homomorphism of $U$-modules.
By Lemma \ref{lem:Ti} we have
\[
T_i^{-1}=\exp_{q_i^{-1}}(-(q_i-q_i^{-1})f_i\otimes e_i)(T_i^{-1}\otimes T_i^{-1})
\]
as an operator on $A\otimes A$.
Write $\exp_{q_i^{-1}}(-(q_i-q_i^{-1})f_i\otimes e_i)=\sum_pb_p\otimes a_p$.
Then we have
\[
Z_{s_i}(\ell_\varphi)(\psi)
=\sum_p\deru_{b_p}(T_i^{-1}\varphi)\deru_{a_p}(\psi),
\]
and hence
\[
Z_{s_i}(\ell_\varphi)=\sum_p\ell_{\deru_{b_p}(T_i^{-1}\varphi)}\deru_{a_p}
\in D
\]
for any $\varphi\in A$.
Similarly, we can show $Z_{s_i}^{-1}(\ell_\varphi)\in D$ for any $\varphi\in A$.
Hence $Z_w$ induces an algebra automorphism of $D$.

It remains to show $Z_w(\ell_{S_1})=\ell_{S_w}$.
It is sufficient to show that for $w\in W, \,\,i\in I$ such that $w(\alpha_i)\in\Delta^+$ we have $Z_{s_i}(\ell_{S_w})=\ell_{S_{ws_i}}$.
Let $\varphi\in S_w$. 
Take $\lambda\in \Lambda^+$ such that $\varphi\in A(\lambda)_{w^{-1}\lambda}\setminus\{0\}$.
Then we have $T_i^{-1}\varphi\in A(\lambda)_{s_iw^{-1}\lambda}\setminus\{0\}$.
By $(s_iw^{-1}\lambda,\alpha_i^\vee)=-(\lambda,w\alpha_i^\vee)\leqq0$ we obtain 
\[
Z_{s_i}(\ell_\varphi)=\sum_p\ell_{\deru_{b_p}(T_i^{-1}\varphi)}\deru_{a_p}
=\ell_{T_i^{-1}\varphi}
\in\ell_{S_{ws_i}}.
\]
\end{proof}

\subsection{Category of $D$-modules}
Applying results in Section \ref{subsec:GS} to $C=D$ we have the localization functor 
\[
\omega^*=\omega_D^*:\Mod_\Lambda(D)\to\CM(D)
:=
\frac{\Mod_{\Lambda}(D)}
{\Tor_\Lambda(D)}
=
\Sigma_{D}^{-1}\Mod_{\Lambda}(D)
\]
and its right adjoint
\[
\omega_*=\omega_{D*}:\CM(D)\to\Mod_\Lambda(D).
\]
Taking the degree zero part of $\omega_*$ we have the global section functor
\begin{equation}
\Gamma:\CM(D)\to\Mod(D(0)).
\end{equation}
Define a functor
\begin{equation}
\CL:\Mod(D(0))\to\CM(D)
\end{equation}
by $\CL(N)=\omega^*(D\otimes_{D(0)}N)$, and we call it the localization functor.
Since $D\otimes_{D(0)}(\bullet):\Mod(D(0))\to\Mod_\Lambda(D)$ is left adjoint to 
$\Mod_\Lambda(D)\ni M\mapsto M(0)\in\Mod(D(0))$, we have the following.
\begin{lemma}
The functor $\CL:\Mod(D(0))\to\CM(D)$ is left adjoint to $\Gamma:\CM(D)\to\Mod(D(0))$.
\end{lemma}

For $\lambda\in\Lambda$ we denote by $\Mod_{\Lambda,\lambda}(D)$ the full subcategory of $\Mod_\Lambda(D)$ consisting of $M\in\Mod_\Lambda(D)$ satisfying $\sigma_\mu|M(\xi)=q^{(\mu,\lambda+\xi)}\id$ for any $\mu, \xi\in\Lambda$.
Let
\begin{equation}
\CM_\lambda(D)=
\frac{\Mod_{\Lambda,\lambda}(D)}
{\Tor_\Lambda(D)\cap\Mod_{\Lambda,\lambda}(D)}
=
\Sigma_{D,\lambda}^{-1}\Mod_{\Lambda,\lambda}(D),
\end{equation}
where $\Sigma_{D,\lambda}$ is the collection of morphisms in 
$\Mod_{\Lambda,\lambda}(D)$ whose kernel and cokernel belong to 
$\Tor_\Lambda(D)\cap\Mod_{\Lambda,\lambda}(D)$,
and denote by
\[
\omega_\lambda^*:\Mod_{\Lambda,\lambda}(D)\to\CM_\lambda(D)
\]
the localization functor.

We shall regard $\CM_\lambda(D)$ as a full subcategory of $\CM(D)$ by Lemma \ref{lem:FF}, and we often write $\omega^*$ instead of $\omega_\lambda^*$.
\begin{lemma}
\label{lem:cd1}
For $M\in\Mod_{\Lambda,\lambda}(D)$ we have $(R^r\omega_*)(\omega^*M)\in\Mod_{\Lambda,\lambda}(D)$ for any $r$.
\end{lemma}
\begin{proof}
For $\mu\in\Lambda$ and $N\in\Mod_\Lambda(D)$ we define $h^\mu_N:N\to N$ by $h^\mu_N|N(\xi)=q^{-(\mu,\xi)}\sigma_\mu|N(\xi)$ for any $\xi\in\Lambda$.
We see easily that $h^\mu_N$ is a morphism in $\Mod_\Lambda(D)$.
Moreover, $h_N^\mu$ is functorial with respect to $N$ in the sense that for a morphism $f:N_1\to N_2$ in $\Mod_\Lambda(D)$ we have $h_{N_2}^\mu\circ f=f\circ h_{N_1}^\mu$.

Let us show $(R^r\omega_*)(\omega^*h^\mu_N)=h^\mu_{(R^r\omega_*)(\omega^*N)}$ for any $r$.
By the functoriality of  $h_N^\mu$ stated above, it is sufficient to show $\omega_*\omega^*h^\mu_N=h^\mu_{\omega_*\omega^*N}$.
Let $j:N\to\omega_*\omega^*N$ be the canonical morphism.
By the definition of adjoint functors we have $\omega_*\omega^*h_N^\mu\circ j=j\circ h_N^\mu$, and by the functoriality of $h_N^\mu$ we have $h^\mu_{\omega_*\omega^*N}\circ j=j\circ h_N^\mu$.
Hence we obtain $\omega_*\omega^*h^\mu_N=h^\mu_{\omega_*\omega^*N}$ by
\begin{align*}
&\Hom(\omega_*\omega^*N,\omega_*\omega^*N)\simeq
\Hom(\omega^*\omega_*\omega^*N,\omega^*N)\simeq
\Hom(\omega^*N,\omega^*N)\\
\simeq&
\Hom(N,\omega_*\omega^*N).
\end{align*}
Here we have used $\omega^*\circ\omega_*\simeq\Id$ and the fact that $\omega^*$ is left adjoint to $\omega_*$.

By $M\in\Mod_{\Lambda,\lambda}(D)$ we have 
$h^\mu_M=q^{(\mu,\lambda)}\id$.
Hence 
$h^\mu_{(R^r\omega_*)(\omega^*M)}=(R^r\omega_*)(\omega^*h^\mu_M)
=q^{(\mu,\lambda)}(R^r\omega_*)(\omega^*\id)
=q^{(\mu,\lambda)}\id$, 
from which we obtain the desired result.
\end{proof}
By Lemma \ref{lem:cd1} the restriction of $\omega_*:\CM(D)\to\Mod_\Lambda(D)$ to $\CM_\Lambda(D)$ gives a left exact functor
\begin{equation}
\omega_{\lambda*}:\CM_\lambda(D)\to\Mod_{\Lambda,\lambda}(D).
\end{equation}
We have
\begin{align*}
&\Hom_{\CM_\lambda(D)}(\omega_\lambda^*M,\omega_\lambda^*N)
\simeq
\Hom_{\CM(D)}(\omega^*M,\omega^*N)\\
\simeq&
\Hom_{\Mod_\Lambda(D)}(M,\omega_*\omega^*N)
\simeq
\Hom_{\Mod_{\Lambda,\lambda}(D)}(M,\omega_{\lambda*}\omega_\lambda^*N),
\end{align*}
and hence the functor $\omega_{\lambda*}$ is right adjoint to $\omega_\lambda^*$.

\subsection{Beilinson-Bernstein correspondence}
Define $D_\lambda\in\Mod_{\Lambda,\lambda}(D)$ by
\begin{equation}
D_\lambda
=D/\sum_{\mu\in\Lambda}D(\sigma_\mu-q^{(\mu,\lambda)})
\end{equation}
Since $\sigma_\mu$ for $\mu\in\Lambda$ belongs to the center of $D(0)$ we have a natural $\BF$-algebra structure on 
 \[
 D_\lambda(0)=D(0)/\sum_{\mu\in\Lambda}
 D(0)(\sigma_\mu-q^{(\mu,\lambda)}).
 \]
 Since $D(0)$ is generated by the elements of the form $\deru_u, \sigma_\mu$ for $u\in U, \mu\in\Lambda$, we have  a natural surjective algebra homomorphism $U\to D_\lambda(0)$.
Set
\[
J_\lambda=\sum_{z\in\Gz}U(z-\zeta_\lambda(z)).
\]
By \eqref{equ:deruz} we have $J_\lambda\subset\Ker(U\to D_\lambda(0))$ and hence we obtain a surjective algebra homomorphism
\begin{equation}
U/J_\lambda\to D_\lambda(0).
\end{equation}

By Lemma \ref{lem:cd1} we obtain a left exact functor
\begin{equation}
\Gamma_\lambda:\CM_\lambda(D)\to\Mod(D_\lambda(0))
\end{equation}
as the restriction of $\Gamma:\CM(D)\to\Mod(D(0))$.
By restricting the functor $\CL:\Mod(D(0))\to\CM(D)$ to $\Mod(D_\lambda(0))$ we obtain a right exact functor
\begin{equation}
\CL_\lambda:\Mod(D_\lambda(0))\to\CM_\lambda(D).
\end{equation}
We see easily the following.
\begin{lemma}
\label{lem:adjoint}
The functor $\CL_\lambda:\Mod(D_\lambda(0))\to\CM_\lambda(D)$ is left adjoint to $\Gamma_\lambda:\CM_\lambda(D)\to\Mod_\lambda(D(0))$.
\end{lemma}
The rest of this paper is devoted to proving the following theorems.
\begin{theorem}
\label{thm:A}
\begin{itemize}
\item[\rm(i)]
If $\lambda+\rho\in\Lambda^+$, then 
$\Gamma_\lambda:\CM_\lambda(D)\to\Mod(D_\lambda(0))$ is exact.
\item[\rm(ii)]
Assume $\lambda\in\Lambda^+$.
If $\Gamma_\lambda(M)=0$ for $M\in\CM_\lambda(D)$, then we have $M=0$.
\end{itemize}
\end{theorem}
\begin{theorem}
\label{thm:B}
\begin{itemize}
\item[\rm(i)]
The canonical homomorphism $U/J_\lambda\to D_\lambda(0)$ is an isomorphism for any $\lambda\in\Lambda$.
\item[\rm(ii)]
Assume $\lambda+\rho\in\Lambda^+$.
Then the canonical homomorphism $D_\lambda(0)\to\Gamma_\lambda(\omega^*D_\lambda)$ is an isomorphism.
\end{itemize}
\end{theorem}

By a standard argument Theorem \ref{thm:A} and Theorem \ref{thm:B} imply the following.
\begin{theorem}
\label{thm:C}
Assume $\lambda\in\Lambda^+$.
Then $\Gamma_\lambda:\CM_\lambda(D)\to\Mod(U/J_\lambda)$ gives an equivalence of categories, and its inverse is given by $\CL_\lambda$.
\end{theorem}

For the sake of completeness we give a proof of Theorem \ref{thm:C} assuming Theorem \ref{thm:A} and Theorem \ref{thm:B} in the rest of this subsection.

Let $M\in\Mod(D_\lambda(0))$.
We can take an exact sequence of the form
\[
D_\lambda(0)^{\oplus J_1}
\to D_\lambda(0)^{\oplus J_2}
\to M\to 0.
\]
By Theorem \ref{thm:A} (i) the functor $\Gamma_\lambda\circ\CL_\lambda$ is right exact.
By Theorem \ref{thm:B} we have $\Gamma_\lambda\CL_\lambda(D_\lambda(0))\simeq D_\lambda(0)$.
Thus we obtain a commutative diagram
\[
\begin{CD}
D_\lambda(0)^{\oplus J_1}
@>>>
D_\lambda(0)^{\oplus J_2}
@>>>
M
@>>>
0
\\
@| @| @VVV
\\
D_\lambda(0)^{\oplus J_1}
@>>>
D_\lambda(0)^{\oplus J_2}
@>>>
\Gamma_\lambda\CL_\lambda(M)
@>>>
0
\end{CD}
\]
whose rows are exact.
Hence the canonical morphism $M\to\Gamma_\lambda\CL_\lambda(M)$ is an isomorphism.
It follows that $\Id\to\Gamma_\lambda\circ\CL_\lambda$ is an isomorphism.

It remains to show that $\CL_\lambda\circ\Gamma_\lambda\to\Id$ is an isomorphism.
By Lemma \ref{lem:adjoint} the composition of 
\[
\Gamma_\lambda\to
(\Gamma_\lambda\circ\CL_\lambda)\circ\Gamma_\lambda
=\Gamma_\lambda\circ(\CL_\lambda\circ\Gamma_\lambda)
\to\Gamma_\lambda
\]
coincides with $\Id$.
Since $\Id\to\Gamma_\lambda\circ\CL_\lambda$ is an isomorphism,
the canonical morphism 
$\Gamma_\lambda\circ(\CL_\lambda\circ\Gamma_\lambda)
\to\Gamma_\lambda$ is an isomorphism.
Let $N\in\CM_\lambda(D)$.
Setting $K_1=\Ker(\CL_\lambda\Gamma_\lambda(N)\to N),\,\,K_2=\Coker(\CL_\lambda\Gamma_\lambda(N)\to N)$ we have an exact sequence
\[
0\to K_1\to\CL_\lambda\Gamma_\lambda(N)\to N\to K_2\to 0
\]
By applying the exact functor $\Gamma_\lambda$ we obtain
\[
0\to \Gamma_\lambda(K_1)\to
\Gamma_\lambda\CL_\lambda\Gamma_\lambda(N)\to 
\Gamma_\lambda(N)\to \Gamma_\lambda(K_2)\to 0.
\]
Since $\Gamma_\lambda\CL_\lambda\Gamma_\lambda(N)\to 
\Gamma_\lambda(N)$ is an isomorphism, we have $\Gamma_\lambda(K_1)=\Gamma_\lambda(K_2)=0$, which implies $K_1=K_2=0$ by Theorem \ref{thm:A} (ii).
Hence $\CL_\lambda\circ\Gamma_\lambda\to\Id$ is an isomorphism.

\subsection{The key lemma}
\label{subsec:key-lemma}
Sections \ref{subsec:key-lemma} and \ref{subsec:pf-of-A} are devoted to the proof of Theorem \ref{thm:A}.
The arguments are identical with those in Lunts-Rosenberg \cite{LR} except for the usage of Proposition \ref{prop:filt-E1}, which was a conjecture in \cite{LR}.
We shall reproduce the arguments in \cite{LR} for the convenience of readers.

Let $\lambda\in\Lambda$ and $M\in\Mod_{\Lambda,\lambda}(D)$.
Let $\mu\in\Lambda^+$.
Regarding $M$ as an object of $\Mod_{\Lambda}(\tU)$ via \eqref{eq:tUD} we have a monomorphism
\[
i_\mu:\omega^*M\to\omega^*(E^\mu\otimes_AM[-w_0\mu])
\]
and an epimorphism
\[
p_\mu:\omega^*(E^\mu\otimes_AM)\to\omega^*M[\mu]
\]
in $\CM(\tU)$ (see Section \ref{subsec:filt}).
Taking $\Gamma$ we obtain 
a monomorphism
\[
\Gamma(i_\mu):\Gamma(\omega^*M)\to\Gamma(\omega^*(E^\mu\otimes_AM[-w_0\mu]))
\]
and a morphism
\[
\Gamma(p_\mu):\Gamma(\omega^*(E^\mu\otimes_AM))\to\Gamma(\omega^*M[\mu])
\]
of $U$-modules (note that $\tU(0)=U$).
\begin{lemma}
\label{lem:KeyLemma}
\begin{itemize}
\item[\rm(i)]
Assume $\lambda+\rho\in\Lambda^+$.
Then there exists a splitting of $\Gamma(i_\mu)$ as $U$-modules, which is functorial with respect to $M$.
\item[\rm(ii)]
Assume $\lambda\in\Lambda^+$.
Then $\Gamma(p_\mu)$ is surjective.
\end{itemize}
\end{lemma}
\begin{proof}
By Lemma \ref{lem:cd1} we may assume that the canonical morphism $M\to\omega_*\omega^*M$ is an isomorphism.
By Section \ref{subsec:filt} we have a filtration 
\begin{equation}
\label{eq:key-filt}
\{0\}=\tilde{N}_0\subset
\tilde{N}_1\subset\cdots\subset
\tilde{N}_n\subset
\omega^*(E^\mu\otimes_AM)
\end{equation}
such that $\tilde{N}_k/\tilde{N}_{k-1}\simeq\omega^*M[\nu_k]^{\oplus m_k}$.
Here $\nu_k$ and $m_k$ are as in Section \ref{subsec:filt}.
We have
\[
\omega_*\omega^*(E^\mu\otimes_AM)
\simeq
E^\mu\otimes_A\omega_*\omega^*M
\simeq
E^\mu\otimes_AM
\]
by Proposition \ref{prop:filt-E1}, and hence by taking $\omega_*$ in 
\eqref{eq:key-filt} we obtain a filtration 
\[
\{0\}={N}_0\subset
{N}_1\subset\cdots\subset
{N}_n\subset
E^\mu\otimes_AM
\]
of $E^\mu\otimes_AM\in\Mod_\Lambda(\tU)$
with the exact sequence
\[
0\to N_{k-1}\to N_k\to M[\nu_k]^{\oplus m_k}
\to R^1\omega_*\tilde{N}_{k-1}.
\]
In particular, the quotient $N_k/N_{k-1}$ is isomorphic to a subobject of $M[\nu_k]^{\oplus m_k}$.
Hence we have
\[
z|(N_k/N_{k-1})(\xi)=\zeta_{\lambda+\nu_k+\xi}(z)\id
\qquad(z\in\Gz,\,\,\xi\in\Lambda)
\]
by \eqref{equ:deruz}.

In general for a $U$-module $V$ and $\xi\in\Lambda$ we set
\[
V{[\Gz,\zeta_\xi]}=
\{
v\in V\mid
\forall z\in\Gz\,\, \exists m\geqq0 \text{ such that } (z-\zeta_\xi(z))^mv=0\}.
\]
We have $(N_k/N_{k-1})(\xi)=((N_k/N_{k-1})(\xi)){[\Gz,\zeta_{\lambda+\nu_k+\xi}]}$.

(i) Set $\mu'=w_0\mu$.
Note
\begin{align*}
&\Gamma(\omega^*(E^\mu\otimes_AM[-\mu']))
=E^\mu\otimes_AM(-\mu')
=N_n(-\mu'),\\
&\Gamma(\omega^*M)
=\omega_*\omega^*(\overline{E}^\mu_1\otimes_AM)(-\mu')=N_1(-\mu').
\end{align*}
We have the canonical filtration
\[
\{0\}=N_0(-\mu')\subset N_1(-\mu')\subset\cdots\subset
N_n(-\mu')=\Gamma(\omega^*(E^\mu\otimes_AM[-\mu']))
\]
of $\Gamma(\omega^*(E^\mu\otimes_AM[-\mu']))$ consisting of $U$-submodules satisfying
\[
N_k(-\mu')/N_{k-1}(-\mu')=
(N_k(-\mu')/N_{k-1}(-\mu'))
[\Gz,\zeta_{\lambda+\nu_k-\mu'}].
\]
Let us show
\begin{equation}
\label{eq:key2}
\zeta_{\lambda+\nu_k-\mu'}=\zeta_{\lambda}\quad
\Longleftrightarrow\quad
k=1.
\end{equation}
By $\nu_1=\mu'$ we have $\zeta_{\lambda+\nu_1-\mu'}=\zeta_{\lambda}$.
Assume $\zeta_{\lambda+\nu_k-\mu'}=\zeta_{\lambda}$.
By Proposition \ref{prop:center2} (i) there exists some $w\in W$ satisfying 
$w(\lambda+\rho)=\lambda+\rho+\nu_k-\mu'$.
By $\lambda+\rho\in\Lambda^+$ we have $\lambda+\rho-w(\lambda+\rho)\in Q^+$.
Since $\mu'$ is the lowest weight, we have $\nu_k-\mu'\in Q^+$.
Hence we obtain $\nu_k-\mu'=0$, which implies $k=1$.
\eqref{eq:key2} is proved.
From this we obtain the canonical direct sum decomposition
\[
\Gamma(\omega^*(E^\mu\otimes_AM[-\mu']))
=\Gamma(\omega^*M)\oplus 
\sum_{\xi\in\Lambda, \xi\notin W\circ\lambda}
(\Gamma(\omega^*(E^\mu\otimes_AM[-\mu']))[\Gz,\zeta_\xi],
\]
which gives the desired splitting.

(ii) Note
\begin{align*}
&\Gamma(\omega^*(E^\mu\otimes_AM))
=E^\mu\otimes_AM
=N_n(0),\\
&\Gamma(\omega^*M[\mu])
=(\omega_*\omega^*M[\mu])(0)=(M[\mu])(0)=M(\mu).
\end{align*}
Consider the exact sequence
\[
0\to N_{n-1}(0)\to N_n(0)\to M(\mu)\to R^1\omega_*\tilde{N}_{n-1}(0).
\]
We have 
\begin{equation}
M(\mu)=(M(\mu))[\Gz,\zeta_{\lambda+\mu}].
\end{equation}
By the exact sequence
\[
R^1\omega_*\tilde{N}_{k-1}(0)\to
R^1\omega_*\tilde{N}_{k}(0)\to
R^1\omega_*(\tilde{N}_{k}/\tilde{N}_{k})(0)
\]
and $R^1\omega_*(\tilde{N}_{k}/\tilde{N}_{k-1})(0)\simeq
R^1\omega_*(\omega^*M)(\nu_k)^{\oplus m_k}$ we obtain
\[
R^1\omega_*\tilde{N}_{k}(0)=
\sum_{r=1}^{k}
(R^1\omega_*\tilde{N}_{k}(0))[\Gz,\zeta_{\lambda+\nu_r}].
\]
by Lemma \ref{lem:cd1}.
In particular, we have
\begin{equation}
R^1\omega_*\tilde{N}_{n-1}(0)=
\sum_{k=1}^{n-1}
(R^1\omega_*\tilde{N}_{n-1}(0))[\Gz,\zeta_{\lambda+\nu_k}].
\end{equation}
Hence it is sufficient to show 
\begin{equation}
\label{eq:key3}
\zeta_{\lambda+\nu_k}=\zeta_{\lambda+\mu}\quad
\Longleftrightarrow\quad
k=n.
\end{equation}
We have $\nu_n=\mu$ and hence $\zeta_{\lambda+\nu_k}=\zeta_{\lambda+\mu}$.
Assume that $\zeta_{\lambda+\nu_k}=\zeta_{\lambda+\mu}$.
By Proposition \ref{prop:center2} (i) there exists some $w\in W$ satisfying 
$w(\lambda+\rho+\nu_k)=\lambda+\rho+\mu$.
By $\lambda+\rho\in\Lambda^+$ we have $\lambda+\rho-w(\lambda+\rho)\in Q^+$.
Since $\mu$ is the highest weight, we have $\mu-w\nu_k\in Q^+$.
Hence we obtain $\lambda+\rho-w(\lambda+\rho)=0$.
By $\lambda\in\Lambda^+$ we have $w=1$ and hence $mu=\nu_k$, which implies $k=n$.
\eqref{eq:key3} is proved.
\end{proof}
\subsection{Proof of Theorem \ref{thm:A}}
\label{subsec:pf-of-A}
We first show the following.
\begin{lemma}
\label{lem:finite}
Let $M\in\Mod_\Lambda(A)$ and let $m\in\Gamma(\omega^*M)$.
Then there exists a finitely generated graded $A$-submodule $N$ of $M$ such that $m\in\Gamma(\omega^*N)$.
\end{lemma}
\begin{proof}
Let $f:M\to\omega_*\omega^*M$ be the canonical morphism.
Note that we have $m\in\Gamma(\omega^*M)=\omega_*\omega^*M(0)\subset\omega_*\omega^*M$.
By $\Coker(f)\in\Tor_\Lambda(A)$ there exists $\lambda\in\Lambda^+$ such that $A(\lambda)m\subset\Image(f)$.
Take a finite-dimensional subspace $V$ of $M(\lambda)$ such that 
$f(V)=A(\lambda)m$.
Set $N=\sum_{v\in V}Av\subset M$, 
$N'=Am\subset \omega_*\omega^*M$, and let 
$h:N\to N'$ be the morphism induced by $f$.
Note that $\omega^*N$ and $\omega^*N'$ are regarded as subobjects of $\omega^*M$ and $\omega^*(\omega_*\omega^*M)$ respectively.
Since $\omega^*f$ is an isomorphism, $\omega^*h$ is a monomorphism.
Let $\overline{m}$ be the element of $\Coker(h)$ corresponding to $m$.
By definition we have $A(\lambda)\overline{m}=\{0\}$, and hence $\overline{m}$ belongs to $\Tor(\Coker(h))$ by Lemma \ref{lem:grading}.
Since $\Coker(h)$ is generated by $\overline{m}$, we have $\Coker(h)\in\Tor_\Lambda(A)$, and hence $\omega^*h$ is an isomorphism.
It follows that $m\in\Gamma(\omega^*N')=\Gamma(\omega^*N)$.
\end{proof}

Let us give a proof of Theorem \ref{thm:A} (i)

Let $\lambda\in\Lambda$ such that $\lambda+\rho\in\Lambda^+$, and let $\overline{f}:\overline{M}\to\overline{N}$ be an epimorphism in $\CM_\lambda(D)$.
We need to show that 
$\Gamma(\overline{f}):\Gamma(\overline{M})\to\Gamma(\overline{N})$
is an epimorphism.

We first show that there exist $M, N\in\Mod_{\Lambda,\lambda}(D)$ and an epimorphism $f:M\to N$ such that 
\begin{equation}
\label{eq:p1}
\Tor(N)=\{0\},\qquad
\overline{M}=\omega^*M, \qquad
\overline{N}=\omega^*N,\qquad
\overline{f}=\omega^*f.
\end{equation}
Take $M, N\in\Mod_{\Lambda,\lambda}(D)$ such that $\overline{M}=\omega^*M, \overline{N}=\omega^*N$.
We may assume $\Tor(N)=\{0\}$.
Then  there exist $L\in\Mod_{\Lambda,\lambda}(D)$ and
morphisms $s:L\to M, \,\, f:L\to N$ such that
$s\in\Sigma_{D},\,\, \overline{f}=\omega^*f\circ(\omega^*s)^{-1}$.
By replacing $M$ with $L$ we may assume that $\overline{f}=\omega^*f$ for some $f:M\to N$.
Since $\overline{f}$ is an epimorphism, we have $\Tor(\Coker(f))=\Coker(f)$.
Hence by replacing $N$ with $\Image(f)$ we may further assume that $f$ is an epimorphism.
The existence of $M, N, f$ as in \eqref{eq:p1} is proved.

Let $n\in\Gamma(\overline{N})$.
We need to show that $n$ is contained in the image of $\Gamma(\overline{M})\to\Gamma(\overline{N})$.
By Lemma \ref{lem:finite} there exists a finitely generated graded $A$-submodule $N_1$ of 
$N$ such that $n\in\Gamma(\omega^*N_1)$.
Take a finitely generated graded $A$-submodule $M_1$ of $M$ such that $N_1=f(M_1)$.
By Proposition \ref{prop:Serre} there exists some $\nu\in \Lambda^+$ such that
\[
\Gamma(\omega^*M_1[\nu])\to\Gamma(\omega^*N_1[\nu])
\]
is surjective.
Set $\mu=-w_0\nu\in\Lambda^+$.
By Proposition \ref{prop:filt-E1} the vertical arrows of the commutative diagram 
\[
\begin{CD}
E^\mu\otimes_A\Gamma(\omega^*(M_1[\nu]))
@>>>
E^\mu\otimes_A\Gamma(\omega^*(N_1[\nu]))\\
@VVV @VVV\\
\Gamma(\omega^*(E^\mu\otimes_AM_1[\nu]))
@>>>
\Gamma(\omega^*(E^\mu\otimes_AN_1[\nu]))
\end{CD}
\]
are isomorphism, and hence 
\[
\Gamma(\omega^*(E^\mu\otimes_AM_1[\nu]))
\to
\Gamma(\omega^*(E^\mu\otimes_AN_1[\nu]))
\]
is surjective.
By \eqref{eq:i-mu} we can regard 
\[
\begin{CD}
\Gamma(\omega^*M_1)
@>>>
\Gamma(\omega^*N_1)\\
@VVV @VVV\\
\Gamma(\omega^*M)
@>>{k}>
\Gamma(\omega^*N)
\end{CD}
\]
as a subdiagram of the commutative diagram
\[
\begin{CD}
\Gamma(\omega^*(E^\mu\otimes_AM_1[\nu]))
@>{h}>>
\Gamma(\omega^*(E^\mu\otimes_AN_1[\nu]))\\
@VVV @VVV\\
\Gamma(\omega^*(E^\mu\otimes_AM[\nu]))
@>>{\ell}>
\Gamma(\omega^*(E^\mu\otimes_AN[\nu])).
\end{CD}
\]
Consider the commutative diagram
\[
\begin{CD}
\Gamma(\omega^*M)
@>{k}>>
\Gamma(\omega^*N)\\
@V{i_M}VV @VV{i_N}V\\
\Gamma(\omega^*(E^\mu\otimes_AM[\nu]))
@>>{\ell}>
\Gamma(\omega^*(E^\mu\otimes_AN[\nu])).
\end{CD}
\]
Since $n\in\Gamma(\omega^*N_1)$, the surjectivity of $h$ implies $i_N(n)\in\Image(\ell)$.
Then we obtain $n\in\Image(k)$ by the existence of the canonical splitting of $i_M$ and $i_N$ assured in Lemma \ref{lem:KeyLemma} (i).
The proof of Theorem \ref{thm:A} (i) is complete.

We next show Theorem \ref{thm:A} (ii).

Assume $\lambda\in\Lambda^+$, and let $\overline{M}$ be a non-zero object of $\CM_\lambda(D)$.
Set $M=\omega_*\overline{M}$.
By $\omega^*M=\overline{M}\ne0$ we have ${M}\notin\Tor_\Lambda(A)$, 
and hence $\Gamma(\omega^*M[\mu])=\omega_*\overline{M}(\mu)=M(\mu)\ne0$ for sufficiently large $\mu\in\Lambda^+$.
Hence by Proposition \ref{prop:filt-E1} and Lemma \ref{lem:KeyLemma} (ii) we have
\[
V(\mu)\otimes\Gamma(\omega^*M)
\simeq
E^\mu\otimes_A\Gamma(\omega^*M)
\simeq
\Gamma(\omega^*(E^\mu\otimes_AM))\ne0.
\]
This implies $\Gamma(\overline{M})=\Gamma(\omega^*M)\ne0$ as desired.
\subsection{Verma modules}
\begin{lemma}
\label{lem:G1}
For any $w\in W$ the canonical homomorphism
\begin{equation}
\label{eq:DShom2}
{S}_w^{-1}D\to\End_\BF(S_w^{-1}A)
\end{equation}
is injective.
\end{lemma}
\begin{proof}
Let $s\in{S}_w$ and $d\in D$, and assume that $s^{-1}d$ belongs to the kernel of \eqref{eq:DShom2}.
Then $d$ belongs to the kernel of \eqref{eq:DShom2} and hence it is contained in $\Ker(D\to\End_\BF(A))$.
Since $D$ is defined as a subalgebra of $\End_\BF(A)$, we have $d=0$ and hence $s^{-1}d=0$.
\end{proof}

For $\lambda\in\Lambda$ we have
\[
S_w^{-1}D_\lambda
=S_w^{-1}A\otimes_A(D\otimes_{\BF[\Lambda]}\BF)
=(S_w^{-1}A\otimes_AD)\otimes_{\BF[\Lambda]}\BF
=S_w^{-1}D\otimes_{\BF[\Lambda]}\BF,
\]
where $\BF[\Lambda]\to S_w^{-1}D$ and $\BF[\Lambda]\to \BF$ are given by $\e(\mu)\mapsto\sigma_\mu$ and $\e(\mu)\mapsto q^{(\mu,\lambda)}$ respectively.
Hence we have
\[
(S_w^{-1}D_\lambda)(0)
=(S_w^{-1}D)(0)\otimes_{\BF[\Lambda]}\BF.
\]
Since the image of $\BF[\Lambda]\to(S_w^{-1}D)(0)$ is contained in the center of $(S_w^{-1}D)(0)$, 
we have an $\BF$-algebra structure on $(S_w^{-1}D_\lambda)(0)$.
Moreover, since the action of $\sigma_\mu$ on $(S_w^{-1}A)(\lambda)$ is given by $q^{(\mu,\lambda)}\id$, we have a natural left $(S_w^{-1}D_\lambda)(0)$-module structure on $(S_w^{-1}A)(\lambda)$.
Hence we obtain an algebra homomorphism
\begin{equation}
\label{eq:Shom}
(S_w^{-1}D_\lambda)(0)\to\End_\BF((S_w^{-1}A)(\lambda))
\end{equation}

In the rest of this subsection we shall consider the case $w=1$.
\begin{lemma}
\label{lem:rs}
Let $s\in S_1$.
The right multiplication $r_s\in D$ is invertible in $S_1^{-1}D$ and the action of $r_s^{-1}\in S_1^{-1}D$ on $S_1^{-1}A$ is given by the right multiplication of $s^{-1}\in S_1^{-1}A$.
\end{lemma}
\begin{proof}
We have $s\in A(\mu)_\mu\setminus\{0\}$ for some $\mu\in\Lambda^+$.
Then $r_s=\ell_s\deru_{k_\mu}\sigma_{-\mu}$ by Lemma \ref{lem:rl}.
Hence $r_s^{-1}=\sigma_\mu\deru_{k_{-\mu}}\ell_s^{-1}\in S_1^{-1}D$.
For $t\in S_1$ and $\varphi\in A$ we have
\[
r_s(t^{-1}\varphi)
=r_s\ell_t^{-1}(\varphi)
=\ell_t^{-1}r_s(\varphi)
=t^{-1}\varphi s
\]
by $r_s\ell_t=\ell_t r_s$.
Thus the action of $r_s$ on $S_1^{-1}A$ is given by the right multiplication of $s$.
Hence its inverse is given by the right multiplication by $s^{-1}$.
\end{proof}
By definition the image $\overline{D}$ of the canonical injective algebra homomorphism
\begin{equation}
\label{eq:DShom3}
({S}_1^{-1}D)(0)\to\End_\BF(S_1^{-1}A)
\end{equation}
is generated by 
\begin{itemize}
\item[\rm(a)]
$\deru_u$ for $u\in U$, 
\item[\rm(b)]
$\sigma_\mu$ for $\mu\in\Lambda$,
\item[\rm(c)]
$L_\varphi:S_1^{-1}A\to S_1^{-1}A\,\,(\psi\mapsto\varphi\psi)$, where
$\varphi\in(S_1^{-1}A)(0)$.
\end{itemize}
By Lemma \ref{lem:rs} $\overline{D}$ also contains
\begin{itemize}
\item[\rm(d)]
$R_\varphi:S_1^{-1}A\to S_1^{-1}A\,\,(\psi\mapsto \psi\varphi)$, where
$\varphi\in(S_1^{-1}A)(0)$.
\end{itemize}
Since $\overline{D}$ preserves $(S_1^{-1}A)(\xi)$ for each $\xi\in\Lambda$, 
we have
\[
\overline{D}\subset\prod_{\xi\in\Lambda}\End_\BF(F^{\geqq0}(\xi)).
\]
Here, we have identified $(S_1^{-1}A)(\xi)$ with $F^{\geqq0}(\xi)$ 
by Proposition \ref{prop:theta}.
In particular, $F^{\geqq0}(\xi)$ is regarded as a $U$-module.
Note that $F^{\geqq0}(\xi)$ is isomorphic to $T^*(\xi)$ as a $U$-module by Proposition \ref{prop:S-1A}.

Set
\begin{align*}
&\overline{\End}_\BF(F^{\geqq0}(\xi))\\
=&\bigoplus_{\nu\in\Lambda}
\{
f\in\End_\BF(F^{\geqq0}(\xi))\mid
f(F^{\geqq0}(\xi)_\mu)\subset F^{\geqq0}(\xi)_{\mu+\nu}\,\,(\forall \mu\in\Lambda)\},
\\
&\overline{\End}_\BF(F^{\geqq0}(\xi)^\bigstar)\\
=&\bigoplus_{\nu\in\Lambda}
\{
f\in\End_\BF(F^{\geqq0}(\xi)^\bigstar)\mid
f(F^{\geqq0}(\xi)^\bigstar_\mu)\subset F^{\geqq0}(\xi)^\bigstar_{\mu+\nu}\,\,(\forall \mu\in\Lambda)\}.
\end{align*}
Then we have an isomorphism of $\BF$-algebras 
\begin{gather*}
\overline{\End}_\BF(F^{\geqq0}(\xi))\simeq
\overline{\End}_\BF(F^{\geqq0}(\xi)^\bigstar)^{\rm op}\qquad
(h\leftrightarrow {}^th)\\
\langle h(v^*),v\rangle=\langle v^*,{}^t h(v)\rangle\qquad
(v^*\in F^{\geqq0}(\xi), v\in F^{\geqq0}(\xi)^\bigstar).
\end{gather*}
Hence we obtain an embedding of $\BF$-algebras
\begin{gather*}
\Theta:\overline{D}\hookrightarrow
\prod_{\xi\in\Lambda}
\End_\BF(F^{\geqq0}(\xi)^\bigstar)^{\rm op},\qquad
\Theta(d)=(\Theta_\xi(d))_{\xi\in\Lambda},\\
\langle v^*,(\Theta_\xi(d))(v)\rangle=\langle d(v^*),v\rangle
\qquad(v^*\in F^{\geqq0}(\xi), v\in F^{\geqq0}(\xi)^\bigstar).
\end{gather*}
Since $F^{\geqq0}(\xi)^\bigstar$ is isomorphic to $T_{\rm r}(\xi)$ as a right $U$-module, it is a free right $U^+$-module generated by the element $n_\xi\in F^{\geqq0}(\xi)^\bigstar$ given by
\[
\langle\varphi,n_\xi\rangle=\langle\varphi, 1\rangle
\qquad
(\varphi\in F^{\geqq0}(\xi)).
\]
Here, $\langle\,\,,\,\,\rangle$ in the left hand side (resp.\ the right hand side) is the canonical paring 
$F^{\geqq0}(\xi)\times F^{\geqq0}(\xi)^\bigstar\to\BF$ 
(resp.\ $F^{\geqq0}\times U^{\geqq0}\to\BF$).
We shall identify $F^{\geqq0}(\xi)^\bigstar$ with $U^+$ by 
$(F^{\geqq0}(\xi)^\bigstar\ni n_\xi u\leftrightarrow u\in U^+)$.
Hence we have an embedding of $\BF$-algebras
\begin{gather*}
\Theta:\overline{D}\hookrightarrow
\prod_{\xi\in\Lambda}
\End_\BF(U^+)^{\rm op},\qquad
\Theta(d)=(\Theta_\xi(d))_{\xi\in\Lambda},\\
\langle \varphi,(\Theta_\xi(d))(u)\rangle=\langle d(\varphi),n_\xi u\rangle
=\langle d(\varphi),u\rangle
\qquad(\varphi\in F^{\geqq0}(\xi), u\in U^+).
\end{gather*}

For $x\in U^+, \mu\in\Lambda$ we define $M_x, N_\mu\in\End_\BF(U^+)$ by
\[
M_x(u)=ux,\qquad
N_\mu(u)=k_\mu uk_\mu^{-1}
\qquad(u\in U^+).
\]
We define linear maps
\[
F^{\geqq0}(0)\ni\varphi\mapsto P_\varphi\in\End_\BF(U^+),\qquad
F^{\geqq0}(0)\ni\varphi\mapsto Q_\varphi\in\End_\BF(U^+)
\]
by
\begin{align*}
&P_\varphi(u)=\sum_{(u)_1}\langle\varphi,u_{(0)}\rangle u_{(1)}
\qquad
&(\varphi\in F^{\geqq0}(0),\,\, u\in U^+),\\
&Q_\varphi(u)=\sum_{(u)_1}\langle\varphi,u_{(1)}\rangle k_\gamma^{-1}u_{(0)}
\qquad
&(\gamma\in Q^+,\,\, \varphi\in F^{\geqq0}(0)_{-\gamma},\,\, u\in U^+).
\end{align*}
Note that $P_\varphi(u), Q_\varphi(u)\in U^+$ by 
\begin{align}
&\Delta(u)\in\sum_{\gamma\in Q^+}k_\gamma U^+\otimes U^+_\gamma
\qquad(u\in U^+),\\ 
&\gamma\ne\delta
\quad\Rightarrow\quad
\langle F^{\geqq0}(0)_{-\gamma},U^+_\delta\rangle=\{0\}.
\end{align}
For $i\in I$ we define $\varphi_i\in F^{\geqq0}(0)$ by 
\[
\langle\varphi_i,tu\rangle=(u,f_i)\qquad(t\in U^0, \,\,u\in U^+),
\]
where $(\,\,,\,\,)$ in the right hand side is the paring \eqref{eq:Drinfeld-paring}.
For $F\in\End_\BF(U^+)$ we define $\Delta(F)\in\prod_{\xi\in\Lambda}
\End_\BF(U^+)$ as the corresponding diagonal element.
\begin{lemma}
\label{lem:action}
We have
\begin{align}
&\Theta_\xi(\sigma_\mu)=q^{(\mu,\xi)}\id
\qquad\qquad(\mu\in\Lambda),
\label{eq:action1}\\
&\Theta(L_\varphi)=\Delta(P_\varphi)
\qquad\qquad(\varphi\in F^{\geqq0}(0)),
\label{eq:action4}\\
&\Theta(R_\varphi)=\Delta(Q_\varphi)\Theta(\sigma_{\gamma})
\qquad\qquad(\gamma\in Q^+,\,\,\varphi\in F^{\geqq0}(0)_{-\gamma}),
\label{eq:action5}\\
&\Theta(\deru_x)=\Delta(M_x)
\qquad\qquad(x\in U^+),
\label{eq:action2}\\
&\Theta(\deru_{k_\mu})=\Delta(N_{-\mu})\Theta(\sigma_{\mu})
\qquad\qquad(\mu\in\Lambda),
\label{eq:action3}\\
&\Theta(\deru_{f_i})=
q^{-(\alpha_i,\alpha_i)}
\Delta(P_{\varphi_i}N_{\alpha_i})\Theta(\sigma_{-\alpha_i})
-\Delta(Q_{\varphi_i})\Theta(\sigma_{\alpha_i})
\qquad(i\in I).
\label{eq:action6}
\end{align}
\end{lemma}
\begin{proof}
\eqref{eq:action1} is easy.
Let $\psi\in F^{\geqq0}(\xi),\,\,u\in U^+$.
For $\varphi\in F^{\geqq0}(0)$ we have
\[
\langle\psi,(\Theta_\xi(L_\varphi))(u)\rangle
=\langle\varphi\psi,u\rangle
=\sum_{(u)_1}\langle\varphi,u_{(0)}\rangle\langle\psi,u_{(1)}\rangle
=\langle\psi,P_\varphi(u)\rangle
\]
and hence we obtain \eqref{eq:action4}.
Similarly, for $\varphi\in F^{\geqq0}(0)_{-\gamma}$ we have
\begin{align*}
\langle\psi,(\Theta_\xi(R_\varphi))(u)\rangle
&=\sum_{(u)_1}\langle\varphi,u_{(1)}\rangle\langle\psi,u_{(0)}\rangle
=\langle\psi,k_\gamma Q_\varphi(u)\rangle
=\langle\psi k_\gamma,Q_\varphi(u)\rangle\\
&=q^{(\gamma,\xi)}\langle\psi,Q_\varphi(u)\rangle.
\end{align*}
\eqref{eq:action5} is proved.
For $x\in U$ we have 
\[
\langle\psi,(\Theta_\xi(\deru_x))(u)\rangle
=\langle x\psi,n_\xi u\rangle
=\langle \psi,n_\xi ux\rangle.
\]
In case $x\in U^+$ we have $ux\in U^+$, and we obtain \eqref{eq:action2}.
In case $x=k_\mu$ for $\mu\in\Lambda$ we have
\[
n_\xi uk_\mu=n_\xi k_\mu(k_\mu^{-1}uk_\mu)
=q^{(\xi,\mu)}n_\xi N_{-\mu}(u),
\]
and hence \eqref{eq:action3} is proved.
In case $x=f_i$ for $i\in I$ we have
\begin{align*}
&uf_i\\
=&
\sum_{(u)_2}(u_{(0)},f_i)(u_{(2)},k_i)k_i^{-1}u_{(1)}
+\sum_{(u)_2}(u_{(0)},1)(u_{(2)},k_i)f_iu_{(1)}\\
&\qquad\qquad\qquad
-\sum_{(u)_2}(u_{(0)},1)(u_{(2)},f_ik_i)u_{(1)}\\
=&
\sum_{(u)_1}(u_{(0)},f_i)k_i^{-1}u_{(1)}
+\sum_{(u)_1}(u_{(1)},k_i)f_iu_{(0)}
-\sum_{(u)_1}(u_{(1)},f_ik_i)u_{(0)}
\end{align*}
by \cite[Lemma 2.1.2, Lemma 2.1.3]{T}, and hence if $u\in U^+_\delta$, then 
\begin{align*}
n_\xi uf_i=&
n_\xi(\sum_{(u)_1}(u_{(0)},f_i)k_i^{-1}u_{(1)})
-n_\xi(\sum_{(u)_1}(u_{(1)},f_ik_i)u_{(0)})\\
=&q^{(\delta-\alpha_i-\xi,\alpha_i)}
n_\xi(\sum_{(u)_1}\langle\varphi_i,u_{(0)}\rangle u_{(1)})
-
n_\xi(\sum_{(u)_1}\langle\varphi_i,u_{(1)}\rangle u_{(0)})\\
=&q^{(\delta-\alpha_i-\xi,\alpha_i)}
n_\xi(P_{\varphi_i}(u))
-q^{(\xi,\alpha_i)}
n_\xi(Q_{\varphi_i}(u))
\end{align*}
and we obtain \eqref{eq:action6}.
\end{proof}
By \eqref{eq:Shom} we obtain an algebra homomorphism
\begin{equation}
\label{eq:Shom2}
(S_1^{-1}D_\lambda)(0)\to\End_\BF(T_{\rm r}(\lambda))^{\rm op}.
\end{equation}
\begin{proposition}
\label{prop:injectivity}
The algebra homomorphism \eqref{eq:Shom2} is injective.
\end{proposition}
\begin{proof}
Note that 
$(S_1^{-1}D_\lambda)(0)=\overline{D}\otimes_{\BF[\Lambda]}\BF$, where $\BF[\Lambda]\to\overline{D}$ and $\BF[\Lambda]\to\BF$ are given by 
$\e(\mu)\mapsto\sigma_\mu$ and $\e(\mu)\mapsto q^{(\mu,\lambda)}$ respectively.
The algebra homomorphism \eqref{eq:Shom2} is induced from $\Theta_\lambda:\overline{D}\to\End_\BF(U^+)^{\rm op}$ under the identification $T_{\rm r}(\lambda)=U^+$.

Let $\overline{D}'$ be the subalgebra of $\prod_{\xi\in\Lambda}
\End_\BF(U^+)^{\rm op}$ generated by the elements $\Delta(M_x)$, $\Delta(N_\mu)$, $\Delta(P_\varphi)$, $\Delta(Q_\varphi)$ for $x\in U^+$, $\mu\in\Lambda$, $\varphi\in F^{\geqq0}(0)$.
By Lemma \ref{lem:action} we have $\overline{D}'\subset\overline{D}$, and $\overline{D}$ is generated as a subalgebra of  $\prod_{\xi\in\Lambda}\End_\BF(U^+)^{\rm op}$ by $\overline{D}'$ and $\{\sigma_\mu\mid\mu\in\Lambda\}$.
Let us show that the linear map 
\[
\overline{D}'\otimes\BF[\Lambda]
\ni d\otimes\e(\mu)\mapsto d\sigma_\mu\in
\overline{D}
\]
is an isomorphism.
The surjectivity is a consequence of the fact that $\sigma_\mu$ belongs to the center of $\overline{D}$.
Assume that we have $\sum_\mu d_\mu\sigma_\mu=0$ for 
$d_\mu\in \overline{D}'$.
For any $\xi\in\Lambda$ and $u\in U^+$ we have 
$\sum_\mu q^{(\mu,\xi)}d_\mu(u)=0$, from which we obtain $d_\mu=0$ for any $\mu$.
Hence $\overline{D}'\otimes\BF[\Lambda]\to\overline{D}$ is bijective.
It follows that
\[
(S_1^{-1}D_\lambda)(0)
=\overline{D}\otimes_{\BF[\Lambda]}\BF
\simeq \overline{D}',
\]
and hence 
$(S_1^{-1}D_\lambda)(0)\to\End_\BF(T_{\rm r}(\lambda))^{\rm op}$
is injective.
\end{proof}

\subsection{Proof of Theorem \ref{thm:B}}
\label{subsec:{thm:B}}
For $\lambda\in\Lambda$ consider the sequence
\begin{equation}
\label{eq:abcd}
U/J_\lambda
\xrightarrow{\alpha}
D_\lambda(0)
\xrightarrow{\beta}
\Gamma(\omega^*D_\lambda)
\xrightarrow{\gamma}
(S_1^{-1}D_\lambda)(0)
\xrightarrow{\delta}
\End_\BF(T_{\rm r}(\lambda))^{\rm op}
\end{equation}
of algebra homomorphisms.
We shall show that $\alpha$ is an isomorphism for any $\lambda\in\Lambda$ and that $\beta$ is an isomorphism if $\lambda+\rho\in\Lambda^+$.
We need the following result due to Joseph \cite{J:an}.
\begin{proposition}
[Joseph]
\label{prop:Joseph}
The homomorphism 
\[
U/J_\lambda\to\End_\BF(T_{\rm r}(\lambda))^{\rm op}
\]
is injective.
\end{proposition}

We see easily that 
$\alpha$ is surjective by the definition of $D_\lambda$.
By Proposition \ref{prop:Joseph} the composition
$\delta\circ\gamma\circ\beta\circ\alpha$ is injective.
Hence $\alpha$ is an isomorphism and $\beta$ is injective.
It remains to show the surjectivity of $\beta$ in the case $\lambda+\rho\in\Lambda^+$.

Assume that we are given an $\BF$-algebra $E$ and an algebra homomorphism $f:U\to E$.
We have a left $U$-module structure on $E$
\[
ue=f(u)e\qquad(u\in U,\,\, e\in E).
\]
given by the left multiplication.
We have also another left $U$-module structure on $E$
\[
\ad(u)(e)=\sum_{(u)_1}f(u_{(0)})ef(Su_{(1)})\qquad
(u\in U,\,\, e\in E)
\]
called the adjoint action.
Taking the $U$-finite parts
\[
E^\fin=
\{e\in E\mid \dim_\BF \ad(U)(e)<\infty\}
\]
of $E=U$, $D_\lambda(0)$, $\Gamma(\omega^*D_\lambda)$, 
$(S_1^{-1}D_\lambda)(0)$, $\End_\BF(T_{\rm r}(\lambda))^{\rm op}$
we obtain 
\begin{equation}
\label{eq:abcd2}
U^\fin 
\xrightarrow{\overline{\alpha}}
D_\lambda(0)^\fin
\xrightarrow{\overline{\beta}}
\Gamma(\omega^*D_\lambda)^\fin
\xrightarrow{\overline{\gamma}}
(S_1^{-1}D_\lambda)(0)^\fin
\xrightarrow{\overline{\delta}}
(\End_\BF(T_{\rm r}(\lambda))^{\rm op})^\fin
\end{equation}

We need the following result of Joseph \cite[Theorem 8.3.9 (ii)]{J:book} which is a $q$-analogue of a theorem of N. Conze-Berline and M. Duflo.
\begin{proposition}
[Joseph]
\label{prop:CDJ}
The homomorphism 
\[
U^\fin 
\to
(\End_\BF(T_{\rm r}(\lambda))^{\rm op})^\fin
\]
is surjective.
\end{proposition}
As shown above $\beta$ is injective, and hence $\overline{\beta}$ is injective.
Since $\delta$ is injective by Proposition \ref{prop:injectivity}, $\overline{\delta}$ is also injective.
Moreover, we have the injectivity of $\overline{\gamma}$ by Lemma 
\ref{lem:fin-inj}.
Hence Proposition \ref{prop:CDJ} implies that $\overline{\beta}$, $\overline{\gamma}$, $\overline{\delta}$ are isomorphisms.

Now we deduce the surjectivity of $\beta$ from the surjectivity of $\overline{\beta}$.
The assumption $\lambda+\rho\in\Lambda^+$ will be used in the arguments below.

Let $m\in\Gamma(\omega^*D_\lambda)$.
We shall show that $m\in\Image({\beta})$.
By Lemma \ref{lem:finite} there exists a finitely generated graded $A$-submodule $M$ of $D_\lambda$ such that $m\in\Gamma(\omega^*M)$.

Let us show that the canonical morphism
\begin{equation}
\label{eq:Mnu}
M(\nu)\to\Gamma(\omega^*M[\nu])
\end{equation}
is surjective for sufficiently large $\nu\in\Lambda^+$.
Since $M$ is finitely generated, there exists an epimorphism
\[
\bigoplus_{j=1}^nA[\xi_j]\to M
\]
in $\Mod_\Lambda(A)$ for some $\xi_1,\dots,\xi_n\in\Lambda$.
Consider the commutative diagram
\[
\begin{CD}
\bigoplus_{j=1}^nA(\xi_j+\nu)@>>>M(\nu)
\\
@VVV@VVV
\\
\bigoplus_{j=1}^n\Gamma(\omega^*A[\xi_j+\nu])
@>>>
\Gamma(\omega^*M[\nu]).
\end{CD}
\]
Assume that $\nu\in\Lambda^+$ is sufficiently large.
Then the left vertical arrow is an isomorphism by Proposition \ref{prop:Borel-Weil}, and the lower horizontal arrow is surjective by Proposition \ref{prop:Serre}.
Hence \eqref{eq:Mnu} is surjective.

Take $\nu\in\Lambda^+$ such that \eqref{eq:Mnu} is surjective.
Set $\mu=-w_0\nu$, and consider the following commutative diagram (see \eqref{eq:i-mu}):
\[
\begin{CD}
\Gamma(\omega^*M)
@>>>
\Gamma(\omega^*(E^\mu\otimes_AM[\nu]))
\\
@VVV@VV{k}V
\\
\Gamma(\omega^*D_\lambda)
@>{i}>>
\Gamma(\omega^*(E^\mu\otimes_AD_\lambda[\nu])).
\end{CD}
\]
Since $\lambda+\rho\in\Lambda^+$, there exists a homomorphism 
\[
j:\Gamma(\omega^*(E^\mu\otimes_AD_\lambda[\nu]))
\to
\Gamma(\omega^*D_\lambda)
\]
of $U$-modules such that $j\circ i=\id$ by Lemma \ref{lem:KeyLemma}.
Here, the action of $U$ on $E^\mu\otimes_AD_\lambda$ is given by
\begin{equation}
\label{eq:action-left}
u\cdot(e\otimes\overline{d})
=\sum_{(u)_1}u_{(0)}e\otimes\overline{\deru_{u_{(1)}}d}
\qquad
(u\in U,\,\,d\in D),
\end{equation}
where $\overline{d}$ denotes the element of $D_\lambda$ corresponding to $d$, and the one on $\Gamma(\omega^*D_\lambda)$ is given by the left multiplication.
Hence it is sufficient to show $\Image(j\circ k)\subset\Image({\beta})$.
By Proposition \ref{prop:filt-E1} we have
\[
\omega_*\omega^*(E^\mu\otimes_AM[\nu])
\simeq
E^\mu\otimes_A\omega_*\omega^*(M[\nu])
\simeq
V(\mu)\otimes\omega_*\omega^*(M[\nu]),
\]
and hence 
$\Gamma(\omega^*(E^\mu\otimes_AM[\nu]))\simeq
V(\mu)\otimes\Gamma(\omega^*(M[\nu]))$.
Similarly we have 
$\Gamma(\omega^*(E^\mu\otimes_AD_\lambda[\nu]))\simeq
V(\mu)\otimes\Gamma(\omega^*(D_\lambda[\nu]))$.
Consider the commutative diagram
\[
\begin{CD}
V(\mu)\otimes M(\nu)
@>>>
V(\mu)\otimes D_\lambda(\nu)
\\
@V{\ell'}VV@VV{\ell}V
\\
V(\mu)\otimes \Gamma(\omega^*(M[\nu]))
@.
V(\mu)\otimes \Gamma(\omega^*(D_\lambda[\nu]))
\\
@|@|
\\
\Gamma(\omega^*(E^\mu\otimes_AM[\nu]))
@>>{k}>
\Gamma(\omega^*(E^\mu\otimes_AD_\lambda[\nu])).
\end{CD}
\]
Since \eqref{eq:Mnu} is surjective, 
$\ell'$ is also surjective.
Hence it is sufficient to show that the image of the composition of
\[
V(\mu)\otimes D_\lambda(\nu)
\xrightarrow{\ell}
\Gamma(\omega^*(E^\mu\otimes_AD_\lambda[\nu]))
\xrightarrow{j}
\Gamma(\omega^*D_\lambda)
\]
is contained in $\Image({\beta})$.
We regard $V(\mu)\otimes D_\lambda(\nu)$, 
$\Gamma(\omega^*(E^\mu\otimes_AD_\lambda[\nu]))$, 
$\Gamma(\omega^*D_\lambda)$ as right $D_\lambda(0)$-modules via the right multiplication of $D(0)$ on $D$.
Then $\ell$ and $\beta$ are homomorphisms of right $D_\lambda(0)$-modules by definition.
Moreover, $j$ is also a homomorphisms of right $D_\lambda(0)$-modules by the following reason.
Recall that $j$ is the projection with respect to the action of the center $\Gz$ of $U$.
Here, the action of $U$ on $E^\mu\otimes_AD_\lambda$ is given by
\eqref{eq:action-left}, and the one on $\Gamma(\omega^*D_\lambda)$ is given by the left multiplication.
Since the action of $U$ and the right action of $D_\lambda(0)$ on 
$\Gamma(\omega^*(E^\mu\otimes_AD_\lambda[\nu]))$ 
commute with each other, $j$ is a homomorphism of right $D_\lambda(0)$-modules.
Let 
\[
r:A(\nu)\to D_\lambda(\nu)
\]
be the composition of $A(\nu)\hookrightarrow D(\nu)\to D_\lambda(\nu)$.
Since $D(\nu)$ is generated by $A(\nu)$ as a right $D(0)$-module, $D_\lambda(\nu)$ is  generated by $r(A(\nu))$ as a right $D_\lambda(0)$-module.
Hence it is sufficient to show 
$(j\circ\ell)(V(\mu)\otimes r(A(\nu)))\subset\Image(\beta)$.

We can regard $V(\mu)\otimes D_\lambda(\nu)$ as a left $U$-module
\begin{equation}
\label{eq:action-left2}
u\cdot(v\otimes\overline{d})
=\sum_{(u)_1}u_{(0)}v\otimes\overline{\deru_{u_{(1)}}d}
\qquad
(u\in U,\,\,d\in D(\nu)),
\end{equation}
We  can also regard $V(\mu)\otimes D_\lambda(\nu)$ as a left $U$-module
by the adjoint action defined by
\[
\ad(u)(v\otimes \overline{d})
=\sum_{(u)_2}u_{(0)}v\otimes\overline{\deru_{u_{(1)}}d\deru_{Su_{(2)}}}
\qquad
(u\in U,\,\,d\in D(\nu)).
\]
Note that $j\circ\ell$ is a homomorphism of right $D_\lambda(0)$-modules as well as a homomorphism of left $U$-modules,
where the $U$-module structure on  $V(\mu)\otimes D_\lambda(\nu)$ is given by \eqref{eq:action-left2} and the one on $\Gamma(\omega^*D_\lambda)$ is given by the left multiplication.
Hence  $j\circ\ell$ also preserves the adjoint actions of $U$.
Since $V(\mu)\otimes r(A(\nu))$ is stable under the adjoint action of $U$, we have $V(\mu)\otimes r(A(\nu))\subset(V(\mu)\otimes D_\lambda(\nu))^\fin$ with respect to the adjoint action.
It follows that
\begin{align*}
(j\circ\ell)(V(\mu)\otimes r(A(\nu)))
&\subset
(j\circ\ell)((V(\mu)\otimes D_\lambda(\nu))^\fin)
\subset
\Gamma(\omega^*D_\lambda)^\fin\\
&\subset
\Image(\overline{\beta})
\subset
\Image({\beta}).
\end{align*}
The proof of Theorem \ref{thm:B} is complete.
\bibliographystyle{unsrt}

\begin{thebibliography}{99}
\bibitem{AZ}
M. Artin, J. J. Zhang,
\newblock
{\it Noncommutative projective schemes},
\newblock 
Adv. in Math. {\bf109} (1994), 228--287.
\bibitem{BB}
A. Beilinson, J. Bernstein,
\newblock
{\it Localisation de $\Gg$-modules},
\newblock 
C. R. Acad. Sci. Paris S\'er. I Math. {\bf 292} (1981), 15--18.
\bibitem{BoB}
W. Borho, J.-L. Brylinski,
\newblock
{\it Differential operators on homogeneous spaces. I},
\newblock 
Invent. Math. {\bf 69} (1982), 437--476.
\bibitem{BGG}
I. N. Bernstein, I. M. Gelfand, S. I. Gelfand,
\newblock
{Differential operators on the base affine space and a study of  $\Gg$-modules},
\newblock
in: {\it Lie groups and their representations}, 
\newblock
I. M. Gelfand (ed.), A. Hilger, London 1975, 22--64.
\bibitem{D}
V. G. Drinfeld, 
\newblock
{Quantum groups},
\newblock 
{\it Proc. Int. Cong. Math.},  Berkeley
1986, 798--820.
\bibitem{GZ}
P. Gabriel, M. Zisman,
\newblock
{\it Calculus of fractions and homotopy theory}, 
Springer-Verlag, 1967.
\bibitem{Hodges}
T. J. Hodges,
\newblock
{\it Ring-theoretical aspects of the Bernstein-Beilinson theorem},
\newblock
LNM {\bf 1448}, 1990, 155--163.
\bibitem{J:an}
A. Joseph, 
\newblock
{Enveloping algebras: Problems old and new},
\newblock
in: Progress in Math. {\bf123}, Birkh\"auser, 1994, 385--413.
\bibitem{J:emb}
A. Joseph, 
\newblock
{Faithfully flat embeddings for minimal primitive quotients of quantized enveloping algebras},
\newblock
in: {\it Quantum Deformations of Algebras and their representation theory}, 
\newblock
A. Joseph, S. Shnider (eds.),
Israel Math. Conf. Proc. 1993, 79--106.
\bibitem{J:book}
A. Joseph,
\newblock
{\it Quantum groups and their primitive ideals}, 
\newblock
Springer-Verlag, 1995.
\bibitem{LRD}
V. A. Lunts, A. L. Rosenberg, 
\newblock
{\it Differential operators on noncommutative rings}, 
\newblock
Selecta Math. New ser. {\bf 3} (1997), 335--359.
\bibitem{LR}
V. A. Lunts, A. L. Rosenberg, 
\newblock
{\it Localization for quantum groups}, 
\newblock
Selecta Math. New ser. {\bf 5} (1999), 123-159.
\bibitem{L1}
G. Lusztig, 
\newblock
{\it Quantum deformations of certain simple modules over enveloping algebras},
\newblock 
Adv. in Math. {\bf 70} (1988), 237--249.
\bibitem{Lbook}
G. Lusztig, 
\newblock
{\it Introduction to quantum groups},
\newblock 
Birkh\"auser, 1993.
\bibitem{Manin}
Yu. I. Manin, 
\newblock
{\it Topics in Noncommutative geometry},
\newblock
Princeton Univ Press, 1991.
\bibitem{P}
N. Popescu,
\newblock
{\it Abelian categories with applications to rings and modules},
\newblock
Academic Press, 1973.
\bibitem{R:book}
A. L. Rosenberg, 
\newblock
{\it Noncommutative algebraic geometry and representations of quantized algebras},
\newblock
Kluwer Academic Publishers, 1995.
\bibitem{R}
A. L. Rosenberg, 
\newblock
{\it Noncommutative schemes},
\newblock
Comp. Math. {\bf112} (1995), 93--125.
\bibitem{T}
T. Tanisaki, 
\newblock
{\it Killing forms, Harish-Chandra isomorphisms, and universal
$R$-matrices for quantum algebras},
\newblock 
Inter. J. Mod. Phys. {\bf A7}, Suppl. 1B 
(1992), 941--961.
\bibitem{V}
A. B. Verevkin,
\newblock
{\it On a noncommutative analogue of the category of coherent sheaves on a projective scheme},
\newblock 
Amer. Math. Soc. Transl. (2) {\bf 151}
(1992).
\end{thebibliography}

\end{document}